\documentclass{gtart}
\usepackage{amsmath}
\usepackage{amssymb}
\usepackage{texdraw}

\newcommand{\Z}{\mathbb{Z}}
\newcommand{\Q}{\mathbb{Q}}
\newcommand{\R}{\mathbb{R}}
\newcommand{\C}{\mathbb{C}}
\newcommand{\Io}{\mathbb I}

\newcommand{\ch}{h^{\vee}}
\newcommand{\cY}{Y^{\vee}}
\newcommand{\ep}{\epsilon}

\newcommand{\I}{\sqrt{-1}}
\newcommand{\la}{\langle}
\newcommand{\npr}{|\Delta_{+}|}
\newcommand{\ra}{\rangle}
\newcommand{\ria}{\rightarrow}
\newcommand{\vep}{\varepsilon}

\newcommand{\inte}{\operatorname{int}}
\newcommand{\sign}{\operatorname{sign}}

\newcommand{\vol}{\operatorname{vol}}

\newcommand{\aW}{W^{\text{\rm aff}}}
\newcommand{\col}{\text{\rm col}}

\newcommand{\Diag}{\text{\rm Diag}}
\newcommand{\dS}{\text{\rm S}}
\newcommand{\End}{\text{\rm End}}
\newcommand{\FWC}{\text{\rm FWC}}
\newcommand{\GL}{\text{\rm GL}}
\newcommand{\Hom}{\text{\rm Hom}}
\newcommand{\id}{\text{\rm id}}
\newcommand{\interior}{\text{\rm int}}
\newcommand{\LHS}{\text{\rm LHS}}
\newcommand{\ns}{\text{\rm n}}
\newcommand{\os}{\text{\rm o}}
\newcommand{\PSL}{\text{\rm PSL}}
\newcommand{\RHS}{\text{\rm RHS}}
\newcommand{\s}{\text{\rm s}}
\newcommand{\SL}{\text{\rm SL}}
\newcommand{\Span}{\text{\rm Span}}
\newcommand{\tO}{\text{\rm O}}

\newcommand{\tr}{\mbox{\rm tr}}

\newcommand{\frg}{\mathfrak g}
\newcommand{\frh}{\mathfrak h}
\newcommand{\frhR}{{\mathfrak h}_{\mathbb{R}}}
\newcommand{\frsl}{\mathfrak sl}

\newcommand{\mA}{\mathcal A}
\newcommand{\mC}{\mathcal C}
\newcommand{\mD}{\mathcal D}
\newcommand{\mM}{\mathcal M}
\newcommand{\mR}{\mathcal R}
\newcommand{\mS}{\mathcal S}
\newcommand{\mT}{\mathcal T}
\newcommand{\mV}{\mathcal V}
\newcommand{\mW}{\mathcal W}

\newcommand{\sm}{\setminus}

\newtheorem{thm}{Theorem}[section]
\newtheorem{conj}[thm]{Conjecture}
\newtheorem{cor}[thm]{Corollary}
\newtheorem{lem}[thm]{Lemma}
\newtheorem{prop}[thm]{Proposition}

\theoremstyle{definition}
\newtheorem{conv}[thm]{Conventions}

\newtheorem{rem}[thm]{Remark}

\newcommand{\refthm}[1]{Theorem~\ref{#1}}
    
    \newcommand{\refconj}[1]{Conjecture~\ref{#1}}
    \newcommand{\refcor}[1]{Corollary~\ref{#1}}
    
    \newcommand{\reflem}[1]{Lemma~\ref{#1}}
    \newcommand{\refprop}[1]{Proposition~\ref{#1}}
    \newcommand{\refrem}[1]{Remark~\ref{#1}}

\newcommand{\HS}{\noindent \hfill{$\sqr55$}}

\begin{document}

\title{Reshetikhin--Turaev invariants of Seifert $3$--manifolds\\ for classical
simple Lie algebras, and their\\ asymptotic expansions}

\shorttitle{Reshetikhin--Turaev invariants of Seifert $3$--manifolds}

\author{S\o ren Kold Hansen}
\address{School of Mathematics, The University of Edinburgh, JCMB, King's Buildings, Edinburgh EH9 3JZ, United Kingdom}
\email{hansen@maths.ed.ac.uk}

\secondauthor{Toshie Takata}
\secondaddress{Department of Mathematics, Faculty of Science, Niigata University, Niigata 950-2181, Japan}
\secondemail{takata@math.sc.niigata-u.ac.jp}

\begin{abstract}
We derive explicit formulas for
the Reshetikhin--Turaev invariants of all oriented Seifert manifolds
associated to an arbitrary complex finite dimensional simple Lie algebra $\mathfrak g$
in terms of the Seifert invariants and standard data for $\mathfrak g$. A main corollary
is a determination of the full asymptotic expansions of these invariants
for lens spaces in the limit of large quantum level.
Our results are in agreement with the asymptotic
expansion conjecture due to J.\ E.\ Andersen \cite{Andersen1}, \cite{Andersen2}.
\end{abstract}




\maketitlepage

\def\undersmile#1{\lower5.7pt\hbox{$\smallsmile$}\kern-0.6em #1}

\section{Introduction}

In 1988 E.\ Witten \cite{Witten} proposed new invariants
$Z_{k}^{G}(X,L) \in \C$ of an
arbitrary $3$--manifold $X$ with
an embedded colored link $L$ by
quantizing the Chern--Simons field theory associated to a
simple and simply connected compact Lie group $G$,
$k$ being an arbitrary positive integer called the (quantum) level.
(Here and in the following a $3$--manifold is
a closed oriented $3$--manifold. In particular
a Seifert manifold is an oriented Seifert manifold.)
The invariant $Z_{k}^{G}(X,L)$ is given by a Feynman path integral 
over the (infinite dimensional) space of gauge equivalence
classes of connections in a $G$ bundle over $X$.
This integral should be understood in a formal way since,
at the moment of
writing, it seems that no mathematically rigorous definition is known,
cf.\ \cite[Sect.~20.2.A]{JohnsonLapidus}.
We will call the invariants $Z_{k}^{G}$ for the quantum
$G$--invariants or the Witten invariants associated to
$G$. 

N.\ Reshetikhin and V.\ G.\ Turaev
\cite{ReshetikhinTuraev} were the first to give a rigorous construction
of quantum invariants of $3$--manifolds with embedded knots.
In fact they constructed invariants $\tau_{r}^{\frsl_{2}(\C)}(X,L) \in \C$
of the pair $(X,L)$ by combinatorial means using surgery presentations of
$(X,L)$ and irreducible representations of the quantum deformations of
$\frsl_{2}(\C)$
at certain roots of unity, $r$ being an integer $\geq 2$ associated
to the order of the root of unity.
Later quantum invariants $\tau_{r}^{\frg}(X,L) \in \C$ associated to
other complex simple Lie algebras $\frg$ were constructed using 
representations of the quantum deformations of $\frg$ at
`nice' roots of unity, see \cite{TuraevWenzl}. 
We call $\tau_{r}^{\frg}$ for the quantum $\frg$--invariants or
the RT--invariants associated to $\frg$.

Both in Witten's approach and in the
approach of Reshetikhin and Turaev the invariants are part of a 
family of topological quantum field theories (TQFT's). This implies that the
invariants are defined for compact oriented $3$--dimensional cobordisms 
(perhaps with some extra structure on the boundary),
and satisfy certain cut-and-paste axioms, see \cite{Atiyah}, \cite{Blanchetetal},
\cite{Quinn}, \cite{Turaev}. The TQFT's of Reshetikhin and Turaev can from
an algebraic point of view be given a more general formulation by using
so-called modular (tensor) categories \cite{Turaev}. The representation theory
of the quantum deformations of $\frg$ at certain roots of unity,
$\frg$ an arbitrary finite dimensional
complex simple Lie algebra, induces such modular categories,
see e.g.\ \cite{Kirillov}, \cite{BakalovKirillov}, \cite{Le2}.

The invariants of Witten are defined by means of a
path integral as stated above. A natural way in physics to obtain
information about quantities defined by means of such path integrals is to
study their perturbative (or asymptotic) expansion for large level $k$.
In fact, by using stationary phase approximation techniques together with path 
integral arguments Witten was able \cite{Witten}
to express the leading large $k$ asymptotics (or the so-called semiclassical approximation) 
of $Z_{k}^{G}(X)$
as a sum over the set
of stationary points for the Chern--Simons functional.
The terms in this sum are expressed by such 
topological/geometric invariants as Chern--Simons invariants, Reidemeister
torsions and spectral flows, so here we see a way to extract topological
information from the invariants. 
A full asymptotic expansion of Witten's invariant
is expected on the basis of a full perturbative 
analysis of the Feynman path integral, see \cite{AxelrodSinger1},
\cite{AxelrodSinger2}, \cite{Axelrod}.

The first rigorous
verifications of the conjectured formula for the semiclassical approximation were
given, partly by Freed and Gompf \cite{FreedGompf}
presenting a large amount of computer calculations for the $SU(2)$--invariants
of lens spaces and some $3$--fibered Seifeirt manifolds,
and about the same time by Jeffrey \cite{Jeffrey2} and Garoufalidis \cite{Garoufalidis} who 
independently gave exact calculations of the semiclassical approximation
of the $SU(2)$--invariants of lens spaces. Jeffrey also verified parts of
the conjecture for the semiclassical approximation of $Z_{k}^{G}(X)$
for $G$ arbitrary and $X$ belonging to a class of mapping tori of the torus.

It is generally believed that the family of TQFT's of Reshetikhin and Turaev
is a mathematical realization of Witten family of TQFT's. 
This belief has together
with the above works on the perturbative expansion of Witten's invariants led
to a detailed conjecture, the asymptotic expansion conjecture (AEC),
which specifies the asymptotic behaviour of the RT--invariants.
The AEC was proposed by Andersen in \cite{Andersen1},
where he proved it for mapping
tori of finite order diffeomorphisms of orientable surfaces of
genus at least two
using the gauge theoretic approach to the quantum invariants.

Let us give an outline of the results obtained in this paper.
Firstly, we determine formulas
for the invariants $\tau_{r}^{\frg}$ of all Seifert manifolds in terms
of the Seifert invariants and standard data for $\frg$, $\frg$ being
an arbitrary complex finite dimensional simple Lie algebra,
cf.\ \refthm{Lie-Seifert}. \refthm{Lie-Seifert} is a generalization
of \cite[Theorem 8.4]{Hansen2}. For a certain subclass of the
Seifert manifolds, containing all the Seifert fibered integral homology spheres,
we simplify the expression for the invariants considerable,
cf.\ \refthm{Lie-Seifert-coprime}.
This result is a generalization of \cite[Formula (4.2)]{LawrenceRozansky}.
Secondly, we analyse more carefully the invariants $\tau_{r}^{\frg}(X)$ for $X$
any lens space. \refthm{Lie-lens} gives the result for any lens space
and any of the invariants $\tau_{r}^{\frg}$. \refprop{Lie-lens-coprime}
gives more compact expressions for $\tau_{r}^{\frg}(L(p,q))$ in case
$r$ and $p$ are coprime. A main corollary of \refthm{Lie-lens},
\refcor{lens-asymp}, is a  determination of
the large $r$ asymptotics of the quantum $\frg$--invariants of the
lens spaces.
The result is in agreement with the AEC, and leads together with the AEC to
a \refconj{Chern--Simons} for the Chern--Simons invariants of the flat $G$
connections on any lens space, $G$ being an arbitrary simply connected, compact simple
Lie group. All the results for lens spaces are generalizations of results
in \cite{Jeffrey2}, which considered the $\frg=\frsl_{2}(\C)$ case. To be precise
\refthm{Lie-lens} generalizes \cite[Theorem 3.4]{Jeffrey2}, \refprop{Lie-lens-coprime}
generalizes \cite[Theorem 3.7]{Jeffrey2}, and \refcor{lens-asymp} is
a generalization of \cite[Formula (5.7)]{Jeffrey2}.

We have via recent private communication learned that J.\ E.\ Andersen has for
the groups $G=SU(n)$ proved the asymptotic expansion conjecture for all 
closed $3$--manifolds via the gauge theoretic approach.
The proof involves asymptotics of Hitchin's connection
over Teichm{\"u}ller space, approximations to all orders, of the boundary
states of handle bodies and techniques similar to the ones presented in 
\cite{Andersen3}. Where Andersen works with the gauge theoretic
definition of the quantum invariant, we work with the definition of
Reshetikhin and Turaev and our proof of the AEC for lens spaces
is very different from Andersen's general proof in the $SU(n)$--case.

A major part of the paper is concerned with studying
a certain family of finite dimensional complex representations $\mR_{r}^{\frg}$
of $\SL (2,\Z)$, one representation for each pair $\frg,r$.
These representations
are known from the study of theta functions and modular forms in
connection with the study of affine Lie algebras, cf.\ \cite{KacPeterson},
\cite[Sect.~13]{Kac}.
They also play a fundamental role in
conformal field theory and (therefore) in the Chern--Simons TQFT's of Witten,
see e.g.\ \cite{GepnerWitten}, \cite{Verlinde}, \cite{Witten}.
In case $\frg=\frsl_{2}(\C)$, Jeffrey \cite{Jeffrey1}, \cite{Jeffrey2} has
determined a nice formula for $\mR_{r}^{\frg}(U)$ in terms of the entries in
$U \in \SL (2,\Z)$ and the integer $r$. \refthm{thm:main1}
is a generalization of Jeffrey's result to arbitrary $\frg$, compare with
\cite[Sect.~2]{Jeffrey2}.
The representations $\mR_{r}^{\frg}$ are of interest when calculating 
the RT--invariants of the Seifert manifolds since certain matrices,
which can be expressed through these representations,
enter into the formulas of the invariants.

The paper is organized as follows. In Sect.~\ref{sec-Formulas-for}
we derive formulas for the representations $\mR_{r}^{\frg}$.
Sect.~\ref{sec-Seifert-manifolds}
is a short section intended to introduce
notation for the Seifert manifolds. Moreover, we recall surgery
presentations for these manifolds due to Montesinos \cite{Montesinos}. In 
Sect.~\ref{sec-The-RT--invariants} we recall the formulas
for the RT--invariants of the Seifert manifolds for an arbitrary modular
category.
These formulas are then used together with the results in
Sect.~\ref{sec-Formulas-for} to
calculate the $\frg$--invariants of the Seifert manifolds.
In Sect.~\ref{sec-The-case}
we analyse the case of lens spaces more carefully.
In the final Sect~\ref{sec-A-rational} we state a rational surgery
formula for the quantum $\frg$--invariants, specializing a rational
surgery formula \cite[Theorem 5.3]{Hansen2}
for the RT--invariants associated to an arbitrary
modular category.
Besides an appendix is added presenting
two related proofs of a reciprocity formula for Gauss sums,
\refprop{prop:gauss}, which plays a vital role in this paper.

After having finished this work, we learned via private communication
that Marino has obtained a similar result as
ours \refthm{Lie-Seifert-coprime} in case the Siefert manifold has
base $S^{2}$ and the Lie algebra $\frg$ is simply laced,
cf.\ \cite[Formula (4.11)]{Marino}.
It seems that the result \refthm{thm:main1} (in the form of the first
formula in \refcor{cor:unitarity}) has been known in the
mathematical physics literature for some time at least for the simply laced case,
cf.\ \cite[Formula (2.5)]{Marino}, \cite[Formula (1.6)]{Rozansky1}.
We have, however, not been able to find any proof of this result in the
literature.

This paper is an extensive expansion of the paper \cite{HansenTakata1}
and gives also the details left out in that paper, in particular
the proof of the main \reflem{lem:main}.

\rk{Acknowledgements} This work were done while the first author was
supported by a Marie Curie Fellowship of the European Commission (CEE
$\text{N}^{o}$ HPMF-CT-1999-00231) and a postdoctoral fellowship
of the European Commision research network EDGE
(European Differential Geometry Endeavour).
He acknowledge hospitality of
l'Institut de Recherche Math\'{e}matique Avanc\'{e}e (IRMA), Universit\'{e}
Louis Pasteur and C.N.R.S., Strasbourg, while being a Marie Curie Fellow,
and the School of Mathematics, University of Edinburgh,
while being an EDGE postdoctoral fellow.
Much of this work were done while the second author visited IRMA
and she also thanks this department for its hospitality.
A part of this work was done while both authors visited the Research Institute
for Mathematical Sciences (RIMS), Kyoto University.
We would like to thank RIMS for hospitality during the special month on 
{\it invariants of knots and $3$--manifolds}, September 2001.
We also thank the organizers of this workshop for financial support
during our visit.
The first author thanks J.\ E.\ Andersen for helpful conversations about quantum
invariants in general and about asymptotics of these invariants in particular,
and both authors would like to thank T.\ T.\ Q.\ Le for some helpful comments.
The first author is partially supported by the Danish Natural Science Research Council.

\section{Formulas for the representations of $\SL(2,\Z)$.}\label{sec-Formulas-for}

In this section we analyse a family of unitary representations of $\SL(2,\Z)$
associated with a complex finite dimensional simple Lie algebra $\frg$.
First let us fix some notation for $\frg$. (For details about standard material
for Lie algebras we refer to \cite{Humphreys1}.) Let $\frh$ be a Cartan subalgebra of
$\frg$, and let $\frhR^{*}$ be the $\R$--vector space spanned by the roots.
We let $X$ and $Y$ be the weight lattice and the root lattice respectively.
Let $\la \cdot,\cdot \ra$ be the standard inner product in $\frhR^{*}$ 
normalized such that the long roots have length $\sqrt{2}$. (That is,
$\la \cdot,\cdot \ra$ is proportional to the inner product induced by the Killing
form of $\frg$ restricted to $\frh$.) Then the short roots have length
$\sqrt{2/m}$, where $m=1$ if $\frg$ is simply laced, i.e.\ belongs to
the series $A$, $D$ or $E$, $m=2$ if $\frg$ belongs to the series $B$ or $C$
or is of type $F_{4}$, and $m=3$ if $\frg$ is of type $G_{2}$.
In the following, the symbol $m$ will be reserved for this number.
For $x \in \frhR^{*} \sm \{0\}$,
we let $x^{\vee}=2x/\la x,x\ra$.
If $\Pi=\{\alpha_{1},\alpha_{2},\ldots,\alpha_{l} \}$ is an arbitrary set of
simple (basis) roots and if
$\{\lambda_{1},\lambda_{2},\ldots,\lambda_{l}\} \subseteq \frhR^{*}$
is the set of fundamental dominant weights relative to $\Pi$, i.e.\
$\la \lambda_{i} , \alpha_{j}^{\vee} \ra = \delta_{ij}$, then
$X$ is the $\Z$--lattice generated by
$\{\lambda_{1},\lambda_{2},\ldots,\lambda_{l}\}$ and $Y$ is the $\Z$--lattice
generated by $\Pi$.
Let $\Delta$ be the set of roots.
Then $\{\;\alpha^{\vee}\;|\;\alpha \in \Delta\;\}$
are the so-called dual roots or coroots
(relative to our inner product $\la\cdot,\cdot\ra$).
(So in this paper a coroot is in $\frhR^{*}$ and not in $\frh$.) 
The coroot lattice $\cY$ is the $\Z$--lattice generated by
$\{\alpha_{1}^{\vee},\ldots,\alpha_{l}^{\vee}\}$
for an arbitrary set of simple roots
$\{\alpha_{1},\ldots,\alpha_{l}\}$. Recall that
\begin{eqnarray}\label{eq:integer1}
\cY &=& \{ x \in \frhR^{*} \; | \; \la x , y \ra \in \Z \hspace{.1in}{\text{\rm{for all}}}\hspace{.1in} y \in X \}, \\
X &=& \{ x \in \frhR^{*} \; | \; \la x , y \ra \in \Z \hspace{.1in}{\text{\rm{for all}}}\hspace{.1in} y \in \cY \}. \nonumber 
\end{eqnarray}
We note that $X$ and $\cY$ are dual to each other. It is obvious from the
above, that the Weyl group $W$ preserves the lattices $X$, $Y$, and $\cY$.
Let us also note the following facts:
For all $x,y \in Y^{\vee}$ we have
\begin{equation}\label{eq:integer2}
\la x , y \ra \in \Z, \hspace{.2in} \la x,x\ra \in 2\Z.
\end{equation}
There exists a (least) positive integer $D$ such that we for all
$\mu,\xi \in X$ have
\begin{equation}\label{eq:integer3}
\la \mu , \xi \ra \in \frac{1}{D} \Z, \hspace{.2in} \la \mu, \mu \ra \in \frac{2}{D} \Z.
\end{equation}

Let us fix a set
$\Pi=\{\alpha_{1},\alpha_{2},\ldots,\alpha_{l} \}$ of
simple roots in the following and let
$\{\lambda_{1},\lambda_{2},\ldots,\lambda_{l}\}$
be the set of fundamental dominant weights relative to $\Pi$.
The (open) fundamental Weyl chamber (relative to $\Pi$) is the set
$$
\FWC=\{ x \in \frhR^{*} \; | \; \la x,\alpha_{i} \ra > 0, i=1,\ldots,l\}.
$$
We let $C=\overline{\FWC}$ be the topological closure of $\FWC$.
For a positive integer $k$, the $k$--alcove (relative to $\Pi$)
is the (closed) set
$$
C_{k}=\{ x \in C \; | \; \langle x,\theta \rangle \leq k\},
$$
where $\theta$ is the highest root of $\frg$ (relative to $\Pi$).
(Note that $\theta$ is a long root). We will also need the following
sets:
\begin{eqnarray*}
Q_{k} &=& \{ c_{1}\lambda_{1}+\ldots+c_{l}\lambda_{l} \; | \; c_{1},\ldots,c_{l} \in [0,k[ \; \}, \\
P_{k} &=& \{ c_{1}\alpha_{1}^{\vee}+\ldots+c_{l}\alpha_{l}^{\vee} \; | \; c_{1},\ldots,c_{l} \in [0,k[ \; \}.
\end{eqnarray*}
Let $\aW_{k}=W \ltimes k \cY$ be the affine Weyl group acting on $\frhR^{*}$
in the usual sense ($k \cY$ acting by translations).
It follows that $P_{k}$ is a fundamental domain of the group $k\cY$. Moreover, it
is well-known (see \cite[Sect.~6]{Kac} for a proof) that $C_{k}$ is
a fundamental domain of $\aW_{k}$.
We will say that a subset $M$ of $\frhR^{*}$ is presicely tiled by $k$--alcoves if there exists
a family $\{ u_{i} \}_{i \in I}$ of elements of $\aW_{k}$ such that
$M=\cup_{i \in I} u_{i}(C_{k})$. If $I$ is finite, we say that $M$ is (precisely)
tiled by $|I|$ $k$--alcoves.
We let $\Delta_{+}$ be the set of
positive roots (relative to $\Pi$). For $\alpha \in \Delta_{+}$ and $n \in \Z$ we let
$$
H_{\alpha,n}^{k}=\{ x \in \frhR^{*} \; | \; \la x,\alpha \ra = nk \},
$$
and we put
$$
H^{k}=\cup_{\alpha \in \Delta_{+}, n \in \Z} H_{\alpha,n}^{k}
  =\{ x \in \frhR^{*} \; | \; \exists \alpha \in \Delta_{+} : \la x,\alpha \ra \in k\Z \; \}.
$$ 
For fixed $\lambda_{0},\lambda_{1} \in X$ and fixed integers $b$ and $a \neq 0$ we let
\begin{eqnarray}\label{eq:gfunction}
g_{\lambda_{0},\lambda_{1}}^{a,b,k}(\lambda) &=& \sum_{\mu \in \cY/a \cY} \sum_{w,w' \in W} \det (w w')
    \exp\left( \frac{\pi\I b}{ak}|\lambda+k \mu|^2 \right.\\
  &&\hspace{.6in} \left. + \frac {2\pi\I }{a k}
       \la \lambda + k\mu, -w(\lambda_{0})-a w'(\lambda_{1}) \ra \right) \nonumber
\end{eqnarray}
for $\lambda \in X$, where $|\cdot|$ is the norm associated to $\la\cdot,\cdot\ra$.
(Note that by (\ref{eq:integer1}) and (\ref{eq:integer2}) the summand in the
expression of $g_{\lambda_{0},\lambda_{1}}^{a,b,k}$ only
depends on $\mu \pmod{a\cY}$.) In particular,
$$
g_{\lambda_{0},\lambda_{1}}^{1,b,k}(\lambda) = \sum_{w,w' \in W} \det (w w') 
     \exp\left( \frac{\pi\I}{k}\left( b|\lambda|^{2}
               + \la \lambda, -2w(\lambda_{0})-2w'(\lambda_{1}) \ra \right)\right).
$$
By using similar arguments as in \cite[Sect.~4]{Jeffrey2}
(see also \cite{Le2}), we have

\begin{prop}\label{prop:inv}
The map $g_{\lambda_{0},\lambda_{1}}^{a,b,k} \co X \to \C$ is invariant under the action of
the affine Weyl group $\aW_{k}$. Moreover, $g_{\lambda_{0},\lambda_{1}}^{a,b,k}(\lambda)=0$
for any $\lambda \in X \cap H^{k}$.
\end{prop}

\begin{proof}
The invariance under the action by an element $u \in W$ follows by the
fact that $u$ is orthogonal together with the identity
$\det (ww')=\det (u^{-1}w u^{-1}w')$
and the fact that $W$ preserves $\cY/a\cY$ (since it preserves $\cY$).
The invariance under the action by an element $k x$, $x \in \cY$,
is obvious. (For $|a| = 1$, use (\ref{eq:integer1}) and (\ref{eq:integer2}).)
To prove the last claim, let $\lambda \in X \cap H_{\alpha,n}^{k}$ for
a positive root $\alpha$ and an integer $n$, and let
$s_\alpha$ be the reflection in $\alpha$. Fix a $w \in W$ and get
\begin{eqnarray*}
\la \lambda+k\mu, -w(\lambda_{0}) - a{s_{\alpha}}^{-1}w'(\lambda_{1}) \ra
 &=& \la \lambda+k\mu, -w(\lambda_{0})-w'(\lambda_{1}) \ra \\
 &&\hspace{.2in} + a\la \lambda+k\mu,\alpha\ra \la \alpha^{\vee}, w'(\lambda_{1}) \ra,
\end{eqnarray*}
where $\la \alpha^{\vee}, w'(\lambda_{1}) \ra \in \Z$, and
$\la \lambda+k\mu,\alpha\ra \in k\Z$,
so
$$
\sum_{w' \in W} \det(w')\exp\left(\frac {2\pi\I }{a k}
    \la \lambda + k\mu, -w(\lambda_{0})-a w'(\lambda_{1}) \ra \right)=0.
$$
\end{proof}

In the calculations to follow, a multi-dimensional reciprocity
formula for Gauss sums plays a crucial role. 
Let $V$ be a real vector space of dimension $l$ with inner product
$\la \cdot, \cdot\ra$, $\Lambda$ a lattice in $V$ and $\Lambda^{*}$
the dual lattice.
For an integer $r$, a self-adjoint automorphism
$f \co  V \to V$, and an element  $\psi \in V$,
we assume
\begin{eqnarray}\label{eq:assumption}
&&\frac{r}{2} \la \lambda,f(\lambda) \ra, \quad
\la \lambda,f(\eta) \ra, \quad
r \la \lambda, \psi \ra \in \Z, \hspace{.2in} \forall \lambda, \eta \in \Lambda,\\
&&\frac{r}{2} \la \mu,f(\mu) \ra, \quad
\la \mu,r \xi \ra, \quad
r \la \mu, \psi \ra \in \Z, \hspace{.2in} \forall \mu, \xi \in \Lambda^*, \nonumber
\end{eqnarray}
and $f(\Lambda^{*}) \subseteq \Lambda^{*}$. Then we have

\begin{prop}(Reciprocity formula for Gauss sums)\label{prop:gauss}
\begin{eqnarray*}
& & \vol (\Lambda^{*}) \sum_{\lambda \in \Lambda / r\Lambda}
     \exp \left(\frac{\pi\I}{r} \la \lambda,f(\lambda) \ra\right)
     \exp \left(2\pi\I \la \lambda,\psi \ra\right) \\
 && \hspace{.3in} =\left( \det \frac f \I \right)^{-1/2} r^{l/2}
     \sum_{\mu \in \Lambda^{*} / f(\Lambda^{*})}
      \exp \left( -\pi r \I \la \mu + \psi, f^{-1}(\mu+\psi) \ra \right).
\end{eqnarray*}
\end{prop}

For a proof, see \cite[Sect.~2]{Jeffrey1}. For the
convenience of the reader, we sketch Jeffrey's proof
in the appendix. Moreover, we present in the appendix
a slightly alternative proof. Both arguments rely on the Poisson
resummation formula.
Below we will use the reciprocity formula, Proposition \ref{prop:gauss},
with $\Lambda =X$, the dual lattice being the
coroot lattice $\cY$.

The group $\SL(2,\Z)$ is generated by two matrices
\begin{equation}\label{eq:generators}
\Xi = \left( \begin{array}{cc}
		0 & -1 \\
		1 & 0
		\end{array}
\right),\hspace{.2in}
\Theta = \left( \begin{array}{cc}
                1 & 1 \\
		0 & 1
		\end{array}
\right).
\end{equation}
We note that
\begin{equation}\label{eq:relationsXiTheta}
\Xi^{2}=(\Xi\Theta)^{3}=-1.
\end{equation}
For a tuple of integers $\mC=(m_{1},\ldots,m_{t})$
we let $\mC_{k}=(m_{1},\ldots,m_{k})$ and
\begin{equation}\label{eq:Bmatrix}
B_{k}^{\mC} = \left( \begin{array}{cc}
                a_{k}^{\mC} & b_{k}^{\mC} \\
		c_{k}^{\mC} & d_{k}^{\mC}
		\end{array}
\right) = \Theta^{m_{k}}\Xi\Theta^{m_{k-1}}\Xi \ldots \Theta^{m_{1}}\Xi
\end{equation}
for $k=1,2,\ldots,t$,
and let $B^{\mC}=B_{t}^{\mC}$. Moreover, we put
$$
a_{0}^{\mC}=d_{0}^{\mC}=1, \qquad b_{0}^{\mC}=c_{0}^{\mC}=0.
$$
We say that $\mC$ has length $|\mC|=t$.
If it is clear from the context what $\mC$ is
we write $a_{k}$ for $a_{k}^{\mC}$ etc.
By \cite[Proposition 2.5]{Jeffrey2} the elements $a_{i},b_{i},c_{i},d_{i}$
satisfy the recurrence relations
\begin{align}\label{eq:abcd}
a_{k}=m_{k} a_{k-1}-c_{k-1},& \quad c_{k}=a_{k-1},\\
b_{k}=m_{k} b_{k-1}-d_{k-1},& \quad d_{k}=b_{k-1} \nonumber
\end{align}
for $k=1,2,\ldots,t$.
Moreover,
\begin{equation}\label{eq:b/a}
\frac{b_{k}}{a_{k}} = - \left( \frac{1}{a_{1}} + \frac{1}{a_{2}a_{1}}+ \cdots
            + \frac{1}{a_{k}a_{k-1}}
              \right)
\end{equation}
and $(m_{1},\ldots,m_{k})$ is a continued fraction expansion of
$a_{k}/b_{k}$, $k=1,2,\ldots,t$, i.e.\
\begin{equation}\label{eq:continued}
\frac{a_{k}}{b_{k}}=m_{k}-\frac{1}{m_{k-1}-\dfrac{1}{\cdots -\dfrac{1}{m_{1}}}}.
\end{equation}
Let in the following $\kappa \in \Z_{>0}$ be fixed.
We have a representation $\mR=\mR_{\kappa}^{\frg}$ of $\SL(2,\Z)$ given by
\begin{eqnarray}\label{eq:mR}
&&\mR(\Xi)_{\lambda \mu} = \frac{\I^{\npr}}{\kappa^{l/2}}
        \left| \frac{\vol (X) }{\vol (\cY) }\right|^{1/2}
        \sum_{w \in W} \det(w)
        \exp \left( -\frac{2\pi\I}{\kappa} \la w(\lambda), \mu\ra \right), \nonumber \\
&&\mR(\Theta)_{\lambda \mu} = \delta_{\lambda \mu}
        \exp \left( \frac{\pi\I}{\kappa} |\lambda|^{2}
           -\frac{\pi\I}{\ch} |\rho|^{2} \right)
\end{eqnarray}
for $\lambda, \mu \in \inte(C_{\kappa}) \cap X$. Here $\rho$
is half the sum of positive roots.
In the following we also write $\tilde{U}$ for $\mR(U)$.
Note that the matrices $\mR(\Xi)$ and $\mR(\Theta)$ are
symmetric (in fact $\mR(\Theta)$ is diagonal).

\begin{rem}\label{rem:unitarity}
It is well-known that the representations $\mR$ are unitary. Let
us here give a few references to this fact. Unitarity of
$\mR(\Theta)$ follows immediately from (\ref{eq:mR}).
By \cite[Theorem 3.3.20]{BakalovKirillov} we have
\begin{equation}\label{eq:unitarity1}
\mR(\Xi)_{\lambda\mu}=\overline{s_{\lambda\mu}},
\end{equation}
where $s$ is defined in \cite[Formula (3.1.16)]{BakalovKirillov},
and $\bar{\cdot}$ is complex conjugation. (In the expression
for $s_{\lambda\mu}$ in \cite[Theorem 3.3.20]{BakalovKirillov} there
is a minor mistake; one has to replace $i^{\npr}$ by $i^{-\npr}$.)
By the proof of \cite[Theorem 3.3.20]{BakalovKirillov} it follows that
\begin{equation}\label{eq:unitarity2}
s_{\lambda\mu}^{2}=\delta_{\lambda^{*}\mu}.
\end{equation}
Here $\lambda^{*}=-w_{0}(\lambda-\rho)+\rho=-w_{0}(\lambda)$, where $w_{0}$
is the longest element in $W$ (relative to our set of simple roots $\Pi$)
and where we use that $w_{0}(\rho)=-\rho$,
see \cite[Sect.~1.8]{Humphreys2}. By \cite[Sect.~1.8]{Humphreys2} we have
$\det(w_{0})=(-1)^{\npr}$. By this and (\ref{eq:mR}) we get
\begin{equation}\label{eq:unitarity3}
\mR(\Xi)_{\lambda^{*}\mu}=\overline{\mR(\Xi)_{\lambda\mu}}.
\end{equation}
Now unitarity of $\mR(\Xi)$ follows from
(\ref{eq:unitarity1}), (\ref{eq:unitarity2}), and (\ref{eq:unitarity3}).
The original reference for the unitarity of $\mR$ seems to be
\cite[Appendix]{GepnerWitten}.
\end{rem}

One should note that the expressions for the entries of $\mR(\Xi)$
and $\mR(\Theta)$ are well-defined for all $\lambda, \mu \in X$.
Note also that if $\lambda$ or $\mu$ is an elements of $X$
belonging to the boundary of $C_{\kappa}$, then $\mR(\Xi)_{\lambda \mu}=0$
(use the same argument as in the proof of the final claim of \refprop{prop:inv}).
This observation allows us to shift between $\inte(C_{\kappa}) \cap X$
and $C_{\kappa} \cap X$ as summation index set in formulas below.
Following Jeffrey \cite[Sect.~2]{Jeffrey1}, \cite[Sect.~2]{Jeffrey2} we consider
$$
\mT^{\mC}_{\lambda_{0},\lambda_{t+1}} = 
\sum_{\lambda_{1},\ldots,\lambda_{t} \in C_{\kappa} \cap X}
\tilde{\Xi}_{\lambda_{t+1} \lambda_{t}}
\tilde{\Theta}_{\lambda_t}^{m_t} \tilde{\Xi}_{\lambda_{t} \lambda_{t-1}}
\tilde{\Theta}_{\lambda_{t-1}}^{m_{t-1}} \tilde{\Xi}_{\lambda_{t-1} \lambda_{t-2}}
\cdots
\tilde{\Theta}_{\lambda_{1}}^{m_{1}} \tilde{\Xi}_{\lambda_{1} \lambda_{0}}
$$
for $\lambda_{0},\lambda_{t+1} \in X$,
where we write $\tilde{\Theta}_{\lambda}$ for $\tilde{\Theta}_{\lambda\lambda}$.
Then we have the following generalization of \cite[Lemma 2.6]{Jeffrey2}:

\begin{lem}\label{lem:main}
Assume that $\mC=(m_{1},\ldots,m_{t})$ is a sequence of integers such that
$a_{k}$ is nonzero for $k=1,\ldots,t$. Then
$$
\mT^{\mC}_{\lambda_{0},\lambda_{t+1}} = K^{\mC}_{\lambda_{0}} \sum_{w\in W} \det(w)
      \sum_{\mu \in \cY/a_{t} \cY}
       \exp\left( -\frac {\pi\I c_{t}}{a_{t} \kappa}
       \left|\lambda_{t+1}+\kappa\mu + \frac {w(\lambda_{0})}{c_{t}}\right|^{2} \right)
$$
for all $\lambda_{0},\lambda_{t+1} \in X$, where
\begin{eqnarray*}
K^{\mC}_{\lambda_{0}} &=& \frac{{\I}^{(t+1) \npr}}{(\kappa|a_{t}|)^{l/2}\vol (\cY)}\; \zeta^{l \; D_{t}}
 \exp \left( -\frac{\pi\I}{\ch} (\sum_{i=1}^{t} m_{i}) |\rho |^{2} \right) \\
 & & \hspace{.3in} \times \exp \left( -\frac{\pi\I}{\kappa}
        (\sum_{i=1}^{t-1} \frac{1}{a_{i-1}a_{i}}) |\lambda_{0} |^{2} \right).
\end{eqnarray*}
Here $\zeta=\exp \frac{\pi \I}{4}$
and $D_{t}=\sign(a_{0}a_{1})+\cdots +\sign(a_{t-1}a_{t})$.
\end{lem}

\begin{proof}
We prove the proposition by induction on the length of $\mC$. First consider
$\mT^{\mC_{1}}_{\lambda_{0},\lambda_{2}}$, $\lambda_{0},\lambda_{2} \in X$.
Since $\tilde{\Xi}$ is symmetric we have
\begin{eqnarray*}
\mT^{\mC_{1}}_{\lambda_{0},\lambda_{2}} &=& \sum_{\lambda_{1} \in C_{\kappa} \cap X } \tilde{\Xi}_{\lambda_{2} \lambda_{1}}
    \tilde{\Theta}_{\lambda_{1}}^{m_{1}} \tilde{\Xi}_{\lambda_{1} \lambda_{0}} \\
 &=& \frac{(-1)^{\npr}}{\kappa^{l}} \frac{\vol (X)}{\vol (\cY)}
     \exp \left( - \frac{\pi \I}{\ch} m_{1} |\rho |^{2} \right) \sum_{\lambda \in C_{\kappa} \cap X}
        g_{\lambda_{0},\lambda_{2}}^{1,m_{1},\kappa}(\lambda),
\end{eqnarray*}
where $g_{\lambda_{0},\lambda_{2}}^{1,m_{1},\kappa}$ is given by (\ref{eq:gfunction}).
The closure $\bar{Q}_{m}$ of $Q_{m}$ is precisely tiled by $1$--alcoves.
Moreover, if $\bar{Q}_{N}$ is tiled by $n$ $1$--alcoves, then $\bar{Q}_{kN}$ is
tiled by $n$ $k$--alcoves. We also have that if $\bar{Q}_{N}$ is tiled by $n$ $1$--alcoves,
then $\bar{Q}_{kN}$ is tiled by $k^{l}n$ $1$--alcoves.
Note that $\vol(\bar{Q}_{N}) = N^{l}\vol(X)$ and $\vol(C_{1})=\vol(\bar{P}_{1})/|W|$.
We use here that $C_{1}$ is a fundamental domain for the action of $\aW_{1}$, while
$P_{1}$ is a fundamental domain for the action of $\cY$.
Therefore, if $\bar{Q}_{N}$ is precisely tiled by $n$ $1$--alcoves, then
$n=\vol(\bar{Q}_{N})/\vol(C_{1})=N^{l}|W|\vol(X)/\vol(\cY)$.
Let $N=mD$, where $D$ is the integer from (\ref{eq:integer3}).
By \refprop{prop:inv} and the above we get
$$
\sum_{\lambda \in X/N\kappa X} g_{\lambda_{0},\lambda_{2}}^{1,m_{1},\kappa}(\lambda)
= N^{l} |W| \vol (X) / \vol (\cY) \sum_{\lambda \in C_{\kappa} \cap X}
g_{\lambda_{0},\lambda_{2}}^{1,m_{1},\kappa}(\lambda).
$$
Therefore
\begin{eqnarray*}
\mT^{\mC_{1}}_{\lambda_{0},\lambda_{2}} &=& \frac{(-1)^{\npr}}{\kappa^{l} N^{l} |W|}
     \exp \left( - \frac{\pi \I}{\ch} m_{1} |\rho |^{2} \right) \sum_{w,w' \in W}\det (ww') \\
 && \hspace{.1in}\times \sum_{\lambda \in X / N\kappa X}
      \exp \left( \frac{\pi \I}{\kappa}
         \left( m_{1}|\lambda|^{2} +
                \la \lambda, -2w(\lambda_{0})-2w'(\lambda_{2}) \ra \right) \right).
\end{eqnarray*}
Let $f=m_{1}N\id_{\frhR^{*}}$,
$\psi=\psi(w,w')=-\frac{1}{\kappa}(w(\lambda_{0})+w'(\lambda_{2}))$, $w,w' \in W$, and
$r=N\kappa$. Then, by (\ref{eq:integer1}), (\ref{eq:integer2}),
and (\ref{eq:integer3}), the assumptions (\ref{eq:assumption}) are satisfied,
so by \refprop{prop:gauss} we get
\begin{eqnarray*}
\mT^{\mC_{1}}_{\lambda_{0},\lambda_{2}} &=& \frac{(-1)^{\npr}}{r^{l}|W|}
     \exp \left( - \frac{\pi \I}{\ch} m_{1} |\rho |^{2} \right) \sum_{w,w' \in W}\det (ww') \\
 && \hspace{.1in} \times \sum_{\lambda \in X / r X}
      \exp \left( \frac{\pi \I}{r} \la \lambda, f (\lambda) \ra \right)
          \exp \left( 2\pi \I\la \lambda, \psi \ra \right) \\
 &=& \vol(\cY)^{-1} \frac{(-1)^{\npr}}{r^{l}|W|} \det \left( \frac{f}{\I} \right)^{-1/2} r^{l/2}
       \exp \left( -\frac{\pi \I}{\ch} m_{1} |\rho |^{2} \right) \\
 && \hspace{.1in}\times \sum_{w,w' \in W} \det (ww')
      \sum_{\mu \in \cY/f(\cY)}
      \exp\left( -\pi \I r \la \mu +\psi, f^{-1}(\mu+\psi)\ra\right) \\
 &=& \frac{(-1)^{\npr}}{(\kappa N)^{l/2}\vol(\cY)|W|} \left( \frac{\I}{m_{1}N} \right)^{l/2}
       \exp \left( -\frac{\pi \I}{\ch} m_{1} |\rho |^{2} \right) \\
 && \hspace{.1in}\times \sum_{w,w' \in W} \det (ww')
      \sum_{\mu \in \cY/m_{1}N\cY}
      \exp\left( -\frac{\pi \I \kappa}{m_{1}} |\mu + \psi|^{2}\right).
\end{eqnarray*}
Here $\exp\left( -\frac{\pi \I \kappa}{m_{1}} |\mu + \psi|^{2} \right)$
is invariant under $\mu \mapsto \mu + m_{1}\alpha$, $\mu,\alpha \in \cY$, by
(\ref{eq:integer1}) and (\ref{eq:integer2}), so
{\allowdisplaybreaks
\begin{eqnarray*}
\mT^{\mC_{1}}_{\lambda_{0},\lambda_{2}} &=&  \frac{(-1)^{\npr}}{(\kappa)^{l/2}\vol(\cY)|W|} \left( \frac{\I}{m_{1}} \right)^{l/2}
       \exp \left( -\frac{\pi \I}{\ch} m_{1} |\rho |^{2} \right) \\*
 && \hspace{.3in}\times \sum_{w,w' \in W} \det (ww')
      \sum_{\mu \in \cY/m_{1}\cY}
      \exp\left( -\frac{\pi \I \kappa}{m_{1}} |\mu +\psi|^{2}\right) \\
 &=& \frac{(-1)^{\npr}}{\kappa^{l/2}\vol (\cY) |W|}
      \left( \frac{\I}{m_{1}} \right)^{l/2}
      \exp \left( -\frac{ \pi \I}{\ch} m_{1}|\rho |^{2} \right) \sum_{w,w' \in W} \det(ww') \\*
 && \hspace{.3in} \times \sum_{\mu \in \cY / m_{1} \cY}
      \exp\left( -\frac{\pi \I}{m_{1} \kappa}
           |\kappa {w'}^{-1}(\mu) -{w'}^{-1}w(\lambda_{0}) -\lambda_{2} |^{2} \right) \\
&=& \frac{(-1)^{\npr}}{\kappa^{l/2}\vol (\cY) |W|}
      \left( \frac{\I}{m_{1}} \right)^{l/2}
      \exp \left( -\frac{ \pi \I}{\ch} m_{1}|\rho|^{2} \right) \\*
 && \hspace{.3in} \times \sum_{w' \in W} \sum_{w \in W} \det ({w'}^{-1}w) \\*
 && \hspace{.3in} \times \sum_{\mu \in \cY / m_{1} \cY}
      \exp\left( -\frac{\pi \I}{m_{1} \kappa}
           |\kappa \mu -{w'}^{-1}w(\lambda_{0}) -\lambda_{2} |^{2} \right) \\
&=& \frac{(\I)^{2\npr}}{\kappa^{l/2}\vol (\cY) }
      \left( \frac{1}{\sqrt{|m_{1}|}} \right)^l \zeta^{l\; \sign(m_{1})}
      \exp \left( -\frac{m_{1} \pi \I}{\ch} |\rho |^{2} \right) \\*
 && \hspace{.3in} \times \sum_{w \in W} \det (w)
      \sum_{\mu \in \cY / m_{1} \cY}
      \exp\left( -\frac{\pi \I}{m_{1} \kappa}
           |\kappa \mu -w(\lambda_{0}) -\lambda_{2} |^{2} \right).
\end{eqnarray*}}\noindent
Finally we replace $\mu$ by $-\mu$, and use that $m_{1}=a_{1}$, $c_{1}=1$.

Assume next by induction that the lemma is true for all
sequences of length $t-1$. Then we get for $\lambda_{0},\lambda_{t+1} \in X$ that
{\allowdisplaybreaks
\begin{eqnarray*}
\mT^{\mC}_{\lambda_{0},\lambda_{t+1}} &=& \sum_{\lambda_{t} \in C_{\kappa} \cap X}
      \tilde{\Xi}_{\lambda_{t+1} \lambda_{t}} \tilde{\Theta}_{\lambda_{t}}^{m_{t}} \mT^{\mC_{t-1}}_{\lambda_{0},\lambda_{t}} \\
&=& K^{\mC_{t-1}}_{\lambda_{0}} \frac{{\I}^{\npr}}{\kappa^{l/2}}
      \left| \frac {\vol (X)}{\vol (\cY)} \right|^{1/2}
      \exp \left( -\frac{\pi \I}{\ch} m_{t} |\rho|^{2} \right) \\*
 && \hspace{.3in} \times \sum_{w,w' \in W} \det(ww') \\
 && \hspace{.3in} \times \sum_{\lambda_{t} \in C_{\kappa} \cap X}
      \exp \left( \frac{\pi \I}{\kappa}
                \left( m_{t}|\lambda_{t}|^{2} +
                \la \lambda_{t}, -2w'(\lambda_{t+1}) \ra \right) \right) \\*
 && \hspace{.3in} \times
      \sum_{\mu \in \cY/a_{t-1} \cY}
      \exp\left( -\frac{\pi\I c_{t-1}}{a_{t-1} \kappa}
      \left|\lambda_{t}+\kappa\mu + \frac {w(\lambda_{0})}{c_{t-1}}\right|^{2} \right).
\end{eqnarray*}}\noindent
By (\ref{eq:integer1}) and (\ref{eq:integer2}) we have
{\allowdisplaybreaks
\begin{eqnarray*}
&&   \sum_{w,w' \in W} \det(ww')
 \sum_{\lambda_{t} \in C_{\kappa} \cap X}
      \exp \left( \frac{\pi \I}{\kappa}
                \left( \la \lambda_{t}, m_{t} \lambda_{t} \ra +
                \la \lambda_{t}, -2w'(\lambda_{t+1}) \ra \right) \right) \\*
 && \hspace{.6in} \times
      \sum_{\mu \in \cY/a_{t-1} \cY}
      \exp\left( -\frac{\pi\I c_{t-1}}{a_{t-1} \kappa}
      \left|\lambda_{t}+\kappa\mu + \frac {w(\lambda_{0})}{c_{t-1}}\right|^{2} \right) \\
 && \hspace{.3in} = \sum_{w,w' \in W} \det (w w')
      \sum_{\lambda_{t} \in C_{\kappa} \cap X} \sum_{\mu \in \cY/a_{t-1} \cY} \\*
 && \hspace{.6in} \times \exp \left( \frac{\pi \I}{\kappa}
      \left(m_{t} |\lambda_{t}+\kappa\mu|^{2} +
      \la \lambda_{t}+\kappa\mu, -2w'(\lambda_{t+1}) \ra \right) \right) \\*
 && \hspace{.6in} \times \exp\left( -\frac{\pi\I c_{t-1}}{a_{t-1} \kappa}
         |\lambda_{t}+\kappa\mu|^{2} \right. \\*
 && \hspace{1.0in} \left. + \frac{2\pi \I}{a_{t-1} \kappa} \la \lambda_{t}+\kappa\mu, -w(\lambda_{0}) \ra -
         \frac{\pi \I}{a_{t-1} c_{t-1}\kappa} |\lambda_{0}|^{2} \right) \\
 && \hspace{.3in} = \exp \left( -\frac{\pi \I}{a_{t-1} c_{t-1}\kappa} |\lambda_{0}|^{2} \right)
      \sum_{\lambda_{t} \in C_{\kappa} \cap X} g_{\lambda_{0},\lambda_{t+1}}^{a_{t-1},a_{t},\kappa}(\lambda_{t}).
\end{eqnarray*}}\noindent
We want to alter the sum
$\sum_{\lambda_{t} \in C_{\kappa} \cap X} g_{\lambda_{0},\lambda_{t+1}}^{a_{t-1},a_{t},\kappa}(\lambda_{t})$
by applying the reciprocity formula, \refprop{prop:gauss}.
Before we can do this we need to alter this sum
using some symmetries.
By \refprop{prop:inv} and the fact that $C_{\kappa}$ and $P_{\kappa}$ are
fundamental domains of respectively $\aW_{\kappa}$ and $\kappa\cY$, we get
$$
|W| \sum_{\lambda \in C_{\kappa} \cap X} g_{\lambda_{0},\lambda_{t+1}}^{a_{t-1},a_{t},\kappa}(\lambda)
 = \sum_{P_{\kappa} \cap X} g_{\lambda_{0},\lambda_{t+1}}^{a_{t-1},a_{t},\kappa}(\lambda).
$$
Since the map $\cY / a_{t-1}\cY \times (P_{\kappa} \cap X) \ria P_{a_{t-1}\kappa} \cap X$,
$(\mu,\lambda) \mapsto \lambda + \kappa\mu$, defines a bijection,
we therefore get
$$
\sum_{\lambda \in C_{\kappa} \cap X} g_{\lambda_{0},\lambda_{t+1}}^{a_{t-1},a_{t},\kappa}(\lambda) = \frac{1}{|W|}
 \sum_{\lambda \in P_{a_{t-1}\kappa} \cap X} h_{\lambda_{0},\lambda_{t+1}}(\lambda),
$$
where
\begin{eqnarray*}
h_{\lambda_{0},\lambda_{t+1}}(\lambda) &=& \sum_{w,w'\in W} \det(ww')\exp\left( \frac{\pi\I a_{t} }{a_{t-1} \kappa} |\lambda|^2 \right. \\
 &&\hspace{.6in} \left. + \frac {2\pi\I }{a_{t-1} \kappa}
       \la \lambda, -w(\lambda_{0})-a_{t-1} w'(\lambda_{t+1}) \ra \right).
\end{eqnarray*}
Exactly as in the proof of \refprop{prop:inv} we get that
$h_{\lambda_{0},\lambda_{t+1}} \co X \to \C$
is invariant under the action of $\aW_{a_{t-1}\kappa}$, and, moreover,
that $h_{\lambda_{0},\lambda_{t+1}}(\lambda)=0$ for all
$\lambda \in X \cap H^{\kappa}$.  
Now let $N=mD$ as in the first part of the induction and get that
$$
\sum_{\lambda \in X/ \kappa a_{t-1} N X} h_{\lambda_{0},\lambda_{t+1}}(\lambda) =\frac{N^{l}|W|\vol(X)}{\vol(\cY)}
\sum_{\lambda \in C_{a_{t-1}\kappa} \cap X} h_{\lambda_{0},\lambda_{t+1}}(\lambda),
$$
where we use that $\bar{Q}_{\kappa a_{t-1}N}$ is precisely tiled by $\kappa a_{t-1}$--alcoves,
in fact by $N^{l}|W|\vol(X)/\vol(\cY)$ of these alcoves. Therefore
$$
\sum_{\lambda \in C_{\kappa} \cap X} g_{\lambda_{0},\lambda_{t+1}}^{a_{t-1},a_{t},\kappa}(\lambda)
= \frac{\vol(\cY)}{N^{l}|W|\vol(X)} \sum_{\lambda \in X/ \kappa a_{t-1} N X} h_{\lambda_{0},\lambda_{t+1}}(\lambda).
$$
Putting everything together we get
\begin{eqnarray*}
\mT^{\mC}_{\lambda_{0},\lambda_{t+1}} &=& K^{\mC_{t-1}}_{\lambda_{0}} \frac{{\I}^{\npr}}{\kappa^{l/2}N^{l}|W|}
      \left| \frac {\vol (\cY)}{\vol (X)} \right|^{1/2}
      \exp \left( -\frac{\pi \I}{\ch} m_{t} |\rho|^{2} \right) \\
 && \hspace{.3in} \times \exp \left( -\frac{\pi \I}{a_{t-1} c_{t-1}\kappa} |\lambda_{0}|^{2} \right)
      \sum_{\lambda \in X/ \kappa a_{t-1} N X} h_{\lambda_{0},\lambda_{t+1}}(\lambda).
\end{eqnarray*}
We are now in a position where we can use the reciprocity formula.
If we let $f=Na_{t}\id_{\frhR^{*}}$, $r=\kappa a_{t-1}N$, and
$\psi=\psi(w,w')=-\frac{1}{\kappa a_{t-1}}(w(\lambda_{0})+a_{t-1}w'(\lambda_{t+1}))$, $w,w' \in W$, then
\begin{eqnarray*}
&& \sum_{\lambda \in X/ \kappa a_{t-1} N X} h_{\lambda_{0},\lambda_{t+1}}(\lambda) = \sum_{w,w' \in W} \det(ww') \\
 &&\hspace{.3in} \times \sum_{\lambda \in X/ r X} \exp \left( \frac{\pi \I}{r} \la \lambda,f(\lambda)\ra \right)
\exp\left( 2\pi \I \la \lambda,\psi\ra\right).
\end{eqnarray*}
By (\ref{eq:integer1}), (\ref{eq:integer2}), and (\ref{eq:integer3}) 
the assumptions (\ref{eq:assumption}) are satisfied, so by
\refprop{prop:gauss} we obtain
{\allowdisplaybreaks
\begin{eqnarray*}
&& \sum_{\lambda \in X/ \kappa a_{t-1}N X} h_{\lambda_{0},\lambda_{t+1}}(\lambda) 
  = \vol(\cY)^{-1} \left( \frac{f}{\I} \right)^{-1/2} r^{l/2} \sum_{w,w' \in W} \det(ww') \\*
 && \hspace{.6in} \times \sum_{\mu \in \cY/f(\cY)} \exp \left( -\pi r \I \la \mu + \psi, f^{-1}(\mu + \psi) \ra \right) \\
 && \hspace{.2in} = \left( \frac{\I}{a_{t}} \right)^{l/2} \frac{(\kappa a_{t-1})^{l/2}}{\vol (\cY)} \sum_{w,w' \in W} \det(ww') \\*
 && \hspace{.6in} \times \sum_{\mu \in \cY / a_{t}N \cY}
      \exp\left( -\frac{\pi \kappa a_{t-1}\I }{a_{t}} | \mu +\psi|^{2} \right).
\end{eqnarray*}}\noindent
Here $\exp\left( -\frac{\pi \kappa a_{t-1}\I}{a_{t}}| \mu +\psi|^{2}\right)$
only depends on $\mu \pmod{a_{t}\cY}$, so in total we get that
\begin{eqnarray*}
\mT^{\mC}_{\lambda_{0},\lambda_{t+1}} &=& K^{\mC_{t-1}}_{\lambda_{0}} \frac{{\I}^{\npr}}{\kappa^{l/2}|W|}
      \exp \left( -\frac{\pi \I}{\ch} m_{t} |\rho|^{2} \right) 
      \exp \left( -\frac{\pi \I}{a_{t-1} c_{t-1}\kappa} |\lambda_{0}|^{2} \right) \\
 && \hspace{.1in} \times \left( \frac{\I}{a_{t}} \right)^{l/2} (\kappa a_{t-1})^{l/2} \sum_{w,w' \in W} \det(ww') \\
 && \hspace{.4in} \times \sum_{\mu \in \cY / a_{t} \cY}
      \exp\left( -\frac{\pi \kappa a_{t-1}\I }{a_{t}} | \mu +\psi |^{2} \right),
\end{eqnarray*}
where we use that $\vol (X) = \vol (\cY)^{-1}$.
In a similar way as in the first step of the induction we finally get
\begin{eqnarray*}
\mT^{\mC}_{\lambda_{0},\lambda_{t+1}} &=& K^{\mC_{t}}_{\lambda_{0}} \sum_{w \in W} \det(w) \\
 && \hspace{.6in} \times \sum_{\mu \in \cY / a_{t} \cY}
      \exp\left( -\frac{\pi c_{t}\I }{a_{t}\kappa} | \kappa \mu + \frac{w(\lambda_{0})}{c_{t}}
         + \lambda_{t+1})|^{2} \right),
\end{eqnarray*}
where we use that $a_{t-1}=c_{t}$. Here
{\allowdisplaybreaks
\begin{eqnarray*}
K^{\mC_{t}}_{\lambda_{0}} &=& K^{\mC_{t-1}}_{\lambda_{0}} \frac{{\I}^{\npr}}{\kappa^{l/2}}
      \exp \left( -\frac{\pi \I}{\ch} m_{t} |\rho|^{2} \right) 
      \exp \left( -\frac{\pi \I}{a_{t-1} c_{t-1}\kappa} |\lambda_{0}|^{2} \right) \\*
 && \hspace{.3in} \times \left( \frac{\I}{a_{t}} \right)^{l/2} (\kappa a_{t-1})^{l/2} \\
 &=& \frac {\I^{t \npr}}{(\kappa|a_{t-1}|)^{l/2}\vol (\cY)}\;
      \zeta^{l \; D_{t-1}} \frac{\I^{\npr}}{\kappa^{l/2}}
      \exp\left( -\frac{m_{t} \pi \I}{\ch} |\rho|^{2} \right) \\*
 && \hspace{.3in} \times \exp \left( -\frac{\pi\I}{\ch} (\sum_{i=1}^{t-1} m_{i}) |\rho |^{2} \right)
      \exp \left( -\frac{\pi\I}{\kappa}
         (\sum_{i=1}^{t-2} \frac{1}{a_{i-1}a_{i}}) |\lambda_{0} |^{2} \right) \\*
 && \hspace{.3in} \times \exp\left( -\frac{\pi \I}{a_{t-1} c_{t-1} \kappa} |\lambda_{0}|^{2} \right)
         \left( \frac{|a_{t-1}|}{|a_{t}|} \right)^{l/2}\kappa^{l/2} \zeta^{l\sign (a_{t} a_{t-1})} \\
&=& \frac{\I^{(t+1) \npr}}{(\kappa |a_{t}|)^{l/2} \vol (\cY)}\;\zeta^{l \; D_{t}}
      \exp \left( -\frac{\pi\I}{\ch} (\sum_{i=1}^{t} m_{i}) |\rho |^{2}  \right) \\*
 && \hspace{.3in} \times \exp \left( -\frac{\pi\I}{\kappa}
      (\sum_{i=1}^{t-1} \frac{1}{a_{i-1}a_{i}}) |\lambda_{0} |^{2} \right),
\end{eqnarray*}}\noindent
where we use that $c_{t-1}=a_{t-2}$.
\end{proof}

We want to use the above result to find a simple expression for the
entries of $\mR(U)$ in terms of the entries of $U$ and data for
the Lie algebra $\frg$. There is, however, a small hurdle to
overcome because of the assumption on the $a_{k}$'s in \reflem{lem:main}.
To this end we need

\begin{lem}\label{lem:matrixdecomposition}
Let $U=\begin{bmatrix}
a & b \\
c & d
\end{bmatrix} \in \PSL (2,\Z)=\SL (2,\Z)/\{\pm 1\}$ with $c \neq 0$. Then we can write
$$
U=V\Theta^{n}
$$
where $n \in \Z$ and $V$ is given in the following way:
If $a=0$ then $V =\Xi$; if $a\neq 0$ then
there exists a sequence of integers $\mC$ such that
$V=B^{\mC}$ and such that $a^{\mC}_{k} \neq 0$,
$k=1,2,\ldots,|\mC|$.
\end{lem}

\begin{proof}
If $a=0$ then $1=\det(U)=-bc$ so $c=-b=\pm 1$. As element of $\PSL (2,\Z)$
we therefore have
$$
U=\begin{bmatrix}
0 & -1 \\
1 & d'
\end{bmatrix}=\Xi\Theta^{d'}
$$
where $d' \in \{ \pm d \}$. Now assume that $a \neq 0$.
By (\ref{eq:relationsXiTheta}) we can find a tuple of
integers $\mC'=(m_{1},\ldots,m_{t})$ such that $U=B^{\mC'}$.
If $a_{i}=a^{\mC'}_{i} \neq 0$, $i=1,2,\ldots,t$, we let
$\mC=\mC'$ and $n=0$. Therefore assume this is not the case. Let
$$
i = \max \{ \; j \in \{1,2,\ldots,t\} \; | \; a_{j}=0 \}.
$$
Since $a_{t} = a \neq 0$ we have $i < t$. As in the case $a=0$
we have that $B^{\mC'_{i}}=\Xi \Theta^{j}$ for some $j \in \Z$.
If $i=t-1$ then
$$
U=\Theta^{m_{t}}\Xi^{2}\Theta^{j}=\Theta^{m_{t}+j}=\begin{bmatrix}
1 & m_{t}+j \\
0 & 1
\end{bmatrix}
$$
contradicting the fact that $c \neq 0$. Therefore $i < t-1$ and
$U=W\Theta^{n'}$ with $n'=m_{i+1}+j$ and
$$
W=\Theta^{m_{t}}\Xi\cdots\Theta^{m_{i+2}}\Xi.
$$
Let $n_{k}=m_{i+k+1}$, $k=1,2,\ldots,t-i-1$, let $\mC''=(n_{1},\ldots,n_{t-i-1})$,
and let $a_{k}'=a_{k}^{\mC''}$ etc.
Then $W=B^{\mC''}$, and
$$
\begin{bmatrix}
a_{i+k+1} & b_{i+k+1} \\
c_{i+k+1} & d_{i+k+1}
\end{bmatrix}=B^{\mC'}_{i+k+1}=B^{\mC''}_{k}\Theta^{n'}=\begin{bmatrix}
a_{k}' & b_{k}' \\
c_{k}' & d_{k}'
\end{bmatrix}\begin{bmatrix}
1 & n' \\
0 & 1
\end{bmatrix}=
\begin{bmatrix}
a_{k}' & n'a_{k}'+b_{k}' \\
c_{k}' & n'c_{k}'+d_{k}'
\end{bmatrix},
$$
so in particular $a_{k}'=a_{i+k+1} \neq 0$, $k=1,2,\ldots,t-i-1$,
by the maximality of $i$. Therefore we can let $n=n'$ and $\mC=\mC''$.
\end{proof}

For the next theorem we need the Rademacher Phi function $\Phi$, which
is defined on $\PSL (2,\Z)$ by
\begin{equation}\label{eq:Phi}
\Phi \left( \left[ \begin{array}{cc}
                  a & b \\
                  c & d
                  \end{array}
\right] \right) = \left\{ \begin{array}{ll}
                  \frac{a+d}{c} - 12(\sign (c))\s (d,c) & ,c \neq 0, \\
                  \frac{b}{d} & ,c=0.
                  \end{array}
\right.
\end{equation}
Here, for $c \neq 0$, the Dedekind sum $\s (d,c)$ is given by
\begin{equation}\label{eq:dedekindsum}
\s (d,c)= \frac{1}{4|c|} \sum_{j=1}^{|c|-1} \cot\frac{\pi j}{c} \cot \frac{\pi d j}{c}
\end{equation}
for $|c|>1$ and $\s (d,\pm 1)=0$. If
$A_{i} = \left( \begin{array}{cc}
                a_{i} & b_{i} \\
		c_{i} & d_{i}
		\end{array}
\right) \in \SL (2,\Z)$
such that $A_{3}=A_{1}A_{2}$ we have
\begin{equation}\label{eq:Phiproduct}
\Phi(A_{3})=\Phi(A_{1})+\Phi(A_{2})-3\sign (c_{1}c_{2}c_{3}).
\end{equation}
We refer to \cite{RademacherGrosswald} for a comprehensive description of
the Rademacher Phi function. We will also need
\begin{equation}\label{eq:Freudenthal}
\frac{|\rho|^{2}}{\ch} = \frac{\dim\frg}{12} = \frac{2 \npr +l}{12},
\end{equation}
where the first identity is Freudenthal's strange formula.
If $c=0$ then $U= \ep \Theta^{b}$ for some $b \in \Z$ and $\ep \in \{ \pm 1 \}$.
By (\ref{eq:mR}) we immediately get
$$
\mR(\Theta^{b})_{\lambda \mu} = \mR(\Theta)_{\lambda \mu}^{b}=\delta_{\lambda \mu}
             \exp \left( b \left(
                  \frac {\pi\sqrt{-1}}{\kappa} |\mu|^{2}
                 -\frac {\pi\sqrt{-1}}{\ch} |\rho |^{2}
                  \right)\right)
$$
for $\lambda,\mu \in \inte(C_{\kappa}) \cap X$.
For the case $U=-\Theta^{b}=\Xi^{2}\Theta^{b}$ (see (\ref{eq:relationsXiTheta}))
we use the identity $\mR(\Xi^{2})_{\lambda \mu}=\delta_{\lambda \mu^{*}}$,
which follows by the unitarity of $\mR$ and (\ref{eq:unitarity3}) (alternatively
use (\ref{eq:unitarity1}) and (\ref{eq:unitarity2}) directly).
This gives 
$$
\mR(-\Theta^{b})_{\lambda \mu} =  \mR(\Theta)_{\lambda^{*}\mu}^{b} = \delta_{\lambda^{*} \mu}
             \exp \left( b \left(
                  \frac {\pi\sqrt{-1}}{\kappa} |\mu|^{2}
                 -\frac {\pi\sqrt{-1}}{\ch} |\rho|^{2}
                  \right)\right)
$$
for $\lambda,\mu \in \inte(C_{\kappa}) \cap X$.
For the case $c \neq 0$ we have

\begin{thm}\label{thm:main1}
Let $U=\left( \begin{array}{cc}
		a & b \\
		c & d
		\end{array}
\right) \in \SL (2,\Z)$ with $c \neq 0$, and let $\lambda,\mu \in \inte(C_{\kappa}) \cap X$.
Then
\begin{eqnarray*}
\mR(U)_{\lambda \mu} &=& \frac{ (\I\sign(c))^{\npr}}{(\kappa|c|)^{l/2} \vol (\cY)}
      \exp \left( -\frac{\pi \I}{\ch} \Phi(U)|\rho|^{2} \right) \\
 && \hspace{.3in} \times \exp \left( \frac{\pi \I}{\kappa} \frac{d}{c} |\mu|^{2} \right) \sum_{\nu \in \cY/c \cY}
      \exp \left( \frac{\pi \I} {\kappa} \frac{a}{c} |\lambda+\kappa\nu|^{2} \right) \\
 && \hspace{.3in} \times \sum_{w\in W} \det (w)
      \exp\left( -\frac{2\pi\I}{\kappa c} \la \lambda + \kappa\nu,w(\mu) \ra \right).
\end{eqnarray*}
If $a \neq 0$ we also have
\begin{eqnarray*}
\mR(U)_{\lambda \mu} &=& \frac{ (\I\sign(c))^{\npr}}{(\kappa|c|)^{l/2} \vol (\cY)}
      \exp \left( -\frac{\pi \I}{\ch} \Phi(U) |\rho|^{2} \right) \\
 && \hspace{.3in} \times \exp \left( \frac{\pi \I}{\kappa} \frac{b}{a} |\mu|^{2} \right) \sum_{w\in W} \det (w) \\
 && \hspace{.3in} \times \sum_{\nu \in \cY/c \cY}
      \exp\left( \frac{\pi\I}{\kappa} \frac{a}{c}
         |\lambda +\kappa\nu-\frac {w(\mu)}{a}|^{2} \right).
\end{eqnarray*}
\end{thm}

\begin{proof}
According to the previous lemma there exists an integer $n$, a sign $\ep \in \{ \pm 1 \}$
and a $V \in \SL (2,\Z)$ as in the \reflem{lem:matrixdecomposition}
such that $U=\ep V\Theta^{n}$. Let us first assume that $\ep=1$.
Assume, moreover, that $a \neq 0$ and that $n=0$, i.e.\ assume that $U=B^{\mC}$,
 where $\mC=(m_{1},\ldots,m_{t})$
and $a_{k} \neq 0$, $k=1,2,\ldots,t$. Let $\mC'=(m_{1},\ldots,m_{t-1})$.
Then by \reflem{lem:main}
{\allowdisplaybreaks
\begin{eqnarray*}
\mR(U)_{\lambda \mu} &=& \tilde{\Theta}_{\lambda\lambda}^{m_{t}} \mT^{\mC'}_{\mu,\lambda} \\*
&=& K^{\mC'}_{\mu} \exp \left( \frac {\pi \I}{\kappa} m_{t} |\lambda|^{2}
       -\frac{\pi \I}{\ch} m_{t}|\rho|^{2} \right) \sum_{w\in W} \det (w) \\*
 && \hspace{.3in} \times \sum_{\nu \in \cY/a_{t-1} \cY}
      \exp\left( -\frac{\pi\I c_{t-1}}{a_{t-1} \kappa}
         \left|\lambda+\kappa\nu+\frac{w(\mu)}{c_{t-1}}\right|^{2} \right) \\
&=& K^{\mC'}_{\mu} \exp \left( -\frac{\pi \I}{\ch} m_{t} |\rho|^{2} \right)
      \exp \left( -\frac{\pi \I}{a_{t-1} c_{t-1}\kappa} |\mu|^{2} \right) \\*
 && \hspace{.3in} \times \sum_{w\in W} \det (w)
      \sum_{\nu \in \cY/a_{t-1} \cY}
      \exp\left( \frac{\pi\I}{\kappa} m_{t} \left|\lambda+\kappa\nu\right|^{2} \right) \\*
 && \hspace{.3in} \times
      \exp\left( -\frac{\pi\I c_{t-1}}{a_{t-1}\kappa}
         \left|\lambda+\kappa\nu\right|^{2} \right) \\*
 && \hspace{.3in} \times \exp\left( -\frac {2\pi\I }{ a_{t-1}\kappa}
      \la \lambda+\kappa\nu, w(\mu) \ra \right) \\
&=&  K^{\mC'}_{\mu} \exp \left( -\frac{\pi \I}{\ch} m_{t} |\rho|^{2} \right)
      \exp \left( -\frac{\pi \I}{a_{t-1} a_{t-2}\kappa} |\mu|^{2} \right) \\*
 && \hspace{.3in} \times \exp \left( -\frac{\pi \I}{a_{t} a_{t-1}\kappa} |\mu|^{2} \right) \sum_{w\in W} \det (w) \\*
 && \hspace{.3in} \times \sum_{\nu \in \cY/c \cY}
      \exp\left( \frac{\pi\I}{\kappa} \frac{a}{c}
         \left|\lambda+\kappa\nu-\frac {w(\mu)}{a}\right|^{2} \right),
\end{eqnarray*}}\noindent
where we have used (\ref{eq:abcd}). Let us calculate the factor
in front of the sum.
By \reflem{lem:main} we get
{\allowdisplaybreaks
\begin{eqnarray*}
K &:=& K^{\mC'}_{\mu} \exp \left( -\frac{\pi \I}{\ch} m_{t} |\rho|^{2} \right) \\*
 && \hspace{.5in} \times \exp \left( -\frac{\pi \I}{a_{t-1}a_{t-2}\kappa} |\mu|^{2} \right) 
       \exp \left( -\frac{\pi \I}{a_{t}a_{t-1}\kappa} |\mu|^{2} \right) \\
&=& \frac{\I^{t \npr} \zeta^{l D_{t-1}}}{(\kappa|c|)^{l/2}\vol (\cY)}
       \exp \left( -\frac{\pi\I}{\ch}(\sum_{i=1}^{t} m_{i}) |\rho |^{2} \right) \\*
 && \hspace{.5in} \times \exp \left( -\frac{\pi\I}{\kappa}(\sum_{i=1}^{t} \frac{1}{a_{i-1}a_{i}})
          |\mu|^2 \right).
\end{eqnarray*}}\noindent
By \cite[Formula (2.20)]{Jeffrey2} we have
$$
\Phi(U)=\sum_{i=1}^{t} m_{i} - 3 \sum_{i=1}^{t-1} \sign(a_{i-1}a_{i}).
$$
This together with (\ref{eq:b/a}) gives
\begin{eqnarray*}
K &=& \frac{\I^{t \npr} \zeta^{l D_{t-1}}}{(\kappa|c|)^{l/2}\vol (\cY)}
         \exp \left( -\frac{\pi\I}{\ch}( \Phi(U)+ 3D_{t-1}) |\rho |^{2} \right) \\
 && \hspace{.5in} \times \exp \left( -\frac{\pi\I}{\kappa} \frac{b}{a}|\mu|^{2} \right).
\end{eqnarray*}
By (\ref{eq:Freudenthal}) we then have
{\allowdisplaybreaks
\begin{eqnarray*}
&& \I^{(t-1) \npr} \zeta^{l D_{t-1}}
      \exp\left( -\frac{\pi\I}{\ch}(\Phi(U)+ 3D_{t-1}) |\rho |^{2} \right) \\
 && \hspace{.3in} = \I^{(t-1) \npr}
      \exp\left( \left( \frac{l\pi\I}{4} - \frac{(2\npr +l)\pi\I}{4} \right) D_{t-1} \right) \\*
 && \hspace{1.0in} \times \exp\left( -\frac{\pi\I}{\ch} \Phi(U) |\rho |^{2} \right) \\
 && \hspace{.3in} = \I^{(t-1-D_{t-1})\npr}
      \exp\left( -\frac{\pi\I}{\ch} \Phi(U)|\rho |^{2}  \right).
\end{eqnarray*}}\noindent
Moreover,
$$
t-1-D_{t-1}=t-1-\sum_{i=1}^{t-1} \sign (a_{i-1}a_{i})
     =\sum_{i=1}^{t-1} (1-\sign (a_{i-1}a_{i}))
$$
and
$$
\I^{1-\sign (a_{i-1}a_{i})}=\sign (a_{i-1}a_{i}),
$$
so $\I^{t-1-D_{t-1}}=\sign(a_{t-1})=\sign(c)$ since $a_{0}=1$.
This proves the last given formula for the entries of $\mR(U)$.
The first formula follows easily from this by observing that
$$
\frac{b}{a} + \frac{1}{ac} = \frac{d}{c}.
$$
Next assume that $U = B^{\mC}\Theta^{n}$ with $n \neq 0$,
where $\mC$ is as above. Then
$$
B^{\mC} = U \Theta^{-n} = \left( \begin{array}{cc}
		a & b \\
		c & d
		\end{array}
\right)\left( \begin{array}{cc}
		1 & -n \\
		0 & 1
		\end{array}
\right)=\left( \begin{array}{cc}
		a & -na+b \\
		c & -nc+d
		\end{array}
\right),
$$
and since the theorem is valid for $U=B^{\mC}$ (by the above) we get
{\allowdisplaybreaks
\begin{eqnarray*}
\mR(U)_{\lambda \mu} &=& \mR(B^{\mC})_{\lambda \mu}\mR(\Theta)^{n}_{\mu \mu} \\*
&=& \frac{ (\I\sign(c))^{\npr}}{(\kappa|c|)^{l/2} \vol (\cY)}
      \exp\left( -\frac{\pi \I}{\ch} \left( \Phi(B^{\mC}) +n\right)|\rho|^{2} \right) \\*
 && \hspace{.3in} \times \exp\left( \frac{\pi \I}{\kappa} \frac{d}{c} |\mu|^{2} \right) \sum_{\nu \in \cY/c \cY}
      \exp\left( \frac{\pi \I}{\kappa} \frac{a}{c} |\lambda+\kappa\nu|^{2} \right) \\*
 && \hspace{.3in} \times \sum_{w\in W} \det (w)
      \exp\left( -\frac{2\pi\I}{\kappa c} \la \lambda +\kappa\nu,w(\mu) \ra \right).
\end{eqnarray*}}\noindent
By (\ref{eq:Phiproduct}), $\Phi(B^{\mC}\Theta^{n})=\Phi(B^{\mC})+\Phi(\Theta^{n})=\Phi(B^{\mC})+n$
and the result follows.
Next consider the case where $a=0$ (and $\ep=1$) so
$U=\Xi \Theta^{n}=
\left( \begin{array}{cc}
		0 & -1 \\
		1 & n
		\end{array}
\right)$. Then the result follows directly by inserting the formulas
for the entries of $\mR(\Xi)$ and $\mR(\Theta)$ into
$$
\mR(U)_{\lambda \mu}=\mR(\Xi)_{\lambda \mu}\mR(\Theta)^{n}_{\mu \mu}
$$
and by using that $\Phi(\Xi\Theta^{n})=\Phi(\Xi)+\Phi(\Theta^{n})=n$
(use (\ref{eq:Phiproduct})).

Let us finally consider the case where $\ep =-1$, so $U=-V\Theta^{n}$ (V being as in
\reflem{lem:matrixdecomposition}). By the remarks just before the
theorem we get
$$
\mR(U)_{\lambda\mu} = \mR(-U)_{\lambda^{*}\mu}.
$$
Since the theorem is valid for $-U$ (by the above) and since
$\lambda^{*}=-w_{0}(\lambda)$ and $\det(w_{0})=(-1)^{\npr}$, see \refrem{rem:unitarity}, we get
{\allowdisplaybreaks
\begin{eqnarray*}
\mR(U)_{\lambda\mu} &=& \frac{ (\I\sign(c))^{\npr}}{(\kappa|c|)^{l/2} \vol (\cY)}
      \exp \left( -\frac{\pi \I}{\ch} \Phi(U) |\rho|^{2} \right) 
      \exp \left( \frac{\pi \I}{\kappa} \frac{d}{c} |\mu|^{2} \right) \\*
 && \hspace{.3in} \times \sum_{\nu \in \cY/c \cY}
      \exp \left( \frac{\pi \I} {\kappa} \frac{a}{c} |-w_{0}(\lambda)+\kappa w_{0}(\nu)|^{2} \right) \\*
 && \hspace{.3in} \times \sum_{w\in W} \det (w_{0}^{-1}w)
      \exp\left( \frac{2\pi\I}{\kappa c} \la -w_{0}(\lambda) + \kappa w_{0}(\nu),w(\mu) \ra \right) \\
 &=& \frac{ (\I\sign(c))^{\npr}}{(\kappa|c|)^{l/2} \vol (\cY)}
      \exp \left( -\frac{\pi \I}{\ch} \Phi(U) |\rho|^{2}\right) \\*
 && \hspace{.3in} \times \exp \left( \frac{\pi \I}{\kappa} \frac{d}{c} |\mu|^{2} \right) \sum_{\nu \in \cY/c \cY}
      \exp \left( \frac{\pi \I} {\kappa} \frac{a}{c} |\lambda+\kappa \nu|^{2} \right) \\*
 && \hspace{.3in} \times \sum_{w\in W} \det (w_{0}^{-1}w)
      \exp\left( -\frac{2\pi\I}{\kappa c} \la \lambda + \kappa \nu,w_{0}^{-1}w(\mu) \ra \right)
\end{eqnarray*}}\noindent
by which the theorem follows. 
\end{proof}

Since $\mR$ is a unitary representation we have
\begin{equation}\label{eq:Uunitarity}
\mR(U)_{\lambda \mu}=\overline{\mR(U^{-1})_{\mu \lambda}}
\end{equation}
for any $U \in \SL (2,\Z)$ and all $\lambda,\mu \in \inte(C_{\kappa}) \cap X$.
By this and the facts that
$\Phi(U^{-1})=-\Phi(U)$ (use (\ref{eq:Phiproduct})) and
$U^{-1}=\left( \begin{array}{cc}
		d & -b \\
		-c & a
		\end{array}
\right)$
we get

\begin{cor}\label{cor:unitarity}
Let $U=\left( \begin{array}{cc}
a & b \\
c & d
\end{array}\right) \in \SL (2,\Z)$ with $c \neq 0$, and let $\lambda,\mu \in \inte(C_{\kappa}) \cap X$.
Then
\begin{eqnarray*}
\mR(U)_{\lambda \mu} &=& \frac{ \left(\I\sign(c)\right)^{\npr}}{(\kappa|c|)^{l/2} \vol (\cY)}
      \exp\left( -\frac{\pi \I}{\ch} \Phi(U)|\rho|^{2} \right) \\
 && \hspace{.3in} \times \exp\left( \frac{\pi \I}{\kappa} \frac{a}{c} |\lambda|^{2} \right) \sum_{\nu \in \cY/c \cY}
      \exp\left( \frac{\pi \I}{\kappa} \frac{d}{c} |\mu+\kappa\nu|^{2} \right) \\
 && \hspace{.3in} \times \sum_{w\in W} \det (w)
      \exp\left( -\frac{2\pi\I}{\kappa c} \la \mu +\kappa\nu,w(\lambda) \ra \right).
\end{eqnarray*}
If $d \neq 0$ we also have
\begin{eqnarray*}
\mR(U)_{\lambda \mu} &=& \frac{ \left(\I\sign(c)\right)^{\npr}}{(\kappa|c|)^{l/2} \vol (\cY)}
      \exp\left( -\frac{\pi \I}{\ch} \Phi(U) |\rho|^{2}\right) \\
 && \hspace{.3in} \times \exp\left( \frac{\pi \I}{\kappa} \frac{b}{d} |\lambda|^{2} \right) \sum_{w\in W} \det (w) \\
 && \hspace{.3in} \times \sum_{\nu \in \cY/c \cY}
      \exp\left( \frac{\pi\I}{\kappa} \frac{d}{c}
      |\mu +\kappa\nu - \frac{w(\lambda)}{d}|^{2} \right).
\end{eqnarray*}\HS
\end{cor}

We are particularly interested in expressions for
$\mR(U)_{\lambda \mu}$ in case $\lambda$
or $\mu$ is equal to $\rho$. Note that
since $\rho^{*}=\rho$, $\Xi^{2}=-1$ and
$\mR(\Xi^{2})_{\lambda \mu}=\delta_{\lambda \mu^{*}}$,
then
\begin{equation}\label{eq:signindependence}
\mR(-U)_{\lambda \rho} = \mR(U)_{\lambda \rho},\hspace{.2in}\mR(-U)_{\rho \lambda} = \mR(U)_{\rho \lambda}
\end{equation}
for all $\lambda \in \inte(C_{\kappa}) \cap X$,
so these entries are in fact functions of $\PSL (2,\Z)$. 
By the Weyl denominator formula we have
\begin{eqnarray*}
&& \sum_{w\in W} \det (w) \exp\left\{ -\frac{2\pi\I}{\kappa c} \la \lambda + \kappa\nu,w(\rho) \ra \right\} \\
 && \hspace{.5in} = \prod_{\alpha \in \Delta_{+}} 2\I \sin\left( -\frac{\pi}{\kappa c} \la \lambda +\kappa\nu,\alpha \ra \right),
\end{eqnarray*}
so we get the alternative formula
\begin{eqnarray}\label{eq:rho1}
\mR(U)_{\lambda \rho} &=& \frac{ \left(2\sign(c)\right)^{\npr}}{(\kappa |c|)^{l/2} \vol (\cY)}
      \exp \left( -\frac{\pi \I}{\ch} \Phi(U)|\rho|^{2} \right) \\
 && \hspace{.1in} \times \exp \left( \frac{\pi \I}{\kappa} \frac{d}{c} |\rho|^{2} \right) \sum_{\nu \in \cY/c \cY}
      \exp \left( \frac{\pi \I}{\kappa} \frac{a}{c} |\lambda+\kappa\nu|^{2} \right) \nonumber \\
 && \hspace{.3in} \times \prod_{\alpha \in \Delta_{+}}
      \sin\left( \frac{\pi}{\kappa c} \la \lambda +\kappa\nu,\alpha \ra \right).\nonumber
\end{eqnarray}
for $U \in \SL (2,\Z)$ as in \refthm{thm:main1}.
By using the Weyl denominator formula and the first expression
in \refcor{cor:unitarity} we get
\begin{eqnarray}\label{eq:rho2}
\mR(U)_{\rho \mu} &=& \frac{ \left(2\sign(c)\right)^{\npr}}{(\kappa |c|)^{l/2} \vol (\cY)}
      \exp \left( -\frac{\pi \I}{\ch} \Phi(U)|\rho|^{2}  \right) \\
 && \hspace{.1in} \times \exp \left( \frac{\pi \I}{\kappa} \frac{a}{c} |\rho|^{2} \right) \sum_{\nu \in \cY/c \cY}
      \exp \left( \frac{\pi \I}{\kappa} \frac{d}{c} |\mu +\kappa\nu|^{2} \right) \nonumber \\
 && \hspace{.3in} \times \prod_{\alpha \in \Delta_{+}}
      \sin\left( \frac{\pi}{\kappa c} \la \mu +\kappa\nu,\alpha \ra \right).\nonumber
\end{eqnarray}

We end this section with some symmetry considerations needed elsewhere.
As mentioned allready below \refrem{rem:unitarity} the expressions for
$\mR(\Xi)_{\lambda\mu}$ and $\mR(\Theta)_{\lambda\mu}$ in (\ref{eq:mR})
are well-defined for all $\lambda,\mu \in X$.
If $U \in \SL (2,\Z)$ we can find a tuple of integers $\mC=(m_{1},\ldots,m_{t})$ such that
$U=B^{\mC}$. We can therefore use the formula
$$
\mR(U)_{\lambda\mu}=\mR(\Theta)_{\lambda\lambda}^{m_{t}}\mT_{\mu,\lambda}^{\mC_{t-1}}
$$
to extend $(\lambda,\mu) \mapsto \mR(U)_{\lambda\mu}$ to all of $X \times X$.
Here $\mT_{\mu,\lambda}^{\mC_{t-1}}$ is defined above \reflem{lem:main}.
We can also make such an extension directly by using the expressions
in \refthm{thm:main1}. An easy inspection of the proof of \refthm{thm:main1}
shows that these two extensions coincide.

\begin{lem}\label{lem:mRsymmetry}
Let $U \in \SL (2,\Z)$. Then $(\lambda,\mu) \mapsto \mR(U)_{\lambda\mu}$,
$X \times X \to \C$, is invariant under the action of $\kappa \cY$
on each factor. Moreover, 
$$
\mR(U)_{w(\lambda)w'(\mu)}=\det(w)\det(w')\mR(U)_{\lambda\mu}.
$$
Finally, $\mR(U)_{\lambda\mu}=0$ if $\lambda$ or $\mu$ belongs to
$X\cap H^{\kappa}$.
\end{lem}

\begin{proof}
Since $\Xi^{4}=1$ we have $U=\Xi B = C \Xi$, where $B=\Xi^{3}U$ and
$C=U\Xi^{3}$. The lemma then follows by the fact that it is true
for $U=\Xi$.
\end{proof}

Since $C_{\kappa}$ is a fundamental domain for the action of $\aW_{\kappa}$
it follows by \reflem{lem:mRsymmetry} that (\ref{eq:Uunitarity}) is
valid for all $\lambda,\mu \in X$. Therefore we have that the
expressions for $\mR(U)_{\lambda\mu}$ stated in \refcor{cor:unitarity}
are valid for all $\lambda,\mu \in X$. Note that (\ref{eq:rho1})
is valid for all $\lambda \in X$ and (\ref{eq:rho2}) is valid for all
$\mu \in X$.

\section{Seifert manifolds}\label{sec-Seifert-manifolds}

For Seifert manifolds we will use the notation introduced by Seifert in his
classification results for these manifolds, see \cite{Seifert1},
\cite{Seifert2}, \cite[Sect.~2]{Hansen2}. That is, 
$(\ep;g\,|\,b;(\alpha_{1},\beta_{1}),\ldots,(\alpha_{n},\beta_{n}))$ is
the Seifert manifold with orientable base of genus $g \geq 0$ if
$\ep=\os$ and non-orientable base of genus $g>0$ if $\ep=\ns$ 
(where the genus of the non-orientable connected sum $\# ^{k} \R \text{P}^{2}$
is $k$). (In \cite{Seifert1}, \cite{Seifert2}
$(\ep;g\,|\,b;(\alpha_{1},\beta_{1}),\ldots,(\alpha_{n},\beta_{n}))$ is
denoted $(\tO,\ep;g\,|\,b;\alpha_{1},\beta_{1};\ldots;\alpha_{n},\beta_{n})$,
but we leave out the $\tO$ since we are only dealing with
oriented Seifert manifolds.)
The pair $(\alpha_{j},\beta_{j})$ of coprime integers is the (oriented)
Seifert invariant of the $j$'th exceptional (or singular) fiber.
We have $0< \beta_{j}<\alpha_{j}$. The integer $-b$ is equal to the
Euler number of the Seifert fibration $(\ep;g\,|\,b)$
(which is a locally trivial $S^{1}$--bundle).
The sign is chosen so that the Euler number of the spherical (or unit) tangent
bundle over an orientable surface $\Sigma$ is equal to the Euler characteristic of
$\Sigma$, see \cite[Chap.~1 and 4]{Montesinos}, \cite[Sect.~3]{Scott}.
More generally the Seifert Euler number of
$(\ep;g\,|\,b;(\alpha_{1},\beta_{1}),\ldots,(\alpha_{n},\beta_{n}))$
is $E=-\left(b+\sum_{j=1}^{n} \beta_{j}/\alpha_{j} \right)$.
We note that lens spaces are Seifert manifolds with base $S^{2}$
and zero, one or two exceptional fibers.
According to \cite[Fig.~12 p.~146]{Montesinos}, the manifold
$( \ep ;g \, | \, b;(\alpha_{1},\beta_{1}),\ldots,(\alpha_{n},\beta_{n}))$
has a surgery
presentation as shown in Fig.~\ref{fig-A1} if $\ep=\os$ and as
shown in Fig.~\ref{fig-A2} if $\ep=\ns$.
The $\undersmile{g}$ indicate $g$ repetitions.

\begin{figure}[h]

\begin{center}
\begin{texdraw}
\drawdim{cm}

\setunitscale 0.6

\linewd 0.02 \setgray 0

\move(0 4)

\move(0 0) \lellip rx:3 ry:1.6

\linewd 0.2 \setgray 1

\move(4 0) \larc r:1.5 sd:218 ed:200
\move(4 0) \larc r:3 sd:213 ed:204
\move(-4 0) \larc r:1.5 sd:-20 ed:150
\move(-4 0) \larc r:2 sd:-23 ed:150
\move(-4 0) \larc r:3.5 sd:-20 ed:155
\move(-4 0) \larc r:3 sd:-23 ed:155

\linewd 0.02 \setgray 0

\move(-4 0) \larc r:1.5 sd:-20 ed:150
\move(-4 0) \larc r:2 sd:-23 ed:150
\move(-4 0) \larc r:3.5 sd:-20 ed:155
\move(-4 0) \larc r:3 sd:-23 ed:155

\move(-5.88 0.67) \clvec(-5.85 0.9)(-5.6 0.95)(-5.45 0.8)
\move(-5.88 -0.67) \clvec(-5.85 -0.9)(-5.6 -0.95)(-5.45 -0.8)

\move(-5.3 0.75) \clvec(-5.35 0.55)(-5.6 0.5)(-5.75 0.7)
\move(-5.3 -0.75) \clvec(-5.35 -0.55)(-5.6 -0.5)(-5.75 -0.7)

\move(-7.3 1.17) \clvec(-7.25 1.4)(-7 1.45)(-6.85 1.3)
\move(-7.3 -1.17) \clvec(-7.25 -1.4)(-7 -1.45)(-6.85 -1.3)

\move(-6.7 1.29) \clvec(-6.75 1.09)(-7 1.04)(-7.15 1.24)
\move(-6.7 -1.29) \clvec(-6.75 -1.09)(-7 -1.04)(-7.15 -1.24)

\move(-4 0) \larc r:1.5 sd:210 ed:320
\move(-4 0) \larc r:2 sd:210 ed:318
\move(-4 0) \larc r:3.5 sd:205 ed:330
\move(-4 0) \larc r:3 sd:205 ed:328
\move(-4 0) \larc r:1.5 sd:160 ed:200
\move(-4 0) \larc r:2 sd:160 ed:200
\move(-4 0) \larc r:3.5 sd:160 ed:200
\move(-4 0) \larc r:3 sd:160 ed:200

\move(4 0) \larc r:1.5 sd:218 ed:200
\move(4 0) \larc r:3 sd:213 ed:204

\move(-6.8 -0.1) \htext{$\cdots$}
\move(1.5 -0.1) \htext{$\cdots$}
\move(-6.85 -0.8) \htext{$\undersmile{g}$}
\move(-6.65 2.8) \htext{$0$}
\move(-5.75 1.6) \htext{$0$}
\move(-7.95 0.3) \htext{$0$}
\move(-6.45 0.15) \htext{$0$}
\move(-3.95 -0.3) \htext{$-b$}
\move(5.1 1) \htext{$\frac{\alpha_{1}}{\beta_{1}}$}
\move(6.2 2.2) \htext{$\frac{\alpha_{n}}{\beta_{n}}$}

\move(0 -4)

\end{texdraw}
\end{center}

\caption{Surgery presentation of $(\os ;g \, | \, b;(\alpha_{1},\beta_{1}),\ldots,(\alpha_{n},\beta_{n}))$}\label{fig-A1}
\end{figure}

\begin{figure}[h]

\begin{center}
\begin{texdraw}
\drawdim{cm}

\setunitscale 0.6

\linewd 0.02 \setgray 0

\move(0 4)

\move(0 0) \lellip rx:3 ry:1.6
\move(-7.5 0) \lellip rx:0.5 ry:0.3
\move(-5.5 0) \lellip rx:0.5 ry:0.3

\linewd 0.2 \setgray 1
\move(-4 0) \larc r:1.5 sd:-20 ed:180
\move(-4 0) \larc r:3.5 sd:-20 ed:180
\move(4 0) \larc r:1.5 sd:215 ed:200
\move(4 0) \larc r:3 sd:212 ed:205

\linewd 0.02 \setgray 0

\move(-4 0) \larc r:1.5 sd:-20 ed:180
\move(-4 0) \larc r:3.5 sd:-20 ed:180

\move(-4 0) \larc r:1.5 sd:200 ed:320
\move(-4 0) \larc r:3.5 sd:188 ed:330

\move(4 0) \larc r:1.5 sd:215 ed:200
\move(4 0) \larc r:3 sd:212 ed:205

\move(-6.8 -0.1) \htext{$\cdots$}
\move(1.5 -0.1) \htext{$\cdots$}
\move(-6.8 -0.9) \htext{$\undersmile{g}$}
\move(-8.3 0.25) \htext{$2$}
\move(-6.3 0.25) \htext{$2$}
\move(-7.4 2.2) \htext{$\frac{1}{2}$}
\move(-5.75 1) \htext{$\frac{1}{2}$}
\move(-3.95 -0.3) \htext{$-b$}
\move(5.1 1) \htext{$\frac{\alpha_{1}}{\beta_{1}}$}
\move(6.2 2.2) \htext{$\frac{\alpha_{n}}{\beta_{n}}$}

\move(0 -4)

\end{texdraw}
\end{center}

\caption{Surgery presentation of $(\ns ; g \, | \, b;(\alpha_{1},\beta_{1}),\ldots,(\alpha_{n},\beta_{n}))$}\label{fig-A2}
\end{figure}

For completeness we will also state
the results in terms of the non-normalized
Seifert invariants due to W.\ D.\ Neumann, see \cite{JankinsNeumann}.
For a Seifert manifold $X$ with non-normalized Seifert invariants
$\{\ep;g;(\alpha_{1},\beta_{1}),\ldots,(\alpha_{n},\beta_{n})\}$
the invariants $\ep$ and $g$ are as above. The
$(\alpha_{j},\beta_{j})$ are here pairs of coprime integers with
$\alpha_{j} >0$ but not necessarily with $0< \beta_{j}<\alpha_{j}$.
These pairs are not invariants of $X$, but can be varied
according to certain rules. 
In fact, $X$ has a surgery presentation as shown in
Fig.~\ref{fig-A1} with $b=0$ if $\ep=\os$ and
as shown in Fig.~\ref{fig-A2} with $b=0$ if $\ep=\ns$.
The Seifert Euler number of $X$ is $-\sum_{j=1}^{n} \beta_{j}/\alpha_{j}$
(which is an invariant of the Seifert fibration $X$).
For more details we refer to 
\cite[Sect.~I.1]{JankinsNeumann}.

\section{The RT--invariants of the Seifert $3$--manifolds}\label{sec-The-RT--invariants}

In this section we will use \refcor{cor:unitarity} together with
results in \cite{Hansen2} to derive expressions for
the quantum $\frg$--invariants of all Seifert manifolds,
$\frg$ being an arbitrary complex finite dimensional simple Lie algebra.
The formulas for the
invariants will be expressed in terms of the Seifert invariants
together with standard data for $\frg$.

\rk{The RT--invariants of the Seifert manifolds for modular categories}
Let us first give some preliminary remarks on modular categories. We use
notation as in \cite{Turaev}.
Let $\left( \mV, \{ V_{i} \}_{i \in I } \right)$
be an arbitrary modular category with braiding $c$ and
twist $\theta$. The ground ring is $K=\End_{\mV}(\Io)$,
where $\Io$ is the unit object.
Let $i \mapsto i^{*}$ be the involution in $I$
determined by the condition that $V_{i^{*}}$ is isomorphic to
the dual of $V_{i}$.
An element $i \in I$ is called self-dual if $i=i^{*}$. For such an element we
have a $K$--module isomorphism $\Hom_{\mV}(V \otimes V, \Io ) \cong K$, $V=V_{i}$.
The map $x \mapsto x(\id_{V} \otimes \theta_{V})c_{V,V}$ is a $K$--module
endomorphism of $\Hom_{\mV}(V \otimes V, \Io )$, so is a multiplication by a certain
$\vep_{i} \in K$. By the definition of the braiding and twist we have
$(\vep_{i})^{2} =1$. In particular $\vep_{i} \in \{ \pm 1 \}$ if $K$ is a field.
We have a distinguished element $0 \in I$ such that $V_{0}=\Io$.

The $S$-- and $T$--matrices of $\mV$ are the matrices $S=(S_{ij})_{i,j \in I}$,
$T=(T_{ij})_{i,j \in I}$ given by 
$S_{ij}=\tr(c_{V_{j},V_{i}} \circ c_{V_{i},V_{j}})$ and
$T_{ij}=\delta_{ij}v_{i}$, where $\tr$ is the categorical trace of $\mV$,
$\delta_{ij}$ is the Kronecher delta equal to $1$ if $i=j$ and zero elsewhere,
and $v_{i} \in K$ such that $\theta_{V_{i}}=v_{i}\id_{V_{i}}$.

Assume that $\mV$ has a rank $\mD$, i.e.\ an element of $K$ satisfying
$$
\mD^{2}=\sum_{i \in I} \dim(i)^{2},
$$
where $\dim(i)=\dim(V_{i})=\tr(\id_{V_{i}})$. We let
$$
\Delta = \sum_{i \in I} v_{i}^{-1} \dim(i)^{2}.
$$
Moreover, let $\tau=\tau_{(\mV,\mD)}$ be the RT--invariant
associated to  $\left( \mV, \{ V_{i} \}_{i \in I },\mD\right)$, cf.\ \cite[Sect.~II.2]{Turaev}.
For a tuple of integers $\mC =(m_{1},\ldots,m_{t})$, let
$$
G^{\mC} = T^{m_{t}}S\cdots T^{m_{1}}S.
$$
We put $a_{\os}=2$ and $a_{\ns}=1$. Moreover, let
$b_{j}^{(\os)}=1$ and $b_{j}^{(\ns)}=\delta_{jj^{*}}$, $j \in I$.
Given pairs $(\alpha_{j},\beta_{j})$ of coprime integers we
let $\mC_{j}=(a_{1}^{(j)},\ldots,a_{m_{j}}^{(j)})$ be a continued
fraction expansion of $\alpha_{j}/\beta_{j}$,
$j=1,2,\ldots,n$.

\begin{thm}[\cite{Hansen2}]\label{invariants}
Let
$M=(\ep;g\;|\;b;(\alpha_{1},\beta_{1}),\ldots,(\alpha_{n},\beta_{n}))$, $\ep=\os,\ns$.
Then
\begin{eqnarray*}
&&\tau(M) = (\Delta\mD^{-1})^{\sigma_{\ep}} \mD^{a_{\ep}g-2-\sum_{j=1}^{n} m_{j}} \\
 && \hspace{1.0in} \times \sum_{j \in I} \left(\vep_{j}\right)^{a_{\ep}g} b_{j}^{(\ep)} v_{j}^{-b} \dim(j)^{2-n-a_{\ep}g} \left( \prod_{i=1}^{n} (SG^{\mC_{i}})_{j0} \right),
\end{eqnarray*}
where
\begin{equation}\label{eq:signature}
\sigma_{\ep}=(a_{\ep}-1)\sign(E) + \sum_{j=1}^{n} \sign(\alpha_{j}\beta_{j}) + \frac{1}{3} \sum_{j=1}^{n} \left( \sum_{k=1}^{m_{j}} a_{k}^{(j)} -\Phi(B^{\mC_{j}}) \right).
\end{equation}
Here $E=-\left( b+\sum_{j=1}^{n} \frac{\beta_{j}}{\alpha_{j}} \right)$ is the Seifert Euler
number.

The RT--invariant $\tau(M)$ of the Seifert manifold $M$ with
non-normalized Seifert invariants
$\{ \ep;g;(\alpha_{1},\beta_{1}),\ldots,(\alpha_{n},\beta_{n})\}$,
is given by the same expression with the exceptions that the factor $v_{j}^{-b}$
has to be removed and $E=-\sum_{j=1}^{n} \frac{\beta_{j}}{\alpha_{j}}$.\HS
\end{thm}

The theorem is also valid in case $n=0$. In this case one just
has to put all sums $\sum_{j=1}^{n}$
equal to zero and all products $\prod_{i=1}^{n}$ equal to $1$.
Note that $\ep_{j}^{g}=1$ if $g$ is even and $\ep_{j}^{g}=\ep_{j}$
if $g$ is odd since $\ep_{j}^{2}=1$. The sum
$\sum_{j=1}^{n} \sign(\alpha_{j}\beta_{j})$ is of course
equal to $n$ for normalized Seifert invariants.

Let us next consider the lens spaces. For $p,q$ a pair of
coprime integers, recall that $L(p,q)$ is given by surgery on
$S^{3}$ along the unknot with surgery coefficient $-p/q$.
In the following corollary we include the possibilities
$L(0,1)=S^{1} \times S^{2}$ and $L(1,q)=S^{3}$, $q \in \Z$.

\begin{cor}[\cite{Hansen2}]\label{cor-lens-spaces}
Let $p,q$ be a pair of coprime integers. If $q \neq 0$ we let
$(a_{1},\ldots,a_{n-1})$ be
a continued fraction expansion of $-p/q$.
If $q=0$, put $n=3$ and $a_{1}=a_{2}=0$. 
Then
\begin{equation}\label{eq:lens spaces}
\tau(L(p,q))=(\Delta\mD^{-1})^{\sigma} \mD^{-n} G^{\mC}_{00},
\end{equation}
where $\mC=(a_{1},\ldots,a_{n-1},0)$ and
$\sigma=\frac{1}{3} \left( \sum_{l=1}^{n-1} a_{l} - \Phi(B^{\mC}) \right)$.\HS
\end{cor}

\rk{The RT--invariants of the Seifert manifolds for the classical Lie algebras}
Let $\frg$ be a complex finite dimensional simple Lie algebra
and let $q=e^{\pi \I/r}$, $r=m \kappa$, where $\kappa$ is an integer $\geq h^{\vee}$.
Let $U_{q}(\frg)$ be the quantum group associated to these data
as defined by Lusztig, see \cite[Part V]{Lusztig}. We follow
\cite[Sect.~1.3 and 3.3]{BakalovKirillov} here but will mostly
use notation from \cite{Turaev} for modular categories as above.
(Note that what we denote
$U_{q}(\frg)$ here is denoted $U_{q}(\frg)|_{q=e^{\pi\I/m \kappa}}$
in \cite{BakalovKirillov}.)
Let $\left(\mV_{r}^{\frg}, \{ V_{i} \}_{i \in I} \right)$ be the modular category
induced by these data, cf.\ \cite[Theorem 3.3.20]{BakalovKirillov}.
In particular the index set for the simple objects
is $I=\inte(C_{\kappa}) \cap X$.
We use here the shifted indexes (shifted by $\rho$) (contrary to \cite{BakalovKirillov}).
Normally the irreducible modules
of $U_{q}(\frg)$ (of type $1$) ($q$ a formal variable) are indexed by
the cone of dominant integer weights $X_{+}$. Here we
denote the irreducible module associated with $\mu \in X_{+}$
by $V_{\mu + \rho}$. For $q$ a root of unity as above $V_{\lambda}$
is an irreducible module of $U_{q}(\frg)$ of non-zero dimension if
$\lambda \in I$.
The involution $I \to I$, $\lambda \mapsto \lambda^{*}$
is given by
$\lambda^{*} = -w_{0}(\lambda)$, where $w_{0}$
is the longest element in $W$, see also \refrem{rem:unitarity}.
The distinguished element $0 \in I$
is equal to $\rho$. According to \cite[Theorem 3.3.20]{BakalovKirillov}
we can use
\begin{equation}\label{eq:rank}
\mD=\kappa^{l/2} \left| \frac{\vol (\cY) }{\vol (X) }\right|^{1/2}
  \left( \prod_{\alpha \in \Delta_{+}} 2\sin \left( \frac{\pi\la \alpha,\rho \ra }{\kappa} \right) \right)^{-1}
\end{equation}
as a rank of $\mV_{r}^{\frg}$.
According to the same theorem we have that
\begin{equation}\label{eq:anomaly}
\Delta \mD^{-1} = \omega^{-3},
\end{equation}
where
\begin{equation}\label{eq:omega}
\omega = e^{\frac{2\pi\I c}{24}} = \exp\left( \frac{\pi\I}{\ch} |\rho|^{2} \right) \exp\left( -\frac{\pi\I}{\kappa} |\rho|^{2} \right),
\end{equation}
where as usual $c=\frac{\kappa-\ch}{\kappa}\dim(\frg)$ is the central
charge and where the last equality follows from Freudenthal's strange formula
(\ref{eq:Freudenthal}).

Let $\tilde{s}$ be as defined in \cite[Sect.~3.1]{BakalovKirillov}
and let $S$ be the $S$--matrix of $\mV_{r}^{\frg}$.
Observe that in the terminology of
\cite{Turaev}, $\tilde{s}$ is the $S$--matrix of the mirror of
$\mV_{r}^{\frg}$. This implies that
$$
S_{\lambda\mu}=\tilde{s}_{\lambda^{*}\mu}=\tilde{s}_{\lambda\mu^{*}}
$$
for $\lambda,\mu \in I$. By \cite[Formula (3.1.16)]{BakalovKirillov}
we have $\tilde{s}=\mD s$, where $s$ is the matrix also considered in \refrem{rem:unitarity}.
The matrix $t$ in \cite{BakalovKirillov} is equal to
the $T$--matrix of $\mV_{r}^{\frg}$. It follows then
from \cite[Theorem 3.3.20]{BakalovKirillov} and (\ref{eq:mR})
that $T=\omega \tilde{\Theta}$. By these remarks and (\ref{eq:unitarity1}) and
(\ref{eq:unitarity3})
we conclude that
\begin{equation}\label{eq:ST}
S_{\lambda \mu} = \mD \tilde{\Xi}_{\lambda \mu},\hspace{.2in} T_{\lambda \mu} = \omega \tilde{\Theta}_{\lambda \mu}
\end{equation}
for $\mu,\lambda \in I$.
Let $\mC=(a_{1},\ldots,a_{n}) \in \Z^{n}$ and let $k \in \{0,1\}$.
By (\ref{eq:ST}) we immediately get
\begin{equation}\label{eq:sumformula}
(S^{k}G^{\mC})_{\lambda\rho} = \mD^{k+n}\omega^{\sum_{j=1}^{n} a_{j}} \left( \tilde{\Xi}^{k} \tilde{\Theta}^{a_{n}}\tilde{\Xi}\tilde{\Theta}^{a_{n-1}}\cdots\tilde{\Theta}^{a_{1}}\tilde{\Xi}\right)_{\lambda\rho}
\end{equation}
for $\lambda \in I$.
We also need the explicit formula for $\dim(\lambda)$, $\lambda \in I$.
In fact
\begin{equation}\label{eq:dim}
\dim(\lambda)=S_{\lambda\rho}=\mD \tilde{\Xi}_{\lambda\rho}=\mD \kappa^{-l/2}
  \left| \frac{\vol (X) }{\vol (\cY) }\right|^{1/2}
  \prod_{\alpha \in \Delta_{+}} 2\sin \left( \frac{\pi\la \alpha,\lambda \ra }{\kappa} \right),
\end{equation}
see also \cite[Formulas (3.3.2) and (3.3.5)]{BakalovKirillov}.

Let $\tau_{r}^{\frg}=\tau_{(\mV_{r}^{\frg},\mD)}$ be the RT--invariant
associated with $(\mV_{r}^{\frg},\mD)$.
Given a pair of coprime integers $(\alpha,\beta)$,
$\alpha>0$, we let $\beta^{*}$ be the inverse of $\beta$ in the multiplicative
group of units of $\Z/\alpha\Z$. For integers $a, b \neq 0$ we let
\begin{equation}\label{eq:dedekindsymbol}
\dS (a/b) = 12\sign(b)\s (a,b),
\end{equation}
where $\s (a,b)$ is given by (\ref{eq:dedekindsum}). This is the
so-called Dedekind symbol. (In particular, the right-hand
side of (\ref{eq:dedekindsymbol}) only depends on the rational
number $a/b$.) Then we have the following generalization
of \cite[Theorem 8.4]{Hansen2}:

\begin{thm}\label{Lie-Seifert}
Let $M=(\ep;g\;|\;b;(\alpha_{1},\beta_{1}),\ldots,(\alpha_{n},\beta_{n}))$,
$\ep \in \{ \os, \ns\}$. Then
\begin{eqnarray*}
\tau_{r}^{\frg}(M) &=& \exp \left( \frac{\pi\I}{\kappa}\left[ 3(a_{\ep}-1)\sign(E) -E - \sum_{j=1}^{n} \dS \left(\frac{\beta_{j}}{\alpha_{j}}\right) \right] |\rho|^{2} \right) \\
 && \hspace{.3in} \times \frac{\I^{n\npr} \kappa^{l(a_{\ep}g/2-1)}}{2^{\npr(n+a_{\ep}g-2)}\vol(\cY)^{2-a_{\ep}g}}
  \frac{1}{\mA^{l/2}} \\
 && \hspace{.3in} \times \exp \left( \frac{3\pi\I}{\ch}(1-a_{\ep})\sign(E) |\rho|^{2} \right) Z_{\ep}^{\frg}(M;r),
\end{eqnarray*}
where $\mA=\prod_{j=1}^{n} \alpha_{j}$ and
\begin{eqnarray*}
Z_{\ep}^{\frg}(M;r) &=& \sum_{\lambda \in I} b_{\lambda}^{(\ep)} \vep_{\lambda}^{a_{\ep}g}
\left( \prod_{\alpha \in \Delta_{+}} 
\sin ^{2-n-a_{\ep}g} \left( \frac{\pi\la \lambda,\alpha \ra}{\kappa} \right)\right)
\exp \left( \frac{\pi\I}{\kappa} E |\lambda|^{2} \right) \\
 & & \hspace{.2in} \times \sum_{w_{1},\ldots,w_{n}\in W} 
  \sum_{\nu_{1} \in \cY/\alpha_{1} \cY} \ldots \sum_{\nu_{n} \in \cY/\alpha_{n} \cY}
  \left(\prod_{j=1}^{n} \det(w_{j})\right) \\
 & & \hspace{0.6in} \times
     \exp \left( -\pi \I \sum_{j=1}^{n} \frac{\beta_{j}^{*}}{\alpha_{j}} \left( \kappa|\nu_{j}|^{2} 
      + 2\la w_{j}(\rho),\nu_{j} \ra \right) \right) \\
 & & \hspace{0.6in} \times  
     \exp \left( - \frac{2 \pi\I}{\kappa} \la \lambda,  \sum_{j=1}^{n} \frac{ \kappa\nu_{j} + w_{j}(\rho)}{\alpha_{j}} \ra \right).
\end{eqnarray*}
The RT--invariant $\tau_{r}^{\frg}(M)$ of the Seifert manifold $M$
with non-normalized Seifert invariants
$\{\ep;g;(\alpha_{1},\beta_{1}),\ldots,$ $(\alpha_{n},\beta_{n})\}$
is given by the same expression.\HS
\end{thm}

The theorem is also valid in case $n=0$. In this case one
just has to put the sum
$\sum_{w_{1},\ldots,w_{n}\in W}
\sum_{\nu_{1} \in \cY/\alpha_{1} \cY} \ldots \sum_{\nu_{n} \in \cY/\alpha_{n} \cY}$
in $Z_{\ep}(M;r)$ equal to $1$,
$\ep=\os,\ns$,
and put $\mA=1$ and
$\sum_{j=1}^{n}\dS (\beta_{j}/\alpha_{j})=0$.

\begin{proof}
The proof follows very closely the proof of \cite[Theorem 8.4]{Hansen2}.
Let $M=(\ep;g\;|\;b;(\alpha_{1},\beta_{1}),\ldots,(\alpha_{n},\beta_{n}))$,
$\ep \in \{ \os, \ns\}$. Choose tuples of integers
$\mC_{j}=(a_{1}^{(j)},\ldots,a_{m_{j}}^{(j)})$ such that
$B^{\mC_{j}}=\left( \begin{array}{cc}
                \alpha_{j} & \rho_{j} \\
		\beta_{j} & \sigma_{j}
		\end{array}
\right)$, $j=1,2,\ldots,n$.
By \refthm{invariants} and (\ref{eq:continued}) we have
\begin{eqnarray*}
\tau_{r}^{\frg}(M) &=& (\Delta\mD^{-1})^{\sigma_{\ep}} \mD^{a_{\ep}g-2-\sum_{j=1}^{n} m_{j}} \\
 && \hspace{.6in} \times \sum_{\lambda \in I} b_{\lambda}^{(\ep)}\vep_{\lambda}^{a_{\ep}g} 
             v_{\lambda}^{-b} \dim(\lambda)^{2-n-a_{\ep}g} \left( \prod_{i=1}^{n} (SG^{\mC_{i}})_{\lambda\rho} \right),
\end{eqnarray*}
where $\sigma_{\ep}$ is given by (\ref{eq:signature}).
By using (\ref{eq:anomaly}), (\ref{eq:ST}), (\ref{eq:sumformula}), and (\ref{eq:dim}) we get
\begin{eqnarray*}
\tau_{r}^{\frg}(M) &=& \alpha_{\ep}(\kappa) \omega^{-3\sigma_{\ep}-b+\sum_{j=1}^{n}\sum_{k=1}^{m_{j}} a_{k}^{(j)}} \sum_{\lambda \in I} b_{\lambda}^{(\ep)}\vep_{\lambda}^{a_{\ep}g} \tilde{\Theta}_{\lambda\lambda}^{-b} \\
 && \hspace{.8in} \times \left( \prod_{\alpha \in \Delta_{+}} \sin^{2-n-a_{\ep}g}\left( \frac{\pi \la \alpha,\lambda \ra}{\kappa} \right) \right) \left( \prod_{i=1}^{n} (\tilde{N}_{i})_{\lambda\rho} \right),
\end{eqnarray*}
where $N_{i}=\Xi B^{\mC_{i}}=\left( \begin{array}{cc}
                -\beta_{i} & -\sigma_{i} \\
		\alpha_{i} & \rho_{i}
		\end{array}
\right)$ and
$$
\alpha_{\ep}(\kappa)=\kappa^{l(n+a_{\ep}g-2)/2} 2^{\npr(2-n-a_{\ep}g)} \vol (\cY)^{n+a_{\ep}g-2}.
$$
By \refcor{cor:unitarity} we get
\begin{eqnarray*}
\prod_{i=1}^{n}(\tilde{N}_{i})_{\lambda\rho} &=& \beta_{\ep}(\kappa) \sum_{w_{1},\ldots,w_{n}\in W}
  \left(\prod_{j=1}^{n} \det(w_{j})\right) \sum_{\nu_{1} \in \cY/\alpha_{1} \cY} \ldots \sum_{\nu_{n} \in \cY/\alpha_{n} \cY} \\
 && \hspace{.6in} \times
     \exp \left( \frac{\pi \I}{\kappa} \sum_{j=1}^{n} \frac{\rho_{j}}{\alpha_{j}} |\rho+\kappa\nu_{j}|^{2} \right) \\
 & & \hspace{.6in} \times
     \exp \left( - \frac{2 \pi\I}{\kappa} \sum_{j=1}^{n} \frac{1}{\alpha_{j}} \la \rho + \kappa\nu_{j}, w_{j}(\lambda) \ra \right),
\end{eqnarray*}
where
\begin{eqnarray*}
\beta_{\ep}(\kappa) &=& \frac{\I^{n\npr}}{\kappa^{nl/2}\mA^{l/2}\vol(\cY)^{n}}
    \exp \left( -\frac{\pi \I}{\ch} \left(\sum_{j=1}^{n} \Phi(N_{j})\right) |\rho|^{2}  \right) \\
 && \hspace{.5in} \times \exp \left( -\frac{\pi \I}{\kappa} \left( \sum_{j=1}^{n} \frac{\beta_{j}}{\alpha_{j}} \right) |\lambda|^{2}\right).
\end{eqnarray*}
In particular we have
\begin{eqnarray*}
\alpha_{\ep}(\kappa)\beta_{\ep}(\kappa)\tilde{\Theta}_{\lambda\lambda}^{-b} &=& \frac{ \I^{n\npr} \kappa^{l(a_{\ep}g-2)/2}}{2^{\npr(n+a_{\ep}g-2)}\mA^{l/2}\vol(\cY)^{2-a_{\ep}g}} \exp \left(\frac{\pi \I}{\kappa} E |\lambda|^{2} \right) \\
 && \hspace{.8in} \times \exp\left( \frac{\pi \I}{\ch} \left( b -\sum_{j=1}^{n} \Phi(N_{j}) \right)|\rho|^{2} \right).
\end{eqnarray*}
By (\ref{eq:Phiproduct}), $\Phi(N_{i})=\Phi(B^{\mC_{i}})-3\sign(\alpha_{i}\beta_{i})$.
By this and (\ref{eq:signature}) we get
$$
3\sigma_{\ep} = 3(a_{\ep}-1)\sign(E) - \sum_{j=1}^{n} \Phi(N_{j}) + \sum_{j=1}^{n}\sum_{k=1}^{m_{j}} a_{k}^{(j)}.
$$
Therefore
{\allowdisplaybreaks
\begin{eqnarray*}
&&\tau_{r}^{\frg}(M) = \frac{\I^{n\npr} \kappa^{l(a_{\ep}g/2-1)}}{2^{\npr(n+a_{\ep}g-2)}\vol(\cY)^{2-a_{\ep}g}}
\frac{1}{\mA^{l/2}} \omega^{\sum_{j=1}^{n} \Phi(N_{j})-3(a_{\ep}-1)\sign(E)-b} \\*
 && \hspace{1.0in} \times \exp\left( \frac{\pi \I}{\ch} \left( b -\sum_{j=1}^{n} \Phi(N_{j}) \right)|\rho|^{2} \right) \\*
 && \hspace{1.0in} \times \exp\left(\frac{\pi\I}{\kappa} \left( \sum_{j=1}^{n} \frac{\rho_{j}}{\alpha_{j}} \right) |\rho|^{2} \right) \\
 && \hspace{.6in} \times \sum_{\lambda \in I} b_{\lambda}^{(\ep)}\vep_{\lambda}^{a_{\ep}g}
         \left( \prod_{\alpha \in \Delta_{+}} 
         \sin ^{2-n-a_{\ep}g} \left( \frac{\pi\la \lambda,\alpha \ra}{\kappa} \right)\right)
         \exp \left(\frac{\pi \I}{\kappa} E |\lambda|^{2} \right) \\
 && \hspace{.6in} \times \sum_{w_{1},\ldots,w_{n}\in W}
         \sum_{\nu_{1} \in \cY/\alpha_{1} \cY} \ldots \sum_{\nu_{n} \in \cY/\alpha_{n} \cY} 
         \left(\prod_{j=1}^{n} \det(w_{j})\right) \\
 && \hspace{1.0in} \times
     \exp \left( \pi \I \sum_{j=1}^{n} \frac{\rho_{j}}{\alpha_{j}} (\kappa|\nu_{j}|^{2} +2\la w_{j}(\rho),\nu_{j}\ra \right) \\*
 && \hspace{1.0in} \times
     \exp \left( - \frac{2 \pi\I}{\kappa} \la \lambda, \sum_{j=1}^{n} \frac{\kappa\nu_{j}+w_{j}(\rho)}{\alpha_{j}} \ra \right).
\end{eqnarray*}}\noindent
By (\ref{eq:omega}) we get
\begin{eqnarray*}
&&\omega^{\sum_{j=1}^{n} \Phi(N_{j})-3(a_{\ep}-1)\sign(E)-b}
         \exp\left( \frac{\pi \I}{\ch} \left( b -\sum_{j=1}^{n} \Phi(N_{j}) \right)|\rho|^{2} \right) \\
 && \hspace{.8in} \times \exp\left(\frac{\pi\I}{\kappa} \left( \sum_{j=1}^{n} \frac{\rho_{j}}{\alpha_{j}} \right) |\rho|^{2} \right) \\
 && \hspace{.2in} = \exp\left( \frac{3\pi\I}{\ch} (1-a_{\ep}) \sign(E)|\rho|^{2} \right) \\
 && \hspace{.3in} \times
     \exp\left( \frac{\pi\I}{\kappa} \left[ 3(a_{\ep}-1)\sign(E) + b + \sum_{j=1}^{n} \left( \frac{\rho_{j}}{\alpha_{j}} - \Phi(N_{j}) \right)\right] |\rho|^{2} \right).
\end{eqnarray*}
The theorem now follows by using (\ref{eq:Phi}) together with the
facts that $\s (a,b)=\s (a',b)$ if $a'a \equiv 1 \pmod{b}$
and $\s (-a,b)=-\s (a,b)$,
cf.\ \cite[Chap.~3]{RademacherGrosswald}. (The identity
$\s (-a,b)=-\s (a,b)$ follows immediately from (\ref{eq:dedekindsum}).)
The case with non-normalized Seifert invariants follows as above by
letting $b$ be equal to zero everywhere.
\end{proof}

By the above proof we get the following compact formula for the invariant
of the Seifert manifold 
$M=(\ep;g\;|\;b;(\alpha_{1},\beta_{1}),\ldots,(\alpha_{n},\beta_{n}))$:
\begin{eqnarray}\label{eq:compact}
\tau_{r}^{\frg}(M) &=& \gamma_{\ep}(\kappa)\sum_{\lambda \in I} b_{\lambda}^{(\ep)}\vep_{\lambda}^{a_{\ep}g}
        \exp\left( -\frac{\pi\I}{\kappa}b|\lambda|^{2}\right) \\
 && \hspace{.6in} \times \left(\tilde{\Xi}_{\lambda\rho}\right)^{2-n-a_{\ep}g}
        \left( \prod_{i=1}^{n} (\tilde{N}_{i})_{\lambda\rho} \right), \nonumber
\end{eqnarray}
where $N_{i} \in \SL (2,\Z)$ with first column equal to
$\left( \begin{array}{c}
                -\beta_{i} \\
		\alpha_{i}
		\end{array}
\right)$, $i=1,\ldots,n$, and
\begin{eqnarray*}
\gamma_{\ep}(\kappa) &=&
  \exp\left( \frac{3\pi\I}{\ch}\left[(1-a_{\ep})\sign(E)+\sum_{j=1}^{n} \Phi(N_{j}) \right] |\rho|^{2}\right) \\
 && \hspace{.3in} \times \exp \left( \frac{\pi\I}{\kappa}\left[ 3(a_{\ep}-1)\sign(E) +b - \sum_{j=1}^{n} \Phi(N_{j}) \right] |\rho|^{2} \right).
\end{eqnarray*}

\rk{The coprime case} In this section we will derive a particularly
nice expression for the quantum $\frg$--invariants of the Seifert manifolds
$M=(\ep;g\;|\;b;(\alpha_{1},\beta_{1}),\ldots,$ $(\alpha_{n},\beta_{n}))$
with the integers $\alpha_{1},\ldots,\alpha_{n}$ mutually coprime.
The subclass of these Seifert manifolds additionally satisfying
that $g=0$ contains the Seifert fibered integral homology spheres,
cf.\ \cite[Corollary 6.2 and pp.~36--37]{JankinsNeumann}.
The results in this subsection generalize the results
obtained in \cite[Sect.~4.1]{LawrenceRozansky}, where the case $\frg=\frsl_{2}(\C)$
is considered. 
Let us first observe that by (\ref{eq:compact}) and \reflem{lem:mRsymmetry}
we have
\begin{eqnarray*}\label{eq:compact1}
\tau_{r}^{\frg}(M) &=& \frac{\gamma_{\ep}(\kappa)}{|W|}
        \sum_{\lambda \in \tilde{P}_{\kappa}\cap X} b_{\lambda}^{(\ep)}\vep_{\lambda}^{a_{\ep}g}
        \exp\left( -\frac{\pi\I}{\kappa}b|\lambda|^{2}\right) \\
 && \hspace{0.9in} \times \left(\tilde{\Xi}_{\lambda\rho}\right)^{2-n-a_{\ep}g}
        \left( \prod_{i=1}^{n} (\tilde{N}_{i})_{\lambda\rho} \right), \nonumber
\end{eqnarray*}
where $\tilde{P}_{\kappa} = P_{\kappa} \sm H^{\kappa}$.
To establish this identity we also extended the function $g \co I \to \C$,
$g(\lambda)=b_{\lambda}^{(\ep)}\vep_{\lambda}^{a_{\ep}g}$, to all of $X$ by letting it
be constantly $1$ if $\ep=\os$ and by forcing it to satisfy the following
symmetry properties if $\ep=\ns$:
\begin{eqnarray}\label{eq:gsymmetry}
g(w(\lambda)) &=& \det(w)^{a_{\ep}g}g(\lambda), \hspace{.2in} w \in W, \, \lambda \in X, \\
g(\lambda +x) &=& g(\lambda), \hspace{.2in} x \in \kappa \cY, \, \lambda \in X, \nonumber \\
g(\lambda) &=& 0, \hspace{.2in} \lambda \in X \cap \cup_{\alpha \in \Delta_{+}, b \in \Z} H_{\alpha,b}^{\kappa} \nonumber.
\end{eqnarray}

\begin{rem}
If $\ep =\ns$ we have
$$
\mD^{-2}\sum_{i \in I} \dim(i)F(L_{i\lambda})=b_{\lambda}^{(\ns)}\vep_{\lambda},
$$
by \cite[Lemma 4.2]{Hansen2},
where $L_{i\lambda}$ is a certain 2--component link, one component
colored by $V_{i}$ and the other
by $V_{\lambda}$. Recall that $b_{\lambda}^{(\ns)}=\delta_{\lambda\lambda^{*}}$.
In \cite{Le1} it was observed that one in a natural way can extend
the invariant $F(L_{i\lambda})$ to be defined for all $\lambda \in X$ in such a way
that  $F(L_{iw(\lambda)})=\det(w)F(L_{i\lambda})$, $w \in W$.
The thus extended function $\lambda \to F(L_{i\lambda})$
is automatically invariant under a change $\lambda \to \lambda + x$, $x \in \kappa \cY$.
(One has to be a little carefull here. In \cite{Le1},
$F(L_{i\lambda})$ is denoted $J_{L}(i,\lambda)$. It is shown
in \cite{Le1} that this function is component-wise invariant
under the action by the affine Weyl group $\aW_{\kappa}$
up to a sign comming from the sign of the Weyl group element, see
the {\it refined} first symmetry principle \cite[Theorem 2.11]{Le1}.
The lattice denoted $Y'$ in \cite{Le1} is equal to $\frac{1}{m}\cY$ here ($m$
is denoted $d$ in \cite{Le1}).
Since $r$ in \cite{Le1} is the same as $r$ here, i.e.\ $r=m\kappa$,
we have that $W \ltimes rY'$ in \cite{Le1} is
equal to $\aW_{\kappa}$ here.) We see, that the symmetries
(\ref{eq:gsymmetry}) actually are consequences of the symmetry results in
\cite{Le1}.
\end{rem}

By \refthm{Lie-Seifert} and (\ref{eq:compact1})
(use also $E= -b-\sum_{j=1}^{n} \beta_{j}/\alpha_{j}$) we have
\begin{eqnarray*}
Z_{\ep}^{\frg}(M;r) &=& \alpha(\kappa)\sum_{\lambda \in \tilde{P}_{\kappa} \cap X} b_{\lambda}^{(\ep)} \vep_{\lambda}^{a_{\ep}g}
     \left( \prod_{\alpha \in \Delta_{+}} 
     \sin ^{2-n-a_{\ep}g} \left( \frac{\pi\la \lambda,\alpha \ra}{\kappa} \right)\right) \\
 && \hspace{.3in} \times \exp\left(-\frac{\pi\I}{\kappa}b|\lambda|^{2}\right) \sum_{w_{1},\ldots,w_{n}\in W} \\
 && \hspace{.3in} \times \sum_{\nu_{1} \in \cY/\alpha_{1} \cY} \ldots \sum_{\nu_{n} \in \cY/\alpha_{n} \cY} 
     \left( \prod_{j=1}^{n} F_{j}(\lambda,\nu_{j},w_{j}) \right),
\end{eqnarray*}
where
$$
\alpha(\kappa)=\frac{1}{|W|}\exp\left( -\frac{\pi\I}{\kappa} \left( \sum_{j=1}^{n} \frac{\beta_{j}^{*}}{\alpha_{j}}\right) |\rho|^{2}\right)
$$
and
\begin{eqnarray*}
&&F_{j}(\lambda,\nu,w) = \det(w) \\
 &&\hspace{.6in} \times \exp \left( -\frac{\pi \I}{\kappa\alpha_{j}} 
 \left[\beta_{j}|\lambda|^{2} + 2\la \lambda,\kappa\nu+ w(\rho) \ra + \beta_{j}^{*} |\kappa\nu+w(\rho)|^{2} \right] \right).
\end{eqnarray*}
By using (\ref{eq:integer1}) and (\ref{eq:integer2}) one finds the
following symmetry result:

\begin{lem}\label{lem:Fsymmetries}
Let $j \in \{1,2,\ldots,n\}$. The map $F_{j} \co X \times \cY \times W \to \C$ is invariant under the transformations
\begin{description}
\item[(a)] $(\lambda,\nu,w) \to (\lambda \pm \kappa x, \nu \mp \beta_{j}x,w)$,
\item[(b)] $(\lambda,\nu,w) \to (\lambda,\nu \pm \alpha_{j}x,w)$,
\item[(c)] $(\lambda,\nu,w) \to (\lambda \pm \kappa \alpha_{j}x,\nu,w)$,
\end{description}
for any $x \in \cY$.\HS
\end{lem}

By (\ref{eq:integer1}) we get that
$\prod_{\alpha \in \Delta_{+}} \sin ^{2-n-a_{\ep}g} \left( \frac{\pi\la \lambda,\alpha \ra}{\kappa} \right)$
is invariant under a change $\lambda \to \lambda + \kappa x$, $x \in \cY$.

\begin{cor}\label{cor:Hsymmetries}
For an arbitrary but fixed $j_{0} \in \{1,2,\ldots,n\}$ and
an arbitrary $x \in \cY$, the expression
\begin{eqnarray*}
H(\lambda,\nu_{1},\ldots,\nu_{n}) &:=& b_{\lambda}^{(\ep)} \vep_{\lambda}^{a_{\ep}g}
     \left( \prod_{\alpha \in \Delta_{+}} 
     \sin ^{2-n-a_{\ep}g} \left( \frac{\pi\la \lambda,\alpha \ra}{\kappa} \right)\right) \\
 && \hspace{.1in} \times \exp\left(-\frac{\pi\I}{\kappa}b|\lambda|^{2}\right) \sum_{w_{1},\ldots,w_{n}\in W}
     \left( \prod_{j=1}^{n} F_{j}(\lambda,\nu_{j},w_{j}) \right)
\end{eqnarray*}
is invariant under the transformation
\begin{eqnarray*}
\lambda &\to& \lambda \pm \kappa \frac{\mA}{\alpha_{j_{0}}} x, \\
\nu_{j} &\to& \nu_{j} \mp \beta_{j_{0}} \frac{\mA}{\alpha_{j_{0}}} \delta_{jj_{0}} x, \hspace{.2in} j=1,2,\ldots,n,
\end{eqnarray*}
where $\mA = \prod_{j=1}^{n} \alpha_{j}$.\HS
\end{cor}

Choose integers $k_{1},\ldots,k_{n}$, $l_{1},\ldots,l_{n}$,
$a_{1},\ldots,a_{n}$, and $b_{1},\ldots,b_{n}$, such that
\begin{eqnarray*}
k_{j}\beta_{j} + l_{j} \alpha_{j} &=& 1, \\
a_{j}\alpha_{j} + b_{j}\frac{\mA}{\alpha_{j}} &=& 1
\end{eqnarray*}
for $j=1,2,\ldots,n$.
Let $(\nu_{1},\ldots,\nu_{n}) \in \cY/\alpha_{j}\cY \times \ldots \cY/\alpha_{n}\cY$.
Let us show, that for any $\lambda \in \tilde{P}_{\kappa} \cap X$ there exists a 
$\lambda ' \in X$ such that $H(\lambda,\nu_{1},\ldots,\nu_{n})=H(\lambda',0,\ldots,0)$.
To see this, let $x \in \cY$ and $j \in \{1,2,\ldots,n\}$ be arbitrary but fixed. First observe
that by \reflem{lem:Fsymmetries}, $H$ is invariant under the transformation
$\nu_{j} \to \nu_{j} \pm \beta_{j}\alpha_{j}x$, $\lambda$ and $\nu_{k}$, $k \neq j$,
unchanged. By \refcor{cor:Hsymmetries},
$H$ is invariant under the transformation $\nu_{j} \to \nu_{j} \pm (\beta_{j} \mA/\alpha_{j})x$,
$\lambda \to \lambda \mp (\kappa\mA/\alpha_{j})x$, $\nu_{k}$, $k \neq j$, unchanged.
If we use the first transformation $a_{j}$ times and the second transformation
$b_{j}$ times we see that $H$ is invariant under the transformation
\begin{eqnarray*}
\lambda &\to& \lambda \mp (\kappa b_{j} \mA/\alpha_{j})x, \\
\nu_{k} &\to& \nu_{k} \pm \beta_{j}\delta_{jk}x, \hspace{.2in} k = 1,2,\ldots,n.
\end{eqnarray*}
By using this transformation $k_{j}$ times and by using the transformation
$\nu_{k} \to \nu_{k} \pm \alpha_{j}\delta_{jk}x$, $k=1,2,\ldots,n$,
$\lambda$ unchanged, $l_{j}$ times we see that $H$ is invariant under the
transformation
\begin{eqnarray*}
\lambda &\to& \lambda \mp (\kappa k_{j}b_{j} \mA/\alpha_{j})x, \\
\nu_{k} &\to& \nu_{k} \pm \delta_{jk}x, \hspace{.2in} k=1,2,\ldots,n.
\end{eqnarray*}
In particular we can change $\nu_{j}$ to $0$ and keep $H$ unchanged if we at the same time 
change $\lambda$ to $\lambda + (\kappa k_{j}b_{j}\mA/\alpha_{j})\nu_{j}$, so
$H$ is invariant under the transformation
\begin{eqnarray*}
\lambda &\to& \lambda + \kappa \sum_{j=1}^{n}k_{j}b_{j}\frac{\mA}{\alpha_{j}}\nu_{j}, \\
\nu_{k} &\to& 0, \hspace{.2in} k=1,2,\ldots,n.
\end{eqnarray*}
By using the above result, and the fact that $H$ is invariant
under the transformation $\lambda \to \lambda +\kappa\mA x$, $x \in \cY$, the
$\nu_{k}$ unchanged, we can always arrange it so that
$\lambda$ is an element of
$$
J := \tilde{P}_{\kappa} \cap X + \kappa \left\{ \; \sum_{i=1}^{l} m_{i} \alpha_{i}^{\vee} \; | \; m_{1},\ldots,m_{l} \in \{0,1,\ldots,\mA -1\} \;\right\}.
$$
We have a bijection $(P_{\kappa} \cap X) \times \cY / \mA \cY \to P_{\kappa\mA} \cap X$, $(\lambda,\mu) \mapsto \lambda+\kappa\mu$,
and $J$ is equal to the the image of $(\tilde{P}_{\kappa} \cap X) \times \cY / \mA \cY$ which can also be
identified with
$(P_{\kappa\mA} \sm H^{\kappa}) \cap X = X/\kappa\mA\cY \sm H^{\kappa}$.

\begin{lem}
Let $b_{j}$ and $k_{j}$ be as above. The map
\begin{eqnarray*}
&& \\
&&(\tilde{P}_{\kappa} \cap X) \times \cY/\alpha_{j}\cY \times \ldots \times \cY/\alpha_{n}\cY \hspace{.3in} \longrightarrow \hspace{.3in} J \\
&& \\
&& \hspace{.4in}(\lambda,\nu_{1},\ldots,\nu_{n}) \hspace{.1in}\mapsto \hspace{.1in} \lambda + \kappa\sum_{j=1}^{n} k_{j}b_{j} \frac{\mA}{\alpha_{j}} \nu_{j} + \kappa\mA \sum_{i=1}^{l} n_{i}\alpha_{i}^{\vee}
\end{eqnarray*}
is a bijection, where the $n_{i}$ are the unique integers such that the right-hand side
is an element of $J$.
\end{lem}

\begin{proof}
Assume that
$$
\lambda + \kappa\sum_{j=1}^{n} k_{j}b_{j} \frac{\mA}{\alpha_{j}} \nu_{j} + \kappa\mA \sum_{i=1}^{l} n_{i}\alpha_{i}^{\vee}=0,
$$
where $\nu_{j} \in \cY/\alpha_{j}\cY$, $j=1,2,\ldots,n$, and $\lambda$ is the difference
of two elements in $\tilde{P}_{\kappa} \cap X$.
Then we immediately get that $\lambda =0$ and
$$
\sum_{j=1}^{n} k_{j}b_{j} \frac{\mA}{\alpha_{j}} \nu_{j} + \mA \sum_{i=1}^{l} n_{i}\alpha_{i}^{\vee}=0.
$$
Write $\nu_{j}=\sum_{i=1}^{l} m_{i}^{(j)}\alpha_{i}^{\vee}$,
$m_{i}^{(j)} \in \{0,1,\ldots,\alpha_{j} -1\}$, and get
$$
\sum_{j=1}^{n} k_{j}b_{j} \frac{\mA}{\alpha_{j}} m_{i}^{(j)} + \mA n_{i} =0, \hspace{.2in} i=1,2,\ldots,l,
$$
since $\alpha_{1}^{\vee},\ldots,\alpha_{l}^{\vee}$ is a basis for $\frhR^{*}$. Now let $j_{0}$
be arbitrary but fixed. Then we have
$$
k_{j_{0}}b_{j_{0}}\frac{\mA}{\alpha_{j_{0}}} m_{i}^{(j_{0})} = -\sum_{j=1,j \neq j_{0}}^{n} k_{j}b_{j} \frac{\mA}{\alpha_{j}}m_{i}^{(j)} - \mA n_{i}
$$
for all $i \in \{1,2,\ldots,l\}$.
But $\alpha_{j_{0}}$ is a divisor of the right-hand side so is also a divisor of
$m_{i}^{(j_{0})}$, since $\alpha_{j_{0}}$
and $k_{j_{0}}b_{j_{0}}\mA/\alpha_{j_{0}}$ are coprime.
Therefore
$m_{i}^{(j_{0})}=0$, $i=1,2,\ldots,l$, so $\nu_{j_{0}}=0$.
It follows that
$\nu_{j}=0$ for all $j \in \{1,2,\ldots,n\}$ (and
$\sum_{i=1}^{l} n_{i}\alpha_{i}^{\vee}=0$, so $n_{1}=\ldots=n_{l}=0$).
The surjectivity follows now by the fact that $J$ and
$(\tilde{P}_{\kappa} \cap X) \times \cY/\alpha_{j}\cY \times \ldots \times \cY/\alpha_{n}\cY$
contain the same number of elements, namely $\mA^{l}|\tilde{P}_{\kappa} \cap X|$ elements. 
\end{proof}

By the above results we get that $Z_{\ep}^{\frg}(M;r)$ is given
by an expression obtained in the following way: Take the
original expression for $Z_{\ep}^{\frg}(M;r)$ as stated
in \refthm{Lie-Seifert}, divide it by $|W|$,
replace the summation index set $I$ by $J$,
remove the sum
$\sum_{\nu_{1} \in \cY/\alpha_{1} \cY} \ldots \sum_{\nu_{n} \in \cY/\alpha_{n} \cY}$ and put
all $\nu_{1},\ldots,\nu_{n}$ equal to zero, i.e.\
\begin{eqnarray*}
Z_{\ep}^{\frg}(M;r) &=& \frac{1}{|W|} \sum_{\lambda \in J} b_{\lambda}^{(\ep)} \vep_{\lambda}^{a_{\ep}g}
     \left( \prod_{\alpha \in \Delta_{+}} 
     \sin ^{2-n-a_{\ep}g} \left( \frac{\pi\la \lambda,\alpha \ra}{\kappa} \right)\right) \\*
 && \hspace{.3in} \times \exp\left( \frac{\pi\I}{\kappa} E|\lambda|^{2} \right) 
     \sum_{w_{1},\ldots,w_{n}\in W} \\*
 && \hspace{.3in} \times \left( \prod_{j=1}^{n} \det(w_{j} )
     \exp \left( -\frac{2\pi \I}{\kappa\alpha_{j}} \la \lambda,w_{j}(\rho) \ra \right)\right).
\end{eqnarray*}
Here
{\allowdisplaybreaks\begin{eqnarray*}
&&\sum_{w_{1},\ldots,w_{n}\in W} \left( \prod_{j=1}^{n} \det(w_{j} )
     \exp \left( -\frac{2\pi \I}{\kappa\alpha_{j}} \la \lambda,w_{j}(\rho) \ra \right) \right) \\
 && \hspace{.3in} = \prod_{j=1}^{n} \left(\sum_{w_{j} \in W} \det(w_{j})
     \exp \left( -\frac{2\pi \I}{\kappa\alpha_{j}} \la \lambda,w_{j}(\rho) \ra \right) \right) \\
 && \hspace{.3in} = \prod_{j=1}^{n} \prod_{\alpha \in \Delta_{+}} 
     2\I\sin \left( -\frac{\pi}{\kappa\alpha_{j}}\la \lambda,\alpha \ra \right),
\end{eqnarray*}}\noindent
so we have shown

\begin{thm}\label{Lie-Seifert-coprime}
Let $M=(\ep;g\;|\;b;(\alpha_{1},\beta_{1}),\ldots,(\alpha_{n},\beta_{n}))$,
$\ep \in \{ \os, \ns\}$, and assume that the $\alpha_{j}$ are
mutually coprime. Then
\begin{eqnarray*}
\tau_{r}^{\frg}(M) &=&\exp \left( \frac{\pi\I}{\kappa}\left[ 3(a_{\ep}-1)\sign(E) -E - \sum_{j=1}^{n} \dS \left(\frac{\beta_{j}}{\alpha_{j}}\right) \right] |\rho|^{2} \right) \\
  &&\hspace{.3in} \times \frac{\kappa^{l(a_{\ep}g/2-1)}}{2^{\npr(a_{\ep}g-2)}\vol(\cY)^{2-a_{\ep}g}|W|}
  \frac{1}{\mA^{l/2}} \\
  &&\hspace{.3in} \times \exp\left( \frac{3\pi\I}{\ch}(1-a_{\ep})\sign(E)|\rho|^{2}\right) W_{\ep}^{\frg}(M;r),
\end{eqnarray*}
where $\mA=\prod_{j=1}^{n} \alpha_{j}$ and
\begin{eqnarray*}
W_{\ep}^{\frg}(M;r) &=& \sum_{\lambda \in J} b_{\lambda}^{(\ep)} \vep_{\lambda}^{a_{\ep}g}
\left( \prod_{\alpha \in \Delta_{+}} 
\sin ^{2-n-a_{\ep}g} \left( \frac{\pi\la \lambda,\alpha \ra}{\kappa} \right)\right) \\
 & & \hspace{.2in} \times \exp \left( \frac{\pi\I}{\kappa} E |\lambda|^{2} \right) \prod_{j=1}^{n} \prod_{\alpha \in \Delta_{+}} 
     \sin \left( \frac{\pi}{\kappa\alpha_{j}}\la \lambda,\alpha \ra \right).
\end{eqnarray*}
The RT--invariant $\tau_{r}^{\frg}(M)$ of the Seifert manifold $M$
with non-normalized Seifert invariants
$\{\ep;g;(\alpha_{1},\beta_{1}),\ldots,$ $(\alpha_{n},\beta_{n})\}$
is given by the same expression.\HS
\end{thm}

\section{The case of lens spaces}\label{sec-The-case}

In this section we present different expressions for
the invariant $\tau_{r}^{\frg}(L(p,q))$, $p,q$ being an
arbitrary but fixed pair of coprime integers.
Let $b,d$ be any integers such that
$U=\left( \begin{array}{cc}
                q & b \\
		p & d
		\end{array}
\right) \in \SL (2,\Z)$. Assume $q \neq 0$, let
$$
V=-\Xi U =\left( \begin{array}{cc}
                p & d \\
		-q & -b
		\end{array}
\right),
$$
and let $C'=(a_{1},a_{2},\ldots,a_{n-1}) \in \Z^{n-1}$ such that
$B^{\mC'}=V$. Then $\mC'$ is a continued fraction expansion of $-p/q$ and
$U=\Xi V=B^{\mC}$ where $\mC=(a_{1},a_{2},\ldots,a_{n-1},0)$.
By \refcor{cor-lens-spaces}, (\ref{eq:anomaly}) and (\ref{eq:sumformula})
we therefore get
\begin{equation}\label{eq:RTlens}
\tau_{r}^{\frg}(L(p,q))= \omega^{\Phi(U)} \tilde{U}_{\rho\rho},
\end{equation}
where $\omega$ is given by (\ref{eq:omega}).
If $q=0$ we have $p=1$ and $L(p,q)=S^{3}$. In this case we have
$\tau_{r}^{\frg}(L(p,q))=\mD^{-1}$. We also have
$U=\Xi \Theta^{d}$ so by using
(\ref{eq:mR}), (\ref{eq:rank}) and (\ref{eq:omega})
we find that the right-hand side of (\ref{eq:RTlens})
is also equal to $\mD^{-1}$.
The identity (\ref{eq:RTlens}) coincides with \cite[Formula (3.7)]{Jeffrey2}
for $\frg=\frsl_{2}(\C)$, see also \cite[Formula (49)]{Hansen2}.

By elaborating on the expression (\ref{eq:RTlens})
along the same lines
as in \cite[Sect.~3]{Jeffrey2} we can now easily derive an explicit expression for
$\tau_{r}^{\frg}(L(p,q))$.
If $q \neq 0$ we get by \refthm{thm:main1} that
{\allowdisplaybreaks
\begin{eqnarray*}
\mR(U)_{\rho\rho} &=& \frac{\left(\I\sign(p)\right)^{\npr}}{(\kappa |p|)^{l/2} \vol (\cY)}
       \exp \left( -\frac{\pi \I}{\ch} \Phi (U) |\rho|^{2} \right) \\*
 && \hspace{.3in} \times \exp \left( \frac{\pi \I}{\kappa} \frac{b}{q} |\rho|^{2} \right) \sum_{w \in W} \det(w) \\*
 && \hspace{.3in} \times \sum_{\nu \in \cY/p\cY}
       \exp\left( \frac{\pi \I}{pq\kappa} |q\rho+\kappa q\nu-w(\rho)|^{2} \right) \\
 &=& \frac {\left(\I\sign(p)\right)^{\npr}}{(\kappa |p|)^{l/2} \vol (\cY)}
       \exp \left( -\frac{\pi \I}{\ch} \Phi (U)|\rho|^{2} \right) \\*
 && \hspace{.6in} \times \exp \left( \frac{\pi \I}{\kappa} \frac{b}{q} |\rho|^{2} \right) \\*
 && \hspace{.3in} \times \sum_{w \in W} \det (w) \exp\left( \frac{\pi \I}{pq\kappa} |q\rho-w(\rho) |^{2} \right) \\*
 && \hspace{.3in} \times \sum_{\nu \in \cY/p\cY} \exp\left( \frac{\pi \I q\kappa}{p} |\nu|^{2} \right) \\*
 && \hspace{.6in} \times \exp\left( \frac{2\pi \I}{p} \la \nu, q\rho-w(\rho) \ra \right).
\end{eqnarray*}}\noindent
By (\ref{eq:omega}) and (\ref{eq:RTlens}) we then get
{\allowdisplaybreaks
\begin{eqnarray*}
\tau_{r}^{\frg}(L(p,q)) &=& \frac {\left(\I\sign(p)\right)^{\npr}}{(\kappa |p|)^{l/2} \vol (\cY)}
       \exp \left( -\frac{\pi \I}{\kappa} \Phi (U)|\rho|^{2} \right) \\*
 && \hspace{.6in} \times \exp \left( \frac{\pi \I}{\kappa} \frac{b}{q} |\rho|^2 \right) \\
 && \hspace{.3in} \times \sum_{w \in W} \det (w) \exp\left( \frac{\pi \I}{pq\kappa} |q\rho-w(\rho) |^{2} \right) \\
 && \hspace{.3in} \times \sum_{\nu \in \cY/p\cY} \exp\left( \frac{\pi \I q\kappa}{p} |\nu|^{2} \right) \\*
 && \hspace{.6in} \times \exp\left( \frac{2\pi \I}{p} \la \nu, q\rho-w(\rho) \ra \right).
\end{eqnarray*}}\noindent
Here
{\allowdisplaybreaks
\begin{eqnarray*}
&& \exp\left( \frac{\pi \I}{\kappa} \frac{b}{q} |\rho|^{2} \right)
     \exp\left( \frac{\pi \I}{pq\kappa} |q\rho-w(\rho) |^{2} \right)\\
 && \hspace{.3in} =\exp \left( \frac{\pi \I}{pq\kappa}
       \left(pb |\rho |^{2} + q^{2} |\rho |^{2} + |\rho |^{2} - 2q \la \rho,w(\rho) \ra \right) \right) \\
 && \hspace{.3in} =\exp \left( \frac{\pi \I}{p\kappa}
       \left( d |\rho |^{2} + q |\rho |^{2} - 2 \la \rho,w(\rho) \ra \right) \right),
\end{eqnarray*}}\noindent
where we use that $pb+1=qd$. Moreover,
\begin{eqnarray*}
&& \exp \left( \frac{\pi \I}{p\kappa} \left( d |\rho |^{2} + q |\rho |^{2} - 2 \la \rho,w(\rho) \ra \right) \right)
      \exp \left( -\frac{\pi \I}{\kappa} \Phi (U) |\rho |^{2} \right) \\
 && \hspace{.3in} = \exp \left( -\frac{2\pi \I}{p\kappa} \la \rho,w(\rho) \ra \right)
       \exp \left( \frac{\pi \I}{\kappa} \left( \frac{d+q}{p} - \Phi(U) \right) |\rho |^{2} \right).
\end{eqnarray*}
Here
\begin{eqnarray*}
\Phi(U) &=& \frac{d+q}{p} - 12\sign(p)\s (d,p), \\
qd &\equiv& 1 \pmod{p}
\end{eqnarray*}
so
$$
\exp \left( \frac{\pi \I}{\kappa} \left( \frac{d+q}{p} - \Phi(U) \right) |\rho |^{2} \right)
 = \exp \left( \frac{\pi \I}{\kappa} 12 \sign(p)\s (q,p) |\rho |^{2} \right).
$$
By putting all the pieces together and using (\ref{eq:dedekindsymbol}) we get

\begin{thm}\label{Lie-lens}
The RT--invariants associated to $\frg$ of the lens space $L(p,q)$, $p \neq 0$,
are given by
\begin{eqnarray*}
&& \tau_{r}^{\frg}(L(p,q)) = \frac{\left(\I\sign(p)\right)^{\npr}}{(\kappa |p|)^{l/2} \vol(\cY)}
       \exp \left( \frac{\pi \I}{\kappa} S \left( \frac{q}{p} \right) |\rho|^{2} \right) \\
 && \hspace{1.0in} \times \sum_{w \in W} \det(w)
       \exp \left( -\frac{2\pi \I}{p\kappa} \la \rho, w(\rho) \ra \right) \\
 && \hspace{0.6in} \times \sum_{\nu \in \cY/p\cY}
       \exp \left( \pi \I \frac{q}{p} \kappa |\nu|^{2} \right)
       \exp \left( 2\pi \I \frac{1}{p} \la \nu, q\rho - w(\rho) \ra \right).
\end{eqnarray*}
We also have
\begin{eqnarray*}
&& \tau_{r}^{\frg}(L(p,q)) = \frac{\left(2\sign(p)\right)^{\npr}}{(\kappa |p|)^{l/2} \vol(\cY)}
       \exp \left( \frac{\pi \I}{\kappa} S \left( \frac{q}{p} \right) |\rho|^{2} \right) \\
 && \hspace{.1in} \times \sum_{\nu \in \cY/p\cY}
       \exp \left( \pi \I \frac{q}{p} \kappa |\nu|^{2} \right) \exp\left(2\pi\I \frac{q}{p}\la \nu , \rho \ra \right) \\
 && \hspace{.5in} \times \prod_{\alpha \in \Delta_{+}} \sin \left( \frac{\pi}{p\kappa} \la \rho + \kappa\nu,\alpha\ra\right).
\end{eqnarray*}\HS
\end{thm}

The second formula simply follows by using the Weyl denominator formula on
the first formula.
(A direct check shows, that the above theorem is also true for $q=0$, in which
case  $p=\pm 1$ and $L(p,q)=S^{3}$. If $p=0$, then $q=\pm 1$ and $L(p,q)=S^{1} \times S^{2}$,
and $\tau_{r}^{\frg}(S^{1}\times S^{2})=1$ in the normalization used here.)

\rk{The coprime case}
In this subsection we consider the coprime case, i.e.\ the case $(r,p)=1$, $r=m\kappa$.
In particular $p \neq 0$. From \refthm{Lie-lens} we have
$$
\tau_{r}^{\frg}(L(p,q)) = \frac{\left(\I\sign(p)\right)^{\npr}}{(\kappa |p|)^{l/2} \vol(\cY)}
       \exp \left( \frac{\pi \I}{\kappa} S \left( \frac{q}{p} \right) |\rho|^{2} \right)\Sigma,
$$
where
\begin{eqnarray*}
\Sigma &=& \sum_{w \in W} \det(w) \exp \left( -\frac {2\pi \I}{p\kappa} \la \rho,w(\rho) \ra \right) 
       \exp \left( -\pi \I \frac{q}{p} \kappa \left| \frac{q\rho-w(\rho)}{q\kappa} \right|^{2} \right) \\
 && \hspace{.3in} \times \sum_{\nu \in \cY/p\cY}
       \exp \left( \pi \I \frac{q}{p} \kappa \left|\nu + \frac{q\rho-w(\rho)}{q\kappa} \right|^{2} \right).
\end{eqnarray*}
First assume that $p$ is odd. Let $\rho_{w}=q\rho-w(\rho)$. Since $p$ and $4qr$ are coprime, there exist
integers $c$ and $a$ such that $pc+4qra=1$.
By definition of $\rho$, we have $\rho \in \frac{1}{2}Y$. Moreover, $mY \subseteq \cY$,
so $2m\rho_{w} \in \cY$. Therefore
{\allowdisplaybreaks
\begin{eqnarray*}
\Sigma &=& \sum_{w \in W} \det(w) \exp \left( -\frac{2\pi \I}{p\kappa} \la \rho,w(\rho) \ra \right) 
       \exp \left( -\frac{\pi \I}{pq\kappa} |\rho_{w}|^{2} \right) \\*
 && \hspace{.3in} \times \sum_{\nu \in \cY/p\cY}
       \exp \left( \pi \I \frac{q}{p} \kappa \left| \nu + \rho_{w} \left(4 a m + \frac{pc}{q\kappa} \right) \right|^{2} \right) \\
 &=& \sum_{w \in W} \det(w) \exp \left( -\frac{2\pi \I}{p\kappa} \la \rho,w(\rho) \ra \right) 
       \exp \left( -\frac{\pi \I}{pq\kappa} |\rho_{w}|^{2} \right) \\*
 && \hspace{.3in} \times \sum_{\nu \in \cY/p\cY}
       \exp \left( \pi \I \frac{q}{p} \kappa \left| \nu + \rho_{w} \frac{pc}{q\kappa} \right|^{2} \right) \\
 &=& \Sigma' \sum_{\nu \in \cY/p\cY} \exp \left( \pi \I \frac{q}{p} \kappa | \nu |^{2} \right), 
\end{eqnarray*}}\noindent
where
$$
\Sigma' = \sum_{w \in W} \det(w) \exp \left( -\frac{2\pi \I}{p\kappa} \la \rho,w(\rho) \ra \right)
       \exp \left( \frac{\pi \I}{pq\kappa} |\rho_{w}|^{2} ( p^{2} c^{2} - 1 ) \right).
$$
We used here (\ref{eq:integer1}) and the fact that $\rho_{w} \in X$.
Since $p^{2} c^{2}-1=4qra(4qra-2)$ we get
\begin{eqnarray*}
\Sigma' &=& \sum_{w \in W} \det(w)
       \exp \left( -\frac{2\pi \I}{p\kappa} \la \rho,w(\rho) \ra \right) \\*
 && \hspace{.3in} \times \exp \left( \frac{\pi \I}{p}
        4am(-2+4qr a) [ (q^{2}+1) |\rho|^{2}-2q\la \rho,w(\rho)\ra ] \right) \\
&=& \sum_{w \in W} \det(w) e^{\beta} e^{w},
\end{eqnarray*}
where
\begin{eqnarray*}
e^{\beta} &=& \exp \left( -\frac{4\pi \I}{p} am(1+pc) (q^{2}+1) |\rho|^{2} \right), \\
e^{w} &=& \exp \left( -\frac{2\pi \I}{p\kappa}
        \left[ 1+4qr a(4qr a -2) \right] \la \rho,w(\rho)\ra \right). \\
\end{eqnarray*}
Since $\rho \in X$ and $2m\rho \in \cY$ we have $2m|\rho|^{2} \in \Z$ so
\begin{eqnarray*}
e^{\beta} &=& \exp \left( -\frac{4\pi \I}{p} am  (q^{2}+1) |\rho|^{2} \right) \\
&=& \exp \left( -\frac{2\pi \I}{p} 4^{*}(q+q^{*})r^{*} 2m|\rho|^{2} \right),
\end{eqnarray*}
where $n^{*}$ is the inverse of $n \pmod{p}$ for any integer $n$
coprime to $p$. Moreover,
$$
e^{w} = \exp \left( -\frac{2\pi \I}{\kappa} p c^{2} \la \rho,w(\rho) \ra  \right)
 = \exp \left( -\frac{2\pi \I}{\kappa}c \la \rho,w(\rho) \ra \right),
$$
where we used $pc^{2}=c-4qrac$ and the fact that $2m\la \rho, w(\rho) \ra \in \Z$.
We thus obtain
{\allowdisplaybreaks
\begin{eqnarray*}
\tau_{r}^{\frg}(L(p,q)) &=& \frac{\left(\I\sign(p)\right)^{\npr}}{(\kappa |p|)^{l/2} \vol(\cY)}
       \exp \left( \frac{\pi \I}{\kappa} S \left( \frac{q}{p} \right) |\rho|^{2} \right) \\*
 && \hspace{.5in} \times \exp \left( -\frac{2\pi \I}{p} 4^{*}(q+q^{*})r^{*} 2m|\rho|^{2} \right) \\
 && \hspace{.2in} \times \sum_{w \in W} \det(w) \exp \left( -\frac{2\pi \I}{\kappa}c \la \rho,w(\rho) \ra \right) \\*
 && \hspace{.2in} \times \sum_{\nu \in \cY/p\cY} \exp \left( \pi \I \frac{q}{p} \kappa | \nu |^{2} \right).
\end{eqnarray*}}\noindent
Here
$$
\sum_{w \in W} \det(w) \exp \left( -\frac{2\pi \I}{\kappa}c \la \rho,w(\rho) \ra \right)
 = \prod_{\alpha \in \Delta_{+}} 2\I\sin \left( -\frac{\pi c}{\kappa} \la \rho,\alpha \ra \right)
$$
by the Weyl denominator formula.
Note that $c$ is the inverse of $p\pmod{4r}$.

Next we assume that $p$ is even.
Then $q$ is odd and there exist two integers $c$ and $a$ such that 
$4pc+qra=1$. We put $\rho_{w}=\frac{1}{2}(q\rho-w(\rho)) \in \frac{1}{2}Y$,
so $2m\rho_{w} \in \cY$. Moreover, $2\rho_{w} \in Y \subseteq X$. We therefore
find 
$$
\Sigma = \Sigma'' \sum_{\nu \in \cY/p\cY} \exp\left( \pi \I \frac{q}{p} \kappa |\nu|^{2} \right),
$$
where
{\allowdisplaybreaks
\begin{eqnarray*}
\Sigma'' &=& \sum_{w \in W} \det(w)
       \exp \left( -\frac{2\pi \I}{p\kappa} \la \rho,w(\rho) \ra \right) \\*
 && \hspace{.6in} \times \exp\left( \frac{\pi \I}{pq\kappa} (16p^{2} c^{2} -1) |2\rho_{w}|^{2} \right) \\
&=& \sum_{w \in W} \det(w)
       \exp \left( -\frac{2\pi \I}{p\kappa} \la \rho,w(\rho) \ra \right) \\*
 && \hspace{.3in} \times
       \exp \left( \frac{am\pi \I}{p} (qra-2)
          \left[ (q^2+1) |\rho|^{2} - 2q \la \rho, w(\rho)\ra \right] \right) \\
&=& \sum_{w \in W} \det(w)e^{\beta}e^{w},
\end{eqnarray*}}\noindent
where
\begin{eqnarray*}
e^{\beta} &=& \exp \left( \frac{am\pi\I}{p} (qra-2)(q^{2}+1)|\rho|^{2} \right), \\
e^{w} &=& \exp \left( \frac{2amq\pi\I}{p}(2-qra)\la \rho,w(\rho) \ra - \frac{2\pi\I}{p\kappa} \la \rho,w(\rho) \ra \right).
\end{eqnarray*}
Here $qra-2 = -1-4pc$ so
\begin{eqnarray*}
e^{\beta} &=& \exp \left( -\frac{2\pi\I}{4p} a (q^{2}+1)2m|\rho|^{2} \right) \\
 &=& \exp \left( -\frac{2\pi\I}{4p} (q+q^{*}) r^{*} 2m|\rho|^{2} \right),
\end{eqnarray*}
where $n^{*}$ is the inverse of $n \pmod{4p}$ for any integer $n$ coprime to $p$. Moreover,
$$
e^{w} = \exp \left( \frac{2\pi \I}{p\kappa} (aqr-1) \la \rho, w(\rho)\ra \right)
   = \exp \left( -\frac{2\pi \I}{\kappa} 4c \la \rho, w(\rho)\ra \right).
$$
We have thus shown

\begin{prop}\label{Lie-lens-coprime}
Let $r=m\kappa$ be coprime to $p$. If $p$ is odd we have
\begin{eqnarray*}
\tau_{r}^{\frg}(L(p,q)) &=& \frac{\left(\I\sign(p)\right)^{\npr}}{(\kappa |p|)^{l/2} \vol(\cY)}
       \exp \left( \frac{\pi \I}{\kappa} S \left( \frac{q}{p} \right) |\rho|^{2} \right) \\
 && \hspace{.5in} \times \exp \left( -\frac{2\pi \I}{p} 4^{*}(q+q^{*})r^{*} 2m|\rho|^{2} \right) \\
 && \hspace{.2in} \times \sum_{w \in W} \det(w) \exp \left( -\frac{2\pi \I}{\kappa}c \la \rho,w(\rho) \ra \right) \\
 && \hspace{.2in} \times \sum_{\nu \in \cY/p\cY} \exp \left( \pi \I \frac{q}{p} \kappa | \nu |^{2} \right).
\end{eqnarray*}
Alternatively we have
\begin{eqnarray*}
\tau_{r}^{\frg}(L(p,q)) &=& \frac{\left(2\sign(p)\right)^{\npr}}{(\kappa |p|)^{l/2} \vol(\cY)}
       \exp \left( \frac{\pi \I}{\kappa} S \left( \frac{q}{p} \right) |\rho|^{2} \right) \\
 && \hspace{.1in} \times \exp \left( -\frac{2\pi \I}{p} 4^{*}(q+q^{*})r^{*} 2m|\rho|^{2} \right) \\
 && \hspace{.1in} \times \left(\prod_{\alpha \in \Delta_{+}} \sin \left( \frac{\pi c}{\kappa} \la \rho,\alpha \ra \right)\right)
       \sum_{\nu \in \cY/p\cY} \exp \left( \pi \I \frac{q}{p} \kappa | \nu |^{2} \right).
\end{eqnarray*}
Here $n^{*}$ is the inverse of $n\pmod{p}$ for any integer $n$ coprime to $p$, and
$c$ is the inverse of $p\pmod{4r}$.

If $p$ is even we have
\begin{eqnarray*}
\tau_{r}^{\frg}(L(p,q)) &=& \frac{\left(\I\sign(p)\right)^{\npr}}{(\kappa |p|)^{l/2} \vol(\cY)}
       \exp \left( \frac{\pi \I}{\kappa} S \left( \frac{q}{p} \right) |\rho|^{2} \right) \\
 && \hspace{.5in} \times \exp \left( -\frac{2\pi\I}{4p} (q+q^{*}) r^{*} 2m|\rho|^{2} \right) \\
 && \hspace{.2in} \times \sum_{w \in W} \det(w) \exp \left( -\frac{2\pi \I}{\kappa} 4c \la \rho, w(\rho)\ra \right) \\
 && \hspace{.2in} \times \sum_{\nu \in \cY/p\cY} \exp \left( \pi \I \frac{q}{p} \kappa | \nu |^{2} \right).
\end{eqnarray*}
Alternatively we have
\begin{eqnarray*}
\tau_{r}^{\frg}(L(p,q)) &=& \frac{\left(2\sign(p)\right)^{\npr}}{(\kappa |p|)^{l/2} \vol(\cY)}
       \exp \left( \frac{\pi \I}{\kappa} S \left( \frac{q}{p} \right) |\rho|^{2} \right) \\
 && \hspace{.1in} \times \exp \left( -\frac{2\pi\I}{4p} (q+q^{*}) r^{*} 2m|\rho|^{2} \right) \\
 && \hspace{.1in} \times \left(\prod_{\alpha \in \Delta_{+}} \sin \left( \frac{4\pi c}{\kappa} \la \rho,\alpha \ra \right) \right)
       \sum_{\nu \in \cY/p\cY} \exp \left( \pi \I \frac{q}{p} \kappa | \nu |^{2} \right).
\end{eqnarray*}
Here $n^{*}$ is the inverse of $n\pmod{4p}$ for any integer $n$ coprime to $p$, and $4c$ is the
inverse of $p\pmod{r}$.\HS
\end{prop}

\begin{rem}
a) By a direct check using \refthm{Lie-lens} one finds that $L(64,9)$ and $L(64,25)$ are distinguished by
the $\frsl_{4}(\C)$--invariant, in fact by $\tau_{6}^{\frsl_{4}(\C)}$.
On the other hand these lens spaces can not be distinguished
by the LMO invariant. This follows by \cite[Proposition 5.1]{Bar-NatanLawrence}
and the fact that $S(9/64)=S(25/64)(=-63/32)$, cf.\ \cite[Remark (4.14)]{KirbyMelvin}.
Since the perturbative invariant $\tau^{P\frsl_{n}(\C)}$ \cite{Le3}
can be recovered from the LMO invariant, cf.\ \cite{Ohtsuki}, we see
that $\tau^{P\frsl_{n}(\C)}$ does not separate $L(64,9)$ and $L(64,25)$.
The perturbative invariant $\tau^{P\frsl_{n}(\C)}$ is
determined by the family of quantum $P\frsl_{n}(\C)$--invariants
$\tau_{r}^{P\frsl_{n}(\C)}$, see \cite{Le3}. Oppositely it is expected
that the perturbative invariant $\tau^{P\frsl_{n}(\C)}$
dominates the invariants $\tau_{r}^{P\frsl_{n}(\C)}$, cf.\
\cite[Conjecture 1.8]{Le3}.
From the explicit formulas for
the $P\frsl_{4}$--invariants $\tau_{r}^{P\frsl_{4}(\C)}$ of the lens spaces
in \cite{Takata}, it follows
that these invariants can not distinguish between $L(64,9)$ and
$L(64,25)$ for $r$ a prime. It would be interesting to examine
if this is also the case for the non-prime $r$ for which $\tau_{r}^{P\frsl_{4}(\C)}$
is defined.

b) In a forthcoming paper \cite{HansenTakata2} we make detailed calculations
of the Gaussian sums $\sum_{\nu \in \cY/p\cY} \exp \left( \pi \I \frac{q}{p} \kappa | \nu |^{2} \right)$,
thereby obtaining more detailed separation results.
\end{rem}

\rk{The asymptotic expansion conjecture and lens spaces}
In this section we calculate the large $r$ asymptotics of
$r \mapsto \tau_{r}^{\frg}(L(p,q))$, $r=m\kappa$.
Let us begin by some introductory remarks.
Therefore, let $X$ be an arbitrary closed oriented $3$--manifold and
let $G$ be a simply connected compact simple Lie group with complexified
Lie algebra $\frg$. We are interested in the behaviour of the complex function
$\kappa \mapsto \tau_{r}^{\frg}(X)$ in
the limit of large $\kappa$, i.e.\ $\kappa \rightarrow \infty$.
It is believed that Witten's TQFT associated with $G$ and $k$
coincides with the TQFT of Reshetikhin and Turaev associated with $\frg$ and $k+h^{\vee}$.
In particular it is conjectured that Witten's semiclassical approximation for 
$Z_{k}^{G}(X)$
should be valid for the function $\kappa \mapsto \tau_{r}^{\frg}(X)$,
and furthermore that this function should have a full 
asymptotic expansion  in the limit $\kappa \ria \infty$.
The precise formulation of this is stated in the following conjecture,
called the asymptotic expansion conjecture (AEC).

\begin{conj}[J.\ E.\ Andersen \cite{Andersen1}, \cite{Andersen2}]\label{AEC} 
Let $\{\alpha_{1},\ldots,\alpha_{M} \}$
be the set
of values of the Chern--Simons functional of flat $G$ connections on 
a closed oriented $3$--manifold $X$. Then there exist $d_{j} \in \Q$,
$\tilde{I}_{j} \in \Q / \Z$, $b_{j} \in \R_{+}$ and $c_{n}^{(j)} \in \C$
for $j=1,\ldots,M$ and $n=1,2,3,\ldots$ such that we for all $N=0,1,2,\ldots$ have
\begin{equation}\label{eq:2}
\tau_{m\kappa}^{\frg}(X) = \sum_{j=1}^{M} b_{j}e^{2\pi \sqrt{-1} \kappa \alpha_{j}} \kappa^{d_{j}} e^{\pi \sqrt{-1} \tilde{I}_{j}/4} \left(1+\sum_{n=1}^{N} c_{n}^{(j)} \kappa^{-n} \right) + o(\kappa^{d-N})
\end{equation}
in the limit $\kappa \ria \infty$, where $d=\max\{d_{1},\ldots,d_{M} \}$.
\end{conj}

\noindent Here $f(\kappa)=o(\kappa^{d-N})$ means as usual that
$f(\kappa)/\kappa^{d-N} \ria 0$ as $\kappa \ria \infty$.
The AEC was proposed by Andersen in \cite{Andersen1}, where he proved it
for the mapping tori of finite
order diffeomorphisms of orientable surfaces
of genus at least $2$ and for general $\frg$
using the gauge theoretic approach to the quantum invariants.
These manifolds are Seifert
manifolds with orientable base
and Seifert Euler number
equal to zero, see \cite[Sect.~4]{Andersen1}.
Note that the semiclassical approximation is given by putting $N=0$
in the sum $\sum_{j=1}^{M}...$ in (\ref{eq:2}).
This part of the AEC and some of the conjectures concerning the
topological interpretation of the different parts of the asymptotic formula,
see \cite{Andersen2} for details, are in fact inspired by the works of Witten,
Freed and Gompf, Jeffrey and Rozansky
on the semiclassical approximation of Witten's invariants.

As already stated, Jeffrey \cite{Jeffrey1}, \cite{Jeffrey2}
and Garoufalidis \cite{Garoufalidis} made completely rigorous
calculations of the semiclassical approximation of the $SU(2)$--invariants
of lens spaces. Actually these calculations contain a complete
verification of the AEC
for the lens spaces and $\frg=\frsl_{2}(\C)$.
Jeffrey's calculations also contain a proof of the AEC for
a certain class of mapping tori over the torus 
for an arbitrary complex finite dimensional simple Lie algebra.

In \cite{Rozansky2}, \cite{Rozansky3}
L.\ Rozansky calculated the Witten $SU(2)$--invariants
of all Seifert manifolds with orientable base
and carried through a rather technical
analysis leading to a candidate for
the full asymptotic expansion of these
invariants (for the Seifert manifolds with Seifert Euler number different from zero).
As shown in \cite[Sect.~8]{Hansen2} the invariants calculated
by Rozansky are equal
to the RT--invariants associated to $\frsl_{2}(\C)$.
However, to actually prove that his formula gives the asymptotic
expansion of the invariants one has to incorporate estimations of error terms in
the calculations. This was left
out in Rozansky's calculations. 
In \cite{Hansen1}, \cite{Hansen3} the first author has supplemented
the calculations of Rozansky by making the necessary error estimates
thereby proving, that Rozansky's formula is really the asymptotic expansion of these invariants. 
In \cite{Hansen3} the calculations are carried through for a big class
of functions including
the $\frsl_{2}(\C)$--invariants of all oriented Seifert manifolds with
orientable base or non-orientable base with even genus (also the
ones with Seifert Euler number equal to zero).
Based on results of D.\ Auckly \cite{Auckly}
the Chern--Simons invariants can be identified
in the asymptotic formula thereby proving
the AEC, \refconj{AEC}, for these Seifert manifolds and $\frg=\frsl_{2}(\C)$.
It should be mentioned that Rozansky and Lawrence have
calculated the asymptotics of the $\frsl_{2}(\C)$--invariants
of a subclass of the Seifert manifolds with base $S^{2}$ by a method
which avoids the rather technical analysis of error terms,
cf.\ \cite{LawrenceRozansky}. However, it seems that their method does not
work for arbitrary oriented Seifert manifolds.

From \refthm{Lie-lens}
it is obvious, that the large $\kappa$ asymptotics
of $\tau_{r}^{\frg}(L(p,q))$ is on a form as predicted by
the asymptotic expansion conjecture, \refconj{AEC}. We
note that $L(p,q)$ and $L(p',q')$ are
homeomorphic if and only if $p=p'$ and
$$
q \equiv \pm q' \pmod{p} \hspace{.1in}\text{or}\hspace{.1in} qq' \equiv \pm 1 \pmod{p}.
$$
A homeomorphism preserves orientation if and only if the relevant sign
is $+$. If $q^{*}$ denotes the inverse of $q$ in the group of units of
$Z/pZ$ we therefore have that there is an orientation preserving
homeomorphism between $L(p,q)$ and $L(p,q^{*})$. In particular,
we can exchange $q$ and $q^{*}$ in any of the formulas
for $\tau_{r}^{\frg}(L(p,q))$. For the following discussion this seems to
be an advantage.
We have
\begin{eqnarray*}
\sin\left( \frac{\pi}{p\kappa} \la \rho + \kappa\nu,\alpha\ra\right) &=& 
  \sin\left(\frac{\pi}{p\kappa}\la\rho,\alpha\ra\right)\cos\left(\frac{\pi}{p}\la\nu,\alpha\ra\right) \\ 
 &&\hspace{.2in} + \cos\left(\frac{\pi}{p\kappa}\la\rho,\alpha\ra\right)\sin\left(\frac{\pi}{p}\la\nu,\alpha\ra\right).
\end{eqnarray*}
If we let
$$
\mM_{j}=\{ \; \nu \in \cY/p\cY \; | \; \la \nu , \alpha \ra \in p\Z \hspace{.1in}\text{for exactly}\hspace{.1in} j \hspace{.1in}\text{elements}\hspace{.1in} \alpha \hspace{.1in}\text{in}\hspace{.1in} \Delta_{+}\;\}
$$
for $j=0,1,\ldots,\npr$, then we get

\begin{cor}\label{lens-asymp}
There exists a family of complex numbers $c_{n}^{(\nu)}$, $n=1,2,\ldots$, $\nu \in \cY/p\cY$, depending
directly on $p$ but only on $q$ through $S(q/p)$, such that
\begin{eqnarray*}
&&\tau_{r}^{\frg}(L(p,q)) = \frac{\left(2\sign(p)\right)^{\npr}}{(\kappa |p|)^{l/2} \vol(\cY)}
       \sum_{j=0}^{\npr} \left( \frac{\pi}{p\kappa} \right)^{j}
       \sum_{\nu \in \mM_{j}} b_{\nu} \exp \left( 2 \pi \I \kappa \frac{q^{*}}{2p} |\nu|^{2} \right) \\ 
 && \hspace{1.1in} \times \exp\left(2\pi\I \frac{q^{*}}{p}\la \nu , \rho \ra \right)
            \left(1+\sum_{n=1}^{\infty} c_{n}^{(\nu)}\kappa^{-n}\right),
\end{eqnarray*}
for all $\kappa \in \Z_{\geq h^{\vee}}$, where
$$
b_{\nu} = \prod_{\alpha \in \Delta_{+} \; : \; \la \nu , \alpha \ra \in p\Z } (-1)^{ \la \nu,\alpha \ra /p } \la \rho,\alpha \ra 
          \prod_{\alpha \in \Delta_{+} \; : \; \la \nu , \alpha \ra \notin p\Z } \sin \left( \frac{\pi}{p} \la \nu,\alpha \ra \right).
$$\HS
\end{cor}

Note here that $S(q/p)=S(q^{*}/p)$.
The infinite power series in $1/\kappa$ present in the
above corollary are convergent for all $\kappa \in \Z_{\geq h^{\vee}}$.
Note also that $0 \in \mM_{\npr}$, so this set is non-empty.
Let us in some greater detail look at a few examples.
First assume that $\frg=\frsl_{2}(\C)$. In this case we have
that $\Delta_{+}$ contains one element $\alpha$ of length $\sqrt{2}$ and
$\cY=Y=\Span_{\Z}\{\alpha\}$ so $\vol(\cY)=\sqrt{2}$.
For $n \in \{0,1,\ldots,|p|-1\}$ we have
$$
\la n\alpha,\alpha\ra =2n
$$
so if we identify $n\alpha$ by $n$ we have
$$
\mM_{1} = \left\{ \begin{array}{ll}
 \{ 0 \} & ,\hspace{.2in}\text{if}\hspace{.1in}p\hspace{.1in}\text{is odd} \\
 \{0,|p|/2\} & ,\hspace{.2in}\text{if}\hspace{.1in}p\hspace{.1in}\text{is even}
\end{array}\right.
$$
and $\mM_{0} = \{0,1,\ldots,|p|-1\} \sm \mM_{1}$.
For $p$ odd we therefore have
\begin{eqnarray*}
&&\tau_{r}^{\frsl_{2}(\C)}(L(p,q)) = \frac{\pi}{|p|\kappa}\sqrt{\frac{2}{|p|\kappa}} \left( 1 + \sum_{l=1}^{\infty} c_{l}^{(0)} \kappa^{-l} \right) \\
 && \hspace{.3in} +  \sign(p) \sqrt{\frac{2}{|p|\kappa}}\sum_{n=1}^{|p|-1} \exp\left( 2\pi\I\kappa\frac{q^{*}}{p} n^{2} \right) 
            \exp\left( 2\pi\I\frac{q^{*}}{p}n\right) \\
 && \hspace{1.1in} \times \sin\left( \frac{2\pi n}{p} \right)\left( 1 + \sum_{l=1}^{\infty} c_{l}^{(n)} \kappa^{-l} \right).
\end{eqnarray*}
For $p$ even we find that
\begin{eqnarray*}
&&\tau_{r}^{\frsl_{2}(\C)}(L(p,q)) = \frac{\pi}{|p|\kappa}\sqrt{\frac{2}{|p|\kappa}}
               \left( 1 - (-1)^{q^{*}p/2} + \sum_{l=1}^{\infty} \left(c_{l}^{(0)}+c_{l}^{|p|/2)} \right)\kappa^{-l} \right) \\
 && \hspace{.3in} +  \sign(p) \sqrt{\frac{2}{|p|\kappa}}\sum_{n \in \{1,2,\ldots,|p|-1\}\; : \; n \neq \frac{|p|}{2}} 
            \exp\left( 2\pi\I\kappa\frac{q^{*}}{p} n^{2} \right) \\
 && \hspace{0.6in} \times \exp\left( 2\pi\I\frac{q^{*}}{p}n\right) \sin\left( \frac{2\pi n}{p} \right)\left( 1 + \sum_{l=1}^{\infty} c_{l}^{(n)} \kappa^{-l} \right).
\end{eqnarray*}
For all $p$ we find that the leading large $\kappa$ asymptotics
of $\tau_{r}^{\frsl_{2}(\C)}(L(p,q))$ is
$$
\I\sign(p) \sqrt{\frac{2}{|p|\kappa}}\sum_{n=1}^{|p|-1}
            \exp\left( 2\pi\I\kappa\frac{q^{*}}{p} n^{2} \right)
            \sin\left( \frac{2\pi q^{*}n}{p}\right) \sin\left( \frac{2\pi n}{p} \right),
$$
where we use that
$$
\sum_{n=1}^{|p|-1}
            \exp\left( 2\pi\I\kappa\frac{q^{*}}{p} n^{2} \right)
            \cos\left( \frac{2\pi q^{*}n}{p}\right) \sin\left( \frac{2\pi n}{p} \right)=0.
$$
This result coincides with \cite[Formula (5.7)]{Jeffrey2}.
It is well-known, see e.g.\ \cite[Formula (5.3)]{Jeffrey2},
that the set of values of the Chern--Simons
functional of flat $SU(2)$ connections on $L(p,q)$ is
\footnote{There seems to be a problem with signs here. The set $\mS$ of
Chern--Simons invariants has been calculated by Kirk and Klassen,
cf.\ \cite[Theorem 5.1]{KirkKlassen}. According to their result
all the above stated Chern--Simons values have to be multiplied by $-1$. Note that $L(p,q)$ in
\cite{KirkKlassen} is equal to $L(p,-q)$ here.}
$$
\mS=\left\{ \frac{q^{*}}{p}n^{2} \pmod{1} \; | \; n=0,1,\ldots,|p|-1 \right\},
$$
so here we see the reason for replacing $q$ by $q^{*}$.
(Note that $q^{*}(p-n)^{2}/p \equiv q^{*}n^{2}/p \pmod{1}$.)
If we e.g.\ have $p=k^{2}$, where $k$ is a
positive integer, then
$$
\frac{q^{*}}{p}n^{2} \equiv 0 \pmod{1}
$$
for $n=k$, so we see that two elements belonging to different
sets $\mM_{j}$ can have the same Chern--Simons value.

Let us also examine the type $A_{2}$. Here
$\cY=Y=\Span_{\Z}\{ \alpha_{1},\alpha_{2} \}$, where $\alpha_{i}$
have length $\sqrt{2}$, $i=1,2$.
We have $\Delta_{+} = \{\alpha_{1},\alpha_{2},\alpha_{1}+\alpha_{2}\}$,
and in particular $\npr=3$.
If we identify $\nu=k\alpha_{1}+n\alpha_{2} \in \cY/p\cY$
with $(k,n) \in \{0,1,\ldots,|p|-1\}^{2}$ one finds by an
elementary analysis that
\begin{eqnarray*}
\mM_{3} &=& \left\{ \begin{array}{ll}
 \{ (0,0),(l,2l),(2l,l) \} & ,\hspace{.2in}\text{if}\hspace{.1in}|p|=3l, \;l \in \Z \sm \{0\}, \\
 \{ (0,0) \} & ,\hspace{.2in}\text{otherwise},
\end{array}\right.\\
\mM_{2} &=& \emptyset
\end{eqnarray*}
and $\mM_{0}=\{0,1,\ldots,|p|-1\}^{2} \sm (\mM_{1} \cup \mM_{3})$,
where $\mM_{1}=\emptyset$ for $|p|=3$ and
$$
\mM_{1} =\left\{ (k,n) \in \{0,1,\ldots,|p|-1\}^{2} \; \left| \begin{array}{l}
 \;\text{exactly one of the identities} \\
 \;k=2n, \hspace{.1in}n=2k, \\
 \;k=2n-|p|, \hspace{.1in}n=2k-|p| \\
 \;k+n=|p| \\
 \;\text{is satisfied}
\end{array}\right.\right\}
$$
otherwise. For $|p|=1$ we have $\mM_{3}=\{(0,0)\}$ and $\mM_{j}=\emptyset$, $j=0,1,2$,
and for $|p|=2$ we have $\mM_{3}=\{(0,0)\}$, $\mM_{1}=\{(0,1),(1,0),(1,1)\}$ and
$\mM_{0}=\mM_{2}=\emptyset$.
\refcor{lens-asymp} leads together with \refconj{AEC}
immediately to the following conjecture (use the
uniqueness property of asymptotic expansions of the
form (\ref{eq:2}), i.e.\ the fact that an arbitrary 
function $f\co \Z_{>0} \to \C$ has at
most one asymptotic expansion of the form (\ref{eq:2}) if the
$\alpha_{j}$'s are mutually different and rational).

\begin{conj}\label{Chern--Simons}
The set of values of the Chern--Simons functional of flat $G$ connections on
$L(p,q)$ is given by\footnote{perhaps with the opposite signs dependend on
the choice of conventions}
$$
\left. \left\{ \; \frac{q^{*}}{2p} |\nu|^{2} \pmod{\Z} \;\; \right| \;\; \nu \in \cY/p\cY \;\right\},
$$
for any simply connected, compact simple Lie group $G$. 
\end{conj}

Proving this conjecture will together with \refcor{lens-asymp} give a proof of the AEC
for the invariants $\tau_{r}^{\frg}(L(p,q))$.
In fact, by the uniqueness property of asymptotic expansions of
the form (\ref{eq:2}), \refconj{Chern--Simons}
should for $G=SU(n)$ be a corollary of \refcor{lens-asymp} and
the recent result of Andersen, mentioned in the introduction.
For $G=SU(n)$ \refconj{Chern--Simons} should also follow from results in \cite{Nishi}.

\section{A rational surgery formula for the invariant $\tau_{r}^{\frg}$}\label{sec-A-rational}

In this final section we derive a rational surgery formula for the invariant $\tau_{r}^{\frg}$.
The result follows easily from the surgery formula derived in
\cite{Hansen2} in the setting of a general modular tensor category.
By rational surgery we mean rational surgery on an arbitrary closed oriented
$3$--manifold $M$ along a framed link inside $M$.
Let us briefly recall the situation from \cite{Hansen2}.
Let $\left( \mV,\{V_{i}\}_{i \in I} \right)$
be a fixed modular category with a fixed rank $\mD$,
and let $\tau$ be the RT--invariant associated to these data.
For a link $L$ we let $\col(L)$ be the set of mappings from the set of components of $L$
to the index set $I$.
Before giving the result in the general case, let us first consider
rational surgery along links in $S^{3}$.
If $L \subseteq S^{3}$ is framed and oriented
and $\lambda \in \col(L)$ we let $\Gamma(L,\lambda)$
be the colored ribbon graph induced by $L$ with the $i$'th component $L_{i}$
of $L$ colored by $V_{\lambda(L_{i})}$.

\begin{thm}[\cite{Hansen2}]\label{surgery-formula3sphere}
Let $L$ be a link in $S^{3}$ with $n$ components and let $M$ be the
$3$--manifold given by surgery on $S^{3}$ along $L$ with
surgery coefficient $p_{i}/q_{i} \in \Q$ attached to the $i$'th
component, $i=1,2,\ldots,n$ {\em (}so we assume
$q_{i} \neq 0$, $i=1,2,\ldots,n$, see the comments to {\em (\ref{eq:b4}))}.
Moreover, let $\Omega$ be a
colored ribbon graph in $M$ {\em (}also identified with a colored
ribbon graph in $S^{3} \sm L${\em )}.
Let $L_{0}$ be $L$ considered as a framed link
with all components given the framing $0$ and with an arbitrary chosen but
fixed orientation. Finally, let
$\mC_{i}=(a_{1}^{(i)},\ldots,a_{m_{i}}^{(i)})$ be a
continued fraction expansion of $p_{i}/q_{i}$, $i=1,2,\ldots,n$.
Then
\begin{eqnarray*}
\tau(M,\Omega) &=& (\Delta \mD^{-1})^{\sigma+\sum_{i=1}^{n} c_{i}}\mD^{-\sum_{i=1}^{n} m_{i}} \\
 && \hspace{.4in} \times \sum_{\lambda \in \col(L)} \tau(S^{3},\Gamma(L_{0},\lambda) \cup \Omega) \left( \prod_{i=1}^{n} G^{\mC_{i}}_{\lambda(L_{i})0} \right),
\end{eqnarray*}
where
$c_{i}=\frac{1}{3}\left( \sum_{j=1}^{m_{i}} a_{j}^{(i)} - \Phi(B^{\mC_{i}}) \right)$,
$i=1,\ldots,n$, and $\sigma$ is the signature of the linking matrix of $L$
{\em (}with the surgery coefficients $p_{1}/q_{1},\ldots,p_{n}/q_{n}$ on the
diagonal {\em )}.\HS
\end{thm}

In the case of surgery on arbitrary closed oriented $3$--manifolds
along framed links we do not have a preferred
framing as in the case of surgery on $S^{3}$ (or on another integral
homology sphere), i.e.\ we can not identify the framing of a
link component with an integer in a canonical way,
see \cite[Appendix B]{Hansen2}.
Here, by a framed link in a closed oriented $3$--manifold $M$, we mean
a pair $(L,Q)$, where
$Q=\amalg_{i=1}^{n} Q_{i} \co \amalg_{i=1}^{n} (B^{2} \times S^{1}) \to M$
is an embedding (or more precisely an isotopy class of such embeddings)
and $L$ is the image by $Q$ of $\amalg_{i=1}^{n} (0 \times S^{1})$.
For other definitions of framed links in $3$--manifolds and how these
relate to this definition we refer to \cite[Appendix B]{Hansen2}.
The following result is
sensitive to a choice of orientations. We will use the following conventions.

\begin{conv}\label{conventions}
The space $B^{2} \times S^{1}$ is the standard
solid torus in $\R^{3}$ with the orientation induced by the standard right-handed
orientation of $\R^{3}$. Here $S^{1}$
is the standard unit circle in the $xz$--plane with centre $0$ and oriented
counterclockwise, i.e., $e_{3}$ is a positively oriented tangent vector in
the tangent space $T_{e_{1}}S^{1} \subseteq \R^{3}$, $e_{i}$ being
the $i$'th standard unit vector in $\R^{3}$, see Fig.~\ref{fig-torus}.
For a framed link $(L,Q)$ as above we will always assume that
each copy of $B^{2} \times S^{1}$ is this oriented standard solid torus,
and that $Q$ is orientation preserving after giving the
image of $Q$ the orientation induced by that of $M$ (we
can always obtain this by composing some of the $Q_{i}$ by
$g \times \id_{S^{1}}$ if necessarily, where $g \co B^{2} \to B^{2}$ is an orientation
reversing homeomorphism).
Moreover, we orient $L$ so that $Q_{i}$ restricted to
$S^{1} \times \{0\}$ is orientation preserving for each $i$.
The oriented meridian $\alpha$ and longitude $\beta$,
see Fig.~\ref{fig-torus},
represent a basis (over $\Lambda$) of $H_{1}(\Sigma_{(1;)};\Lambda) = \Lambda \oplus \Lambda$,
$\Lambda=\Z,\R$, $\Sigma_{(1;)}=S^{1} \times S^{1}$.
(For the notation $\Sigma_{(1;)}$, see \cite[Chap.~IV]{Turaev}.)
\end{conv}

\begin{figure}[h]

\begin{center}
\begin{texdraw}
\drawdim{cm}

\setunitscale 0.7

\linewd 0.02

\move(0 -0.7) \lellip rx:1 ry:1.25
\move(0 -0.7) \lellip rx:2.6 ry:3.1

\move(-2.75 -0.55) \lvec(-2.6 -0.85)
\move(-2.45 -0.55) \lvec(-2.6 -0.85)

\move(-3.2 -0.85) \htext{$\beta$}

\move(2 -0.9) \lvec(1.7 -0.75)
\move(2 -0.9) \lvec(1.7 -1.05)

\move(1.8 -1.35) \htext{$\alpha$}

\move(5 -1.5) \lvec(7 -1.5)
\lvec(6.7 -1.35)

\move(7 -1.5) \lvec(6.7 -1.65)

\move(5 -1.5) \lvec(5 0.5)
\lvec(4.85 0.2)

\move(5 0.5) \lvec(5.15 0.2)

\move(6.7 -2) \htext{$x$}
\move(5.2 0.2) \htext{$z$}

\move(1 -0.6) \clvec(1 -1)(2.6 -1)(2.6 -0.6)
\lpatt(0.067 0.1) \move(1 -0.6) \clvec(1 -0.2)(2.6 -0.2)(2.6 -0.6)

\end{texdraw}
\end{center}

\nocolon
\caption{}\label{fig-torus}
\end{figure}

Let us recall the notion of rational surgery on $M$ along $(L,Q)$.
Therefore, let $U_{i}=Q_{i}(B^{2} \times S^{1})$
and let $l_{i}=Q_{i}(e_{1} \times S^{1})$ oriented so that
$[l_{i}]=[L_{i}]$ in $H_{1}(U_{i};\Z)$, where
$L_{i}=Q_{i}(0 \times S^{1})$. Moreover, let
$\mu_{i}=Q_{i}(\partial B^{2} \times 1)$ oriented so that
$(\partial Q_{i})_{*}([\alpha])=[\mu_{i}]$ in $H_{1}(\partial U_{i};\Z)$,
where $\partial Q_{i}$ is the restriction of $Q_{i}$
to $\partial B^{2} \times S^{1}=\Sigma_{(1;)}$.
Let $(p_{i},q_{i})$ be pairs of coprime integers,
let $h_{i} \co \partial U_{i} \to \partial U_{i}$ be homeomorphisms such that
\begin{equation}\label{eq:b4}
(h_{i})_{*}([\mu_{i}])= \pm (p_{i}[\mu_{i}]+q_{i}[l_{i}])
\end{equation}
in $H_{1}(\partial U_{i};\Z)$, let $h$ be the union of the $h_{i}$,
and let $U=\amalg_{i=1}^{n} U_{i}$ be the image of $Q$.
Then the $3$--manifold $M' = (M \sm \interior (U)) \cup_{h} U$
is said to be the result of doing surgery on $M$ along the framed link $(L,Q)$ with
surgery coefficients $\{ p_{i}/q_{i} \}_{i=1}^{n}$. If $q_{i}=0$ so
$p_{i}=\pm 1$ we just write $\infty$ for $p_{i}/q_{i}$. Such surgeries
do not change the manifold
(up to an orientation preserving homeomorphism).
If, in (\ref{eq:b4}), $p_{i}=0$ and $q_{i}=\pm 1$ for all $i$, i.e.\ all
surgery coefficients are $0$, then we call $M'$
the result of doing surgery on $M$ along the framed link $(L,Q)$.
We equip $M'$ with the unique orientation extending the
orientation in $M \sm \interior (U)$. The above generalizes
ordinary rational surgery along links in $S^{3}$.
We call a homeomorphism $h$ satisfying (\ref{eq:b4}) an attaching map for the
surgery. We can and will always choose an orientation preserving attaching
map. Up to an orientation preserving homeomorphism the result of doing surgery
on $M$ along the framed link $(L,Q)$ with surgery coefficients
$\{ p_{i}/q_{i} \}_{i=1}^{n}$ is well defined, independent of the choices
of representative $Q$ and attaching map $h$.

For $\lambda \in \col(L)$ we let 
$\Gamma(L,\lambda)=\cup_{i=1}^{n} \Gamma(L_{i},\lambda(L_{i}))$,
where $\Gamma(L_{i},j)$ is
the colored ribbon graph equal to the
directed annulus $Q_{i}(([-1/2,1/2]\times 0) \times S^{1})$ with oriented
core $L_{i}$ and color $V_{j}$, $j \in I$.

\begin{thm}[\cite{Hansen2}]\label{surgery-formula}
Let $\mC_{i}=(a_{1}^{(i)},\ldots,a_{m_{i}}^{(i)}) \in \Z^{m_{i}}$ be
a continued fraction expansion of $p_{i}/q_{i}$, $i=1,\ldots,n$.
Moreover, let $\Omega$ be a
colored ribbon graph in $M'$ {\em (}also identified with a colored
ribbon graph in $M \sm L${\em )}.
Then
\begin{eqnarray*}
\tau(M',\Omega) &=& (\Delta \mD^{-1})^{\mu+\sum_{i=1}^{n} c_{i}}\mD^{-\sum_{i=1}^{n} m_{i}} \\
 & & \hspace{.4in} \times \sum_{\lambda \in \col(L)} \tau(M,\Gamma(L,\lambda) \cup \Omega) \left( \prod_{i=1}^{n} G^{\mC_{i}}_{\lambda(L_{i})0} \right),
\end{eqnarray*}
where $\mu$ is a sum of signs given by {\em (\ref{eq:e7})} and
$c_{i}=\frac{1}{3}\left( \sum_{j=1}^{m_{i}} a_{j}^{(i)} - \Phi(B^{\mC_{i}}) \right)$,
$i=1,\ldots,n$.\HS
\end{thm}

The theorem is proved by using 
the machinery of the $2+1$--dimensional TQFT of Reshetikhin and Turaev
induced by $\left( \mV,\{V_{i}\}_{i \in I}, \mD \right)$, 
see \cite[Sect.~5]{Hansen2}. The integer $\mu$,
present in the above theorem, is given by the following sum of Maslov indices
\begin{equation}\label{eq:e7}
\mu=\sum_{i=1}^{n} \mu((\partial Q_{i})_{*}(\lambda_{0}),(\partial Q_{i})_{*}(\lambda_{i}),N_{i}).
\end{equation}
We refer to \cite[Sect.~IV.3]{Turaev} or \cite[Sect.~5.2]{Hansen2}
for the definitions of Maslov index and Lagrangian subspace.
The spaces $\lambda_{0}$, $\lambda_{i}$, and $N_{i}$ are
Lagrangian subspaces of $H_{1}(\Sigma_{(1;)};\R)$ given by
$\lambda_{0}=\Span_{\R} \{ [ \alpha] \}$, $\lambda_{i}=\Span_{\R} \{ p_{i}[\alpha]+q_{i}[\beta] \}$,
and $N_{i}$ equal to the kernel of the inclusion homomorphism
$H_{1}(\partial U_{i};\R) \to H_{1}(M_{i-1} \sm \interior (U_{i});\R)$, where $M_{i}$
is the manifold obtained by doing surgery on
$M$ along $\left( \amalg_{j=1}^{i} L_{j},\amalg_{j=1}^{i} Q_{j} \right)$
with surgery coefficients $\{ p_{j}/q_{j} \}_{j=1}^{i}$, $i=1,2,\ldots,n$, and
$M_{0}=M$.

We have
the following corollary to  
\refthm{surgery-formula3sphere} and \ref{surgery-formula}:

\begin{cor}\label{surgery-formula2}
Let the situation be as in {\em \refthm{surgery-formula}} and let $B_{i} \in \SL (2,\Z)$
with first column equal to $\pm \left( \begin{array}{c}
		p_{i} \\
		q_{i}
		\end{array}
\right)$, $i=1,2,\ldots,n$. Then
\begin{eqnarray*}
\tau_{r}^{\frg}(M',\Omega) &=& \left( \exp\left( \frac{\pi\I}{\ch} |\rho|^{2} \right) \exp\left( -\frac{\pi\I}{\kappa} |\rho|^{2} \right)\right)^{\sum_{i=1}^{n} \Phi(B_{i}) -3\mu } \\
 & & \hspace{.4in} \times \sum_{\lambda \in \col(L)} \tau_{r}^{\frg}(M,\Gamma(L,\lambda) \cup \Omega) \left( \prod_{i=1}^{n} (\tilde{B}_{i})_{\lambda(L_{i})\rho} \right),
\end{eqnarray*}
where $\mu$ is given by {\em (\ref{eq:e7})}. If all the $q_{i}$ are different from
$0$, a similar formula holds with $M'$, $M$, $L$ and $\mu$ replaced by $M$, $S^{3}$, $L_{0}$, and
$\sigma$ respectively, where $M$, $L_{0}$ and $\sigma$ are as in
{\em \refthm{surgery-formula3sphere}}.\HS
\end{cor}

To see this simply choose tuples of integers $\mC_{i}$ such that
$B_{i}=B^{\mC_{i}}$ and use (\ref{eq:anomaly}), (\ref{eq:omega}), and
(\ref{eq:sumformula}) and the fact that $\mC_{i}$
is a continued fraction expansion of $p_{i}/q_{i}$, $i=1,2,\ldots,n$.

\ppar
\small
\sh{Appendix}

In the proof of the main \reflem{lem:main}, the reciprocity formula, \refprop{prop:gauss},
played a crucial role. The idea of using the reciprocity formula stems from Jeffrey's
calculations in \cite{Jeffrey1}, \cite{Jeffrey2} as already stated in the introduction.
Jeffrey's proof of the reciprocity formula, which is a verbatim generalization
of the argument presented in the proof of \cite[Chap.~IX Theorem 1]{Chandrasekharan} can be found
in her thesis \cite{Jeffrey1}.

In this appendix we will first sketch Jeffrey's proof, which builds on a limiting case of the
Poisson summation formula applied to Gaussian functions. A main ingredient is the
Fourier transformation of a Gaussian function, which rely on an analytic
continuation argument for complex functions of several variables.
Next we will present a slightly different
argument which avoids the direct use of this Fourier transform result.

\rk{Jeffrey's proof of \refprop{prop:gauss}}  According to the Poisson summation formula we have
$$
\sum_{n \in \Z^{l}} \hat{\phi}(an) = (2\pi/a)^{l} \sum_{n \in \Z^{l}} \phi(2\pi n/a), \hspace{.2in} a>0
$$
for $\phi \in \mS (\R^{l})$, see \cite[p.~178]{Hormander}. Here $\mS (\R^{l})$,
is the Schwartz space of smooth ($C^{\infty}$) functions that are rapidly
decreasing at infinity, and $\hat{\phi}$ is the Fourier transformation of $\phi$.
Recall that
$$
\hat{\phi}(\xi) = \int_{\R^{l}} e^{-i \la x,\xi\ra} \phi(x)dx, \hspace{.2in} \phi \in L^{1}(\R^{l})
$$
for $\xi \in \R^{l}$, where $\la\cdot ,\cdot \ra$ is the standard inner product in $\R^{l}$.
(Below $\la\cdot,\cdot\ra$ will also be used to denote the inner product in $V$, but
the meaning of $\la\cdot,\cdot\ra$ will always be clear from the context.)
We recall that $\mS(\R^{l}) \subseteq L^{1}(\R^{l})$.

Let $V_{\C}=V \otimes_{\R} \C$ and extend $\la \cdot,\cdot \ra$ to a hermitian product in
$V_{\C}$ in the usual way. Let $\tau \in \End_{\C}(V_{\C}))$ such that
the imaginary part of $\tau$ is positive definite. By assumption the sum
$$
\sum_{\lambda \in \Lambda} \exp(\pi\I\la\tau(\lambda),\lambda\ra)\exp(2\pi\I\la\lambda,\psi\ra)
$$
is absolutely convergent. In general, if $\Omega$ is a
symmetric complex $l\times l$--matrix with positive definite imaginary part and $z \in \C^{l}$
and if $\phi \co \R^{l} \to \C$ is given by
$$
\phi(x) =\exp(\pi\I x^{t} \Omega x + 2\pi \I x^{t}z)
$$
then
$$
\hat{\phi}(2\pi\xi) =\left( \det \left( \frac{\Omega}{\I} \right)\right)^{-1/2} \exp \left( -\pi\I (\xi-z)^{t} \Omega^{-1} (\xi-z) \right).
$$
Here the square root is positive on the positive real axis with a cut along the
negative real axis. To see this, first assume that $\Omega^{-1}z \in \R^{l}$
and complete the square and use \cite[Theorem 7.6.1]{Hormander}.
The general case then follows by analytic continuation.

By using this result together with the Poisson summation formula we get
\begin{eqnarray}\label{eq:square}
&&\Sigma_{\lambda \in \Lambda} \exp(\pi\I\la\tau(\lambda),\lambda\ra)\exp(2\pi\I\la\lambda,\psi\ra) \\
 &&\hspace{.3in} = \vol(\Lambda)^{-1}\left(\det\left( \frac{\tau}{\I} \right)\right)^{-1/2} \sum_{\mu \in \Lambda^{*}} \exp\left(-\pi\I\la \tau^{-1}(\mu+\psi),\mu+\psi\ra\right) \nonumber
\end{eqnarray}
if we futhermore e.g.\ assume that $\tau$ can be represented by a symmetric
matrix w.r.t.\ a basis of $V_{\C}$ of the form $\{ w_{1} \otimes 1,\ldots,w_{l}\otimes 1 \}$,
where $\{ w_{1},\ldots,w_{l} \}$ is some basis of $V$. This is always the
case in the situations where we will use (\ref{eq:square}) below, since $f \co V \to V$
is self-adjoint.
Let us use the identity (\ref{eq:square}) with $\tau = \frac{1}{r} f_{\C} +\I\vep\id_{V_{\C}}$, $\vep>0$,
where $f_{\C}=f \otimes\id_{\C}$. For $\lambda,\alpha \in \Lambda$ we have
$F(\lambda+r\alpha)=F(\lambda)$ by (\ref{eq:assumption}), where
\begin{equation}\label{eq:F}
F(\lambda)=\exp\left(\frac{\pi\I}{r}\la f(\lambda),\lambda\ra\right)\exp(2\pi\I\la\lambda,\psi\ra),
\end{equation}
so the left-hand side of (\ref{eq:square}) becomes
$$
\LHS (\vep ) = \sum_{\Lambda/ r\Lambda} F(\lambda) \sum_{\alpha \in \Lambda} \exp (-\pi\vep|\lambda+r\alpha|^{2}).
$$
Here
\begin{eqnarray*}
\sum_{\alpha \in \Lambda} \exp (-\pi\vep|\lambda+r\alpha|^{2}) &=& \vol(\Lambda^{*}) \left( \frac{1}{\vep r^{2}} \right)^{l/2} \\
 && \hspace{.1in} \times \sum_{\beta \in \Lambda^{*}} \exp \left( -\pi \la \beta, \frac{1}{\vep r^{2}} \beta \ra \right)
             \exp\left( 2\pi \I \la \beta,\frac{\lambda}{r}\ra \right)
\end{eqnarray*}
by using (\ref{eq:square}) with the roles of $\Lambda$ and $\Lambda^{*}$ reversed and with
$\tau^{-1}=-\I\vep r^{2} \id_{V_{\C}}$ (so
$\tau=\I \frac{1}{\vep r^{2}} \id_{V_{\C}}$ has positive definite imaginary part).
We first observe that
$$
\lim_{\vep \ria 0_{+}} \vep^{l/2} \LHS (\vep) = r^{-l} \vol(\Lambda^{*}) \sum_{\Lambda/r\Lambda} F(\lambda).
$$
This follows by the fact that
$$
\lim_{\vep \ria 0_{+}} \sum_{\beta \in \Lambda^{*}} \exp \left( -\pi\la\beta,\frac{1}{\vep r^{2}}\beta\ra \right)
        \exp\left( 2\pi\I\la\beta,\frac{\lambda}{r}\ra \right) = 1,
$$
which follows by the fact that
\begin{eqnarray*}
&& \left| \sum_{\beta\in\Lambda^{*} \sm \{ 0\} } \exp \left(-\pi\la\beta,\frac{1}{\vep r^{2}} \beta \ra \right) 
  \exp\left( 2\pi\I\la\beta,\frac{\lambda}{r}\ra\right)\right| \\
&& \hspace{.3in} \leq \sum_{\beta\in\Lambda^{*} \sm \{ 0\} } \exp \left(-\frac{\pi}{\vep r^{2}} |\beta|^{2}\right) 
     \leq \exp\left(-\frac{\pi}{2\vep r^{2}} c^{2}\right) \sum_{\beta\in\Lambda^{*}}\exp\left(-\frac{\pi}{2r^{2}}|\beta|^{2}\right)
\end{eqnarray*}
for $\vep \in ]0,1]$, $c=\min\{\; |\beta|\; : \; \beta \in \Lambda^{*} \sm \{0\} \; \} >0$.

Hereafter we calculate $\lim_{\vep \ria 0_{+}} \vep^{l/2}\RHS (\vep)$, where $\RHS (\vep)$ is
the right-hand side of (\ref{eq:square}) with $\tau=\frac{1}{r}f_{\C} + \I\vep\id_{V_{\C}}$.
Note that
$$
\tau^{-1} = rf_{\C}^{-1} - \I \vep (rf_{\C}^{-1})^{2}\left(\id_{V_{\C}} + \I\vep r f_{\C}^{-1}\right)^{-1}.
$$
By using this we find that $\RHS (\vep)$ is equal to
\begin{eqnarray*}
 && \vol(\Lambda)^{-1} \left(\det\left(\frac{\tau}{\I}\right)\right)^{-1/2}
   \sum_{\mu\in \Lambda^{*} / f(\Lambda^{*})} \exp( -\pi\I\la \mu+\psi,rf^{-1}(\mu+\psi)\ra) \\
 &&\hspace{.2in} \times \sum_{\beta \in \Lambda^{*}} \exp\left( -\pi\vep r^{2}\la (\id + \I\vep r f_{\C}^{-1})^{-1} 
         \left( f^{-1}(\mu+\psi) + \beta\right),f^{-1}(\mu+\psi) + \beta\ra \right),
\end{eqnarray*}
where we use that $G(\mu + f(\beta))=G(\mu)$ for $\mu,\beta \in \Lambda^{*}$, where
$G(\mu)=\exp(-\pi\I \la \mu+\psi,rf^{-1}(\mu+\psi)\ra)$. The sum $\sum_{\beta \in \Lambda^{*}}$
can now be calculated by using (\ref{eq:square}) once more with $\psi$ replaced by
$f^{-1}(\mu+\psi)$ and $\tau=\frac{\I}{\vep r^{2}}(\id_{V_{\C}}+\I\vep r f_{\C})$.
By doing this we get
\begin{eqnarray*}
\RHS (\vep) &=& \left(\det\left(\frac{1}{r\I}f_{\C} + \vep \id_{V_{\C}} \right)\right)^{-1/2} \left(\frac{1}{\vep r^{2}}\right)^{l/2}
    \left(\det (\id + \I\vep r f_{\C} )\right)^{1/2} \\
 && \hspace{.2in} \times \sum_{\mu \in \Lambda^{*} /f(\Lambda^{*})} \exp\left( -\pi\I\la\mu+\psi,rf^{-1}(\mu+\psi)\ra\right) \\
 &&\hspace{.3in}\times \sum_{\lambda \in \Lambda} \exp\left(-\frac{\pi}{\vep r^{2}} |\lambda|^{2} \right)
    \exp\left(-\frac{\pi\I}{r}\la f(\lambda),\lambda\ra\right) \\
 && \hspace{1.0in} \times \exp\left(2\pi\I\la f^{-1}(\mu+\psi),\lambda\ra\right).
\end{eqnarray*}
As before we find that the sum over $\Lambda$ converges to $1$ as $\vep \ria 0_{+}$. Therefore
\begin{eqnarray*}
\lim_{\vep \ria 0_{+}} \vep^{l/2} \RHS (\vep) &=& \left( \det\left(\frac{f}{\I}\right)\right)^{-1/2} r^{-l/2}\\
 && \hspace{.1in} \times
 \sum_{\mu \in \Lambda^{*} /f(\Lambda^{*})} \exp\left( -\pi\I\la\mu+\psi,rf^{-1}(\mu+\psi)\ra\right).
\end{eqnarray*}
Since $\lim_{\vep \ria 0_{+}} \vep^{l/2}\RHS (\vep)=\lim_{\vep \ria 0_{+}} \vep^{l/2} \LHS (\vep)$ the result follows.

\rk{A second proof of \refprop{prop:gauss}} This proof builds on the following periodicity result:

\begin{lem}\label{lem:latticeperiodicity}
Let $\Lambda$ be a lattice in $V$ and let
$h \co V \to V$ be a linear map such that $h(\Lambda) \subseteq \Lambda$.
Moreover, let $g_{\vep} \co V \to V$, $\vep \in ]0,a]$, be a curve of self-adjoint
positive definite maps such that $g_{\vep} \ria g_{0}$ in $\End_{\R}(V)$ as
$\vep \ria 0$, where $g_{0}$ is positive definite, $a$ being a fixed positive
number. Finally, let $v_{0} \in V$
be fixed but arbitrary and let $F \co V \to \C$ be a map such that  
$$
F(\lambda + h(\alpha))=F(\lambda)
$$
for all $\lambda,\alpha \in \Lambda$. Then
\begin{eqnarray*}
\sum_{\lambda \in \Lambda/h(\Lambda)} F(\lambda) &=& \vol(\Lambda) |\det(h)|\sqrt{\det(g_{0})} \\
 &&\hspace{.2in}\times \lim_{\vep \ria 0_{+}} \vep^{l/2} \sum_{\lambda \in \Lambda} e^{-\pi\vep \la \lambda + v_{0}, g_{\vep}(\lambda + v_{0}) \ra}
     F(\lambda).
\end{eqnarray*}
\end{lem}

\begin{proof}
By assumption we have
\begin{eqnarray*}
&&\sum_{\lambda \in \Lambda} e^{-\pi\vep \la \lambda + v_{0}, g_{\vep}(\lambda + v_{0}) \ra} F(\lambda) \\
 &=&\hspace{.2in} \sum_{\lambda \in \Lambda/h(\Lambda)} F(\lambda)
         \sum_{\alpha \in \Lambda} e^{-\pi\vep \la \lambda + v_{0} + h(\alpha), g_{\vep}(\lambda + v_{0}+h(\alpha)) \ra}.
\end{eqnarray*}
The lemma will therefore follow if we can show that
$$
\vol(\Lambda) |\det(h)|\sqrt{\det(g_{0})}\lim_{\vep \ria 0_{+}} \vep^{l/2}
   \sum_{\alpha \in \Lambda} e^{-\pi\vep \la v + h(\alpha), g_{\vep}(v+h(\alpha)) \ra}=1
$$
for any $v \in V$. This is done by changing the sum $\sum_{\alpha \in \Lambda}$
to a sum over $\Z^{l}$ (using coordinates) and then use the Poisson summation formula to this sum.
Note that if $\mV=\{v_{1},\ldots,v_{l}\}$ is a basis for $V$ such that $\Lambda$ is generated by
this set over the integers and if $\mW=\{w_{1},\ldots,w_{l}\}$ is an orthonormal basis for $V$
then $\vol(\Lambda)=\det(k)$, where $k\co V \to V$ is the linear isomorphism given by
$k(w_{j})=v_{j}$, $j=1,2,\ldots,l$.
\end{proof}

Now let $F \co \Lambda \to \C$ be given by (\ref{eq:F}),
and let $h=r\id_{V}\co V \to V$. Then
$F(\lambda + h(\alpha))=F(\lambda)$ for $\alpha,\lambda\in\Lambda$ by
(\ref{eq:assumption}), so by the above
lemma we get
$$
\sum_{\lambda \in \Lambda/r\Lambda} F(\lambda) = \vol(\Lambda)r^{l}
 \lim_{\vep \ria 0_{+}} \vep^{l/2} \sum_{\lambda \in \Lambda} e^{-\pi\vep |\lambda|^{2}} F(\lambda).
$$
To continue we use coordinates. Let $\mV$ and $\mW$ be bases for $V$ as in the proof of \reflem{lem:latticeperiodicity},
and let $C$ be the matrix of $f$ w.r.t.\ $\mW$. Moreover, let $D=(d_{ij})_{i,j=1}^{l}$ such that
$v_{j}=\sum_{i=1}^{l}d_{ij}w_{i}$. Then
\begin{eqnarray*}
&&\sum_{\lambda \in \Lambda} e^{-\pi\vep |\lambda|^{2}} F(\lambda) \\
 &&\hspace{.2in} = \sum_{n \in \Z^{l}} e^{-\pi\vep \la Dn,Dn\ra}
          \exp\left(\frac{\pi\I}{r}\la Dn,CDn\ra\right)\exp\left(2\pi\I\la Dn,y\ra\right),
\end{eqnarray*}
where $y$ are the coordinates of $\psi$ w.r.t.\ the basis $\mW$. By the Poisson summation formula
we get
\begin{eqnarray*}
&&\sum_{\lambda \in \Lambda} e^{-\pi\vep |\lambda|^{2}} F(\lambda) = |\det(D)|^{-1} 
        \sum_{m\in \Z^{l}} \int_{\R^{l}} e^{2\pi\I\la m,D^{-1}x\ra} \\
 &&\hspace{.4in} \times e^{-\pi\vep \la x,x\ra }\exp\left(\frac{\pi\I}{r}\la x,Cx\ra\right)\exp\left(2\pi\I\la x,y\ra\right)dx.
\end{eqnarray*}
The summands are all Gaussian integrals and can be calculated by diagonalizing $C$.
In fact, if we choose an orthogonal matrix $Q$ such that 
$Q^{-1}CQ=\Diag(\lambda_{1},\ldots,\lambda_{l})$ and let $\eta=\vep r^{2}$, 
then we arrive at the following identity
\begin{eqnarray*}
&&\sum_{\lambda \in \Lambda/r\Lambda} F(\lambda) = r^{l/2} \left( \det \left(\frac{f}{\I}\right) \right)^{-1/2} 
     \lim_{\eta \ria 0_{+}} \eta^{l/2} \\
 &&\hspace{.3in}\times \sum_{m\in \Z^{l}} 
          \exp\left(-\pi\eta \la m + D^{t}y,D^{-1}C^{-1}QH(\eta)Q^{-1}C^{-1}(D^{-1})^{t}(m+D^{t}y) \ra \right)\\
 &&\hspace{.3in}\times\exp\left(-\pi\I r \la D^{-1}QH(\eta)Q^{-1}C^{-1}(D^{-1})^{t}(m+D^{t}y),m+D^{t}y\ra\right),
\end{eqnarray*}
where $H(\eta)=\Diag(f_{1}(\eta),\ldots,f_{l}(\eta))$, where 
$f_{j}(\eta)=\left( 1 + \left(\frac{\eta}{r\lambda_{j}}\right)^{2}\right)^{-1} \ria 1$ as $\eta \ria 0_{+}$.
Next we use the following technical but straightforward

\begin{lem}
Let $a>0$ and let $A\co ]0,a] \to \GL (l,\R)$ be a curve of positive
definite symmetric matrices such that $A(\vep) \ria A_{0}$ as $\vep \ria 0_{+}$,
where $A_{0}$ is a symmetric positive definite matrix. Moreover, let
$B\co ]0,a] \to \GL (l,\R)$ be a curve and $B_{0}$ a fixed matrix such that
$$
\lim_{\vep \ria 0_{+}} \frac{1}{\vep}\left(B(\vep)-B_{0}\right) =0.
$$
Then we have
\begin{eqnarray*}
&&\lim_{\vep\ria 0_{+}} \vep^{l/2} \sum_{n \in \Z^{l}} e^{-\vep \la n+x_{1},A(\vep)(n+x_{1})\ra }
     \exp\left(\I\la n+x_{2},B(\vep)(n+x_{2})\ra \right) \\
&&\hspace{.1in} = \lim_{\vep\ria 0_{+}} \vep^{l/2} \sum_{n \in \Z^{l}} e^{-\vep \la n+x_{1},A(\vep)(n+x_{1})\ra }
     \exp\left(\I\la n+x_{2},B_{0}(n+x_{2})\ra \right)
\end{eqnarray*}
for any $x_{1},x_{2} \in \R^{l}$ in the sense that if one of the two limits exists then does
the other and they are equal.\HS
\end{lem}

By this lemma (and \reflem{lem:latticeperiodicity}) we get that
\begin{eqnarray*}
\sum_{\lambda \in \Lambda/r\Lambda} F(\lambda) &=& r^{l/2} \left( \det \left(\frac{f}{\I}\right) \right)^{-1/2} 
     \lim_{\eta \ria 0_{+}} \eta^{l/2} \\
 &&\hspace{.1in}\times \sum_{m\in \Z^{l}} 
          \exp\left(-\pi\eta \la m + D^{t}y,D^{-1}C^{-2}(D^{-1})^{t}(m+D^{t}y) \ra \right)\\
 &&\hspace{.1in}\times\exp\left(-\pi\I r \la D^{-1}C^{-1}(D^{-1})^{t}(m+D^{t}y),m+D^{t}y\ra\right) \\
 &=& r^{l/2} \left( \det \left(\frac{f}{\I}\right) \right)^{-1/2} 
     \lim_{\eta \ria 0_{+}} \eta^{l/2} \\
 && \hspace{.2in} \times \sum_{\mu \in \Lambda^{*}} \exp\left(-\pi\eta\la \mu + \psi, f^{-2}(\mu+\psi)\ra\right) \\
 && \hspace{.5in} \times \exp\left(-\pi\I r \la \mu +\psi,f^{-1}(\mu + \psi) \ra \right).
\end{eqnarray*}
Now $f(\Lambda^{*}) \subseteq \Lambda^{*}$ by assumption and $G(\mu+f(\beta))=G(\mu)$
for $\mu,\beta \in \Lambda^{*}$, where
$G(\mu)=\exp\left(-\pi\I r \la \mu +\psi,f^{-1}(\mu + \psi) \ra \right)$, so
by \reflem{lem:latticeperiodicity} we get
\refprop{prop:gauss}.
\end{document}